\documentclass[a4paper]{article}

\usepackage[cp1250]{inputenc}
\usepackage[IL2]{fontenc}

\usepackage{a4wide}
\usepackage[english]{babel}
\usepackage{euscript}
\usepackage{amstext,amsbsy,amscd,amssymb}
\usepackage{amsmath,enumitem}
\usepackage{amsfonts}
\usepackage{graphics}
\usepackage{graphicx}
\usepackage{microtype}
\usepackage{xcolor}
\usepackage{dynkin-diagrams}
\usepackage{tikz,caption,subcaption}
\usetikzlibrary{arrows,decorations.pathmorphing,backgrounds,positioning,fit,petri}
\usepackage{afterpage}

\include{diagxy}

\allowdisplaybreaks

\let\rarr=\rightarrow

\let\veps=\varepsilon
\let\mcal=\mathcal
\let\mfrak=\mathfrak

\let\bra=\langle
\let\ket=\rangle

\def\N{\mathbb{N}}
\def\Z{\mathbb{Z}}
\def\R{\mathbb{R}}
\def\C{\mathbb{C}}

\def\Q{\mathbb{Q}}

\DeclareMathSymbol{\squares}{\mathord}{AMSa}{"03}

\DeclareMathOperator{\Ann}{Ann}

\DeclareMathOperator{\diag}{diag}

\DeclareMathOperator{\Specm}{Specm}
\DeclareMathOperator{\End}{End}
\DeclareMathOperator{\Hom}{Hom}
\DeclareMathOperator{\Res}{Res}
\DeclareMathOperator{\Aut}{Aut}

\DeclareMathOperator{\rank}{rank}
\DeclareMathOperator{\ad}{ad}

\def\ad{\mathop {\rm ad} \nolimits}
\def\gr{\mathop {\rm gr} \nolimits}

\long\def\proof #1{\noindent \emph{Proof.}\ #1 \hfill $\squares$

\medskip}

\newcounter{num}[section]
\numberwithin{equation}{section}
\numberwithin{num}{section}

\long\def\definition #1 {\refstepcounter{num} \noindent {\bf Definition \thenum.} #1

\medskip}

\long\def\theorem #1{\refstepcounter{num} \noindent \ifnum\value{section}>0 {\bf Theorem \thenum.} #1 \else {\bf Theorem \Alph{num}.} #1 \fi

\medskip}

\long\def\lemma #1{\refstepcounter{num}  \noindent {\bf Lemma \thenum.} #1

\medskip}

\long\def\proposition #1{\refstepcounter{num}  \noindent {\bf Proposition \thenum.} #1

\medskip}

\long\def\corollary #1{\refstepcounter{num}  \noindent {\bf Corollary \thenum.} #1

\medskip}

\long\def\example #1{\noindent {\bf Example.}\ #1}

\long\def\remark #1{\noindent {\bf Remark.}\ #1}

\makeatletter
\newcommand{\raisemath}[1]{\mathpalette{\raisem@th{#1}}}
\newcommand{\raisem@th}[3]{\raisebox{#1}{$#2#3$}}
\makeatother

\makeatletter
\newcommand*\if@single[3]{%
  \setbox0\hbox{${\mathaccent"0362{#1}}^H$}%
  \setbox2\hbox{${\mathaccent"0362{\kern0pt#1}}^H$}%
  \ifdim\ht0=\ht2 #3\else #2\fi
  }
\newcommand*\rel@kern[1]{\kern#1\dimexpr\macc@kerna}
\newcommand*\widebar[1]{\@ifnextchar^{{\wide@bar{#1}{0}}}{\wide@bar{#1}{1}}}
\newcommand*\wide@bar[2]{\if@single{#1}{\wide@bar@{#1}{#2}{1}}{\wide@bar@{#1}{#2}{2}}}
\newcommand*\wide@bar@[3]{%
  \begingroup
  \def\mathaccent##1##2{%
    \if#32 \let\macc@nucleus\first@char \fi
    \setbox\z@\hbox{$\macc@style{\macc@nucleus}_{}$}%
    \setbox\tw@\hbox{$\macc@style{\macc@nucleus}{}_{}$}%
    \dimen@\wd\tw@
    \advance\dimen@-\wd\z@
    \divide\dimen@ 3
    \@tempdima\wd\tw@
    \advance\@tempdima-\scriptspace
    \divide\@tempdima 10
    \advance\dimen@-\@tempdima
    \ifdim\dimen@>\z@ \dimen@0pt\fi
    \rel@kern{0.6}\kern-\dimen@
    \if#31
      \overline{\rel@kern{-0.6}\kern\dimen@\macc@nucleus\rel@kern{0.4}\kern\dimen@}%
      \advance\dimen@0.4\dimexpr\macc@kerna
      \let\final@kern#2%
      \ifdim\dimen@<\z@ \let\final@kern1\fi
      \if\final@kern1 \kern-\dimen@\fi
    \else
      \overline{\rel@kern{-0.6}\kern\dimen@#1}%
    \fi
  }%
  \macc@depth\@ne
  \let\math@bgroup\@empty \let\math@egroup\macc@set@skewchar
  \mathsurround\z@ \frozen@everymath{\mathgroup\macc@group\relax}%
  \macc@set@skewchar\relax
  \let\mathaccentV\macc@nested@a
  \if#31
    \macc@nested@a\relax111{#1}%
  \else
    \def\gobble@till@marker##1\endmarker{}%
    \futurelet\first@char\gobble@till@marker#1\endmarker
    \ifcat\noexpand\first@char A\else
      \def\first@char{}%
    \fi
    \macc@nested@a\relax111{\first@char}%
  \fi
  \endgroup
}
\makeatother

\newcommand\rsmraise[1]{%
  \ifx#1\displaystyle .8\else
    \ifx#1\textstyle .8\else
      \ifx#1\scriptstyle .6\else
        .45%
      \fi
    \fi
  \fi}

\title{Admissible representations of simple affine vertex algebras}

\author{Vyacheslav Futorny, Oscar Armando Hern\'andez Morales, Libor K\v{r}i\v{z}ka}

\AtEndDocument{\bigskip{\footnotesize%
  (V.\,Futorny) \textsc{Instituto de Matem\'atica e Estat\'istica, Universidade de S\~ao Paulo,  S\~ao Paulo, Brazil, and  International Center for Mathematics, SUSTech, Shenzhen, China} \par
  \textit{E-mail address}: \texttt{vfutorny@gmail.com} \par
  \addvspace{\medskipamount}
  (O.\,Hern\'andez) \textsc{Instituto de Matem\'atica e Estat\'istica, Universidade de S\~ao Paulo,  S\~ao Paulo, Brazil} \par
  \textit{E-mail address}: \texttt{oscarhm@ime.usp.br} \par
  \addvspace{\medskipamount}
  (L.\,K\v{r}i\v{z}ka) \textsc{Instituto de Matem\'atica e Estat\'stica, Universidade de Sa\~{o} Paulo, Sa\~{o} Paulo,  Brasil} \par
  \textit{E-mail address}: \texttt{krizka.libor@gmail.com} \par
}}

\date{}

\begin{document}

\maketitle

\begin{abstract}
We provide an explicit combinatorial description of highest weights of simple highest weight modules over the simple affine vertex algebra $\mathcal{L}_\kappa(\mathfrak{sl}_{n+1})$ with $n \in \N$ of admissible level $\kappa$. For admissible simple highest weight modules corresponding to the principal, subregular and maximal parabolic nilpotent orbits we give a realization using the Gelfand--Tsetlin theory, which also allows us to obtain a realization of certain classes of simple admissible $\mfrak{sl}_2$-induced modules in these orbits. In particular, simple admissible $\mfrak{sl}_2$-induced modules are fully classified for the principal nilpotent orbit.

\medskip

\medskip


\end{abstract}

\thispagestyle{empty}

\tableofcontents

%%%%%%%%%%%%%%%%%%%%%%%%%%%%%%%%%%%%%%%%%%%%%%%%%%%%%%%%%%%%%%%%%%%%%%%%%%%%%%%%%%%%%%%%%%%
%%%%%%%%%%%%%%%%%%%%%%%%%%%%%%%%%%%%%%%%%%%%%%%%%%%%%%%%%%%%%%%%%%%%%%%%%%%%%%%%%%%%%%%%%%%

\section*{Introduction}
\addcontentsline{toc}{section}{Introduction}

Representation theory of simple affine vertex algebras is actively expanding field with important applications in 2d conformal field theory. Theory of highest weight representations was developed in \cite{Kac-Wakimoto1989} and \cite{Arakawa2016}. More general positive energy representations were studied in \cite{Adamovic-Milas1995}, \cite{Adamovic-Milas2007}, \cite{Adamovic-Perse2008}, \cite{Arakawa-Futorny-Ramirez2017}, \cite{Kawasetsu-Ridout2019}, \cite{Kawasetsu-Ridout2019b}, \cite{Futorny-Krizka2021}, \cite{Futorny-Morales-Ramirez2020}, \cite{Adamovic-Milas2009}, \cite{Adamovic2016}, \cite{Auger-Creutzig-Ridout2018}. By a key result of \cite{Zhu1996}, there is a one-to-one correspondence between simple positive energy representations of the simple affine vertex algebra $\mcal{L}_\kappa(\mfrak{g})$ of level $\kappa$ attached to a simple Lie algebra $\mfrak{g}$ and simple modules over the corresponding Zhu's algebra $A(\mcal{L}_\kappa(\mfrak{g}))$ which is a subalgebra of the universal enveloping algebra $U(\mfrak{g})$ of $\mfrak{g}$. Moreover, the Zhu's algebra $A(\mcal{L}_\kappa(\mfrak{g}))$ is isomorphic to $U(\mfrak{g})/I_\kappa(\mfrak{g})$, where $I_\kappa(\mfrak{g})$ is a two-sided ideal of $U(\mfrak{g})$. Hence, a full classification of simple positive energy representations of $\mcal{L}_\kappa(\mfrak{g})$ is reduced to a classification of simple $\mfrak{g}$-modules such that their annihilators in $U(\mfrak{g})$ contain $I_\kappa(\mfrak{g})$. On the other hand, from \cite{Duflo1977} for any primitive ideal $I$ of $U(\mfrak{g})$ there exists $\lambda \in \mfrak{h}^*$ such that $I = \Ann_{U(\mfrak{g})}\!L^\mfrak{g}_\mfrak{b}(\lambda)$, where $\mfrak{h}$ is a Cartan subalgebra of $\mfrak{g}$ and $L^\mfrak{g}_\mfrak{b}(\lambda)$ is the simple $\mfrak{g}$-module with highest weight $\lambda$ for a Borel subalgebra $\mfrak{b}$ of $\mfrak{g}$ containing $\mfrak{h}$. A complete description of weights $\lambda \in \mfrak{h}^*$ satisfying $I_\kappa(\mfrak{g}) \subset \Ann_{U(\mfrak{g})}\!L^\mfrak{g}_\mfrak{b}(\lambda)$ for an admissible level $\kappa$ was done in \cite{Arakawa2016}. We call such weights admissible of level $\kappa$. Besides, we say that a $\mfrak{g}$-module $E$ is admissible of level $\kappa$ if $\Ann_{U(\mfrak{g})}\!E = \Ann_{U(\mfrak{g})}\!L^\mfrak{g}_\mfrak{b}(\lambda)$ for some admissible weight $\lambda \in \mfrak{h}^*$.

In particular, explicit Gelfand--Tsetlin tableau realizations of certain highest weight modules over the simple Lie algebra $\mfrak{g}$ of type $A$ were constructed in \cite{Arakawa-Futorny-Ramirez2017}. All such simple admissible highest weight $\mfrak{g}$-modules of admissible level $\kappa$ corresponding to the minimal nilpotent orbit were described in \cite{Futorny-Morales-Ramirez2020}. A class of non-highest weight $\mfrak{g}$-modules was also constructed in \cite{Futorny-Krizka2019} and \cite{Futorny-Krizka2019b} by using the Arkhipov's twisting functor. A characterization of which of these $\mfrak{g}$-modules are admissible is given in \cite{Futorny-Krizka2021}.
\medskip

The goal of this paper is threefold. First, we obtain a convenient combinatorial description of admissible weights of admissible level $\kappa$ for a simple Lie algebra $\mfrak{g}$ of type $A$ in Theorem \ref{thm:Pr_k parametrization}. In particular, weights corresponding to the minimal, subregular and principal nilpotent orbits of $\mfrak{g}$ are described. Next, we give an explicit construction of all simple admissible highest weight $\mfrak{g}$-modules in Theorem \ref{thm-h.w.-adm}. Finally, we address the problem of an explicit realization of simple admissible  weight $\mfrak{g}$-modules in the principal, subregular and maximal parabolic nilpotent orbits. We give a Gelfand--Tsetlin realization of simple admissible highest weight $\mfrak{g}$-modules and certain classes of simple admissible $\mfrak{sl}_2$-induced $\mfrak{g}$-modules in these orbits. Constructed modules are \emph{tame}, i.e.\ they admit a diagonalizable action of a certain Gelfand--Tsetlin subalgebra $\Gamma$ of $U(\mfrak{g})$. In addition, they are \emph{strongly tame}, i.e.\ they have a basis consisting of Gelfand--Tsetlin tableaux with the explicit action of $\mfrak{g}$ given by the classical Gelfand--Tsetlin formulas with respect to $\Gamma$. Obtained results provide an evidence for the following conjecture.
\medskip

\noindent
{\bf Conjecture 1.} Let $L^\mfrak{g}_\mfrak{b}(\lambda)$ be a simple admissible highest weight $\mfrak{g}$-module of level $\kappa$ with highest weight $\lambda \in \mfrak{h}^*$. If $\langle \lambda,\alpha^\vee \rangle \notin \Z$ for a simple root $\alpha$, then there exists $w \in \widebar{W}$ such that $\smash{D^\nu_{f_\alpha}\!(L^\mfrak{g}_\mfrak{b}(\lambda))}$ is a strongly tame $\Gamma_{\rm st}(w(\Pi))$-Gelfand--Tsetlin $\mfrak{g}$-module, where $\widebar{W}$ is the subgroup of $\Aut(\mfrak{g})$ generated by the Tits extension of the Weyl group $W$ of $\mfrak{g}$ and by the symmetries of Dynkin diagram of $\mfrak{g}$.
\medskip

We show that Conjecture 1 holds for any $\lambda \in \mfrak{h}^*$ which corresponds to the principal nilpotent orbit (Theorem \ref{prop-diff-hw}), and for the representatives of the subregular nilpotent orbit (Theorem \ref{thm-subreg}) and maximal parabolic nilpotent orbits (Theorem \ref{Pro:maxparabolic}). Moreover, in the case of maximal parabolic and subregular nilpotent orbits Conjecture 1 also holds for suitable non-simple roots. Earlier, this was shown for the minimal nilpotent orbit in \cite{Futorny-Morales-Ramirez2020}.

Further on, since any simple admissible cuspidal $\mfrak{g}$-module is obtained from a highest weight $\mfrak{g}$-module by a sequence of twisted localizations, we dare to make the following conjecture.
\medskip

\noindent
{\bf Conjecture 2.} Let $M$ be a simple admissible weight $\mfrak{g}$-module. Then there exists $w\in \widebar{W}$ such that $D_{f_{\alpha}}^\nu\!(M)$ is a strongly tame $\Gamma_{\rm st}(w(\Pi))$-Gelfand--Tsetlin $\mfrak{g}$-module.
\medskip

Explicit realizations of localized simple admissible $\mfrak{sl}_4$-modules is given in Section \ref{sl4-modules}.
\medskip

We denote by $\C$, $\R$, $\Z$, $\N_0$ and $\N$ the set of complex numbers, real numbers, integers, non-negative integers and positive integers, respectively. All algebras and modules are considered over the field of complex numbers.

%%%%%%%%%%%%%%%%%%%%%%%%%%%%%%%%%%%%%%%%%%%%%%%%%%%%%%%%%%%%%%%%%%%%%%%%%%%%%%%%%%%%%%%%%%
%%%%%%%%%%%%%%%%%%%%%%%%%%%%%%%%%%%%%%%%%%%%%%%%%%%%%%%%%%%%%%%%%%%%%%%%%%%%%%%%%%%%%%%%%%

\section{Admissible modules}

%%%%%%%%%%%%%%%%%%%%%%%%%%%%%%%%%%%%%%%%%%%%%%%%%%%%%%%%%%%%%%%%%%%%%%%%%%%%%%%%%%%%%%%%%%%

\subsection{Induced modules}

Let $\mfrak{g}$ be a complex semisimple finite-dimensional Lie algebra and let $\mfrak{h}$ be a Cartan subalgebra of $\mfrak{g}$. We denote by $\Delta$ the root system of $\mfrak{g}$ with respect to $\mfrak{h}$, by $\Delta_+$ a positive root system in $\Delta$ and by $\Pi \subset \Delta_+$ the set of simple roots. For $\alpha \in \Delta_+$, let $h_\alpha \in \mfrak{h}$ be the corresponding coroot and let $e_\alpha$ and $f_\alpha$ be basis of root subspaces $\mfrak{g}_\alpha$ and $\mfrak{g}_{-\alpha}$, respectively, defined by the requirement $[e_\alpha, f_\alpha] = h_\alpha$. We also set
\begin{align*}
  Q = \sum_{\alpha \in \Pi} \Z \alpha \qquad \text{and} \qquad P = \sum_{\alpha \in \Pi} \Z \omega_\alpha,
\end{align*}
where $\omega_\alpha \in \mfrak{h}^*$ for $\alpha \in \Pi$ is the fundamental weight determined by $\omega_\alpha(h_\gamma) = \delta_{\alpha,\gamma}$ for all $\gamma \in \Pi$. We call $Q$ the root lattice and $P$ the weight lattice. Further, we define the Weyl vector $\rho \in \mfrak{h}^*$ by
\begin{align*}
 \rho =  {1 \over 2} \sum_{\alpha \in \Delta_+} \alpha.
\end{align*}
Let $\kappa_\mfrak{g}$ be the Cartan--Killing form on $\mfrak{g}$ and $(\cdot\,,\cdot)_\mfrak{g}$ the corresponding induced bilinear form on $\mfrak{g}^*$. Whenever $\alpha \in \mfrak{h}^*$ satisfies $(\alpha,\alpha)_\mfrak{g} \neq 0$, we define $s_\alpha \in GL(\mfrak{h}^*)$ by
\begin{align*}
  s_\alpha(\gamma) = \gamma - {2(\alpha,\gamma)_\mfrak{g} \over (\alpha,\alpha)_\mfrak{g}}\, \alpha
\end{align*}
for $\gamma \in \mfrak{h}^*$. The Weyl group  $W$ of $\mfrak{g}$ is the subgroup of
$GL(\mfrak{h}^*)$ generated by $s_\alpha$ for $\alpha \in \Pi$. The dot action of $W$ on $\mfrak{h}^*$ is defined by $w\cdot \lambda=w(\lambda+\rho)-\rho$ for $w\in W$ and $\lambda\in \mfrak{h}^*$.

If $\mfrak{g}$ is a simple Lie algebra, we denote by $\kappa_0$ the $\mfrak{g}$-invariant symmetric bilinear form on $\mfrak{g}$ normalized in such a way that $(\theta,\theta)=2$, where $\theta \in \Delta_+$ is the highest root of $\mfrak{g}$ and $(\cdot\,,\cdot)$ is the corresponding induced bilinear form on $\mfrak{g}^*$. It easily follows that
\begin{align*}
  \kappa_\mfrak{g} = 2h^\vee \kappa_0,
\end{align*}
where $h^\vee$ denotes the dual Coxeter number of $\mfrak{g}$. Besides, we denote by $\theta_s$ the highest short root of $\mfrak{g}$. Then we have $(\theta_s,\theta_s)= 2/r^\vee$, where $r^\vee$ is the lacing number of $\mfrak{g}$, i.e.\ the maximal number of edges in the Dynkin diagram of $\mfrak{g}$. We also define
\begin{align*}
  P^\vee = \bigoplus_{\alpha \in \Pi} \Z\omega_\alpha^\vee \qquad  \text{and} \qquad P^\vee_+ = \bigoplus_{\alpha \in \Pi} \N_0\omega_\alpha^\vee,
\end{align*}
where $\omega_\alpha^\vee \in \mfrak{h}^*$ for $\alpha \in \Pi$ is the fundamental coweight defined by $(\omega_\alpha^\vee, \gamma) = \delta_{\alpha,\gamma}$ for all $\gamma \in \Pi$. We call $P^\vee$ the coweight lattice.
\medskip

The standard Borel subalgebra $\mfrak{b}$ of $\mfrak{g}$ is defined through $\mfrak{b} = \mfrak{h} \oplus \mfrak{n}$ with the nilradical $\mfrak{n}$ and the opposite nilradical $\widebar{\mfrak{n}}$ given by
\begin{align*}
  \mfrak{n} = \bigoplus_{\alpha \in \Delta_+} \mfrak{g}_\alpha \qquad \text{and} \qquad \widebar{\mfrak{n}} = \bigoplus_{\alpha \in \Delta_+} \mfrak{g}_{-\alpha}.
\end{align*}
In addition, we have the corresponding triangular decomposition
\begin{align*}
  \mfrak{g} = \widebar{\mfrak{n}} \oplus \mfrak{h} \oplus \mfrak{n}
\end{align*}
of the Lie algebra $\mfrak{g}$.

For a subset $\Sigma$ of $\Pi$ we denote by $\Delta_\Sigma$ the root subsystem in $\mfrak{h}^*$ generated by $\Sigma$. The standard parabolic subalgebra $\mfrak{p}_\Sigma$ of $\mfrak{g}$ associated to $\Sigma$ is defined as $\mfrak{p}_\Sigma = \mfrak{l}_\Sigma \oplus \widebar{\mfrak{u}}_\Sigma$, where the nilradical $\mfrak{u}_\Sigma$, the opposite nilradical $\widebar{\mfrak{u}}_\Sigma$ and the Levi subalgebra $\mfrak{l}_\Sigma$ are given by
\begin{align*}
  \mfrak{u}_\Sigma = \bigoplus_{\alpha \in \Delta_+ \setminus \Delta_\Sigma} \mfrak{g}_\alpha, \qquad \mfrak{l}_\Sigma = \mfrak{h} \oplus \bigoplus_{\alpha \in \Delta_\Sigma} \mfrak{g}_\alpha, \qquad \widebar{\mfrak{u}}_\Sigma = \bigoplus_{\alpha \in \Delta_+ \setminus \Delta_\Sigma} \mfrak{g}_{-\alpha}.
\end{align*}
We have also the corresponding triangular decomposition $\mfrak{g}= \widebar{\mfrak{u}}_\Sigma \oplus \mfrak{l}_\Sigma \oplus \mfrak{u}_\Sigma$ of the Lie algebra $\mfrak{g}$.
\medskip

We say that a $\mfrak{g}$-module $M$ is a \emph{weight} $\mfrak{g}$-module if $\mfrak{h}$ is diagonalizable on $M$. For $\lambda\in \mfrak{h}^*$, the vector subspace
\begin{align*}
  M_\lambda=\{v\in M;\, (\forall h \in \mfrak{h})\, (h-\lambda(h))v=0\}
\end{align*}
of $M$ is called the weight subspace with weight $\lambda$ provided $M_\lambda \neq \{0\}$. The dimension of $M_{\lambda}$ is the \emph{multiplicity} of the weight $\lambda$. For a weight $\mfrak{g}$-module $M$, we have $M=\oplus_{\lambda\in \mfrak{h}}M_{\lambda}$.

Let $\mfrak{p} = \mfrak{l} \oplus \mfrak{u}$ be a parabolic subalgebra of $\mfrak{g}$ and let $E$  a be a simple weight $\mfrak{l}$-module. Then the induced $\mfrak{g}$-module $M^\mfrak{g}_\mfrak{p}(E)$ is defined by
\begin{align*}
  M^\mfrak{g}_\mfrak{p}(E) = U(\mfrak{g}) \otimes_{U(\mfrak{p})}\! E,
\end{align*}
where $E$ is considered as the $\mfrak{p}$-module on which $\mfrak{u}$ acts trivially. In addition, we denote by $L^\mfrak{g}_\mfrak{p}(E)$ the unique simple quotient of $M^\mfrak{g}_\mfrak{p}(E)$.
If $E$ is a simple finite-dimensional $\mfrak{l}$-module with highest weight $\lambda \in \mfrak{h}^*$, we denote by $M^\mfrak{g}_\mfrak{p}(\lambda)$ the generalized Verma module $M^\mfrak{g}_\mfrak{p}(E)$ and by $L^\mfrak{g}_\mfrak{p}(\lambda)$ its unique simple quotient.

A weight $\mfrak{g}$-module $M$ is \emph{torsion free} if for any weight subspace $M_{\lambda}$ and any root $\alpha$, a nonzero root vector $a \in \mfrak{g}_{\alpha}$ defines an isomorphism
between $M_\lambda$ and $M_{\lambda+\alpha}$. Further, a weight $\mfrak{g}$-module $M$ is \emph {cuspidal} if $M$ is a finitely generated torsion free $\mfrak{g}$-module with finite weight multiplicities.

%%%%%%%%%%%%%%%%%%%%%%%%%%%%%%%%%%%%%%%%%%%%%%%%%%%%%%%%%%%%%%%%%%%%%%%%%%%%%%%%%%%%%%%%%%%

\subsection{Simple affine vertex algebras}

Let $\mfrak{g}$ be a simple Lie algebra and $\kappa$ be a $\mfrak{g}$-invariant symmetric bilinear form on $\mfrak{g}$. Since $\mfrak{g}$ is a simple Lie algebra, we have $\kappa = k\kappa_0$ for some $k \in \C$. The affine Kac--Moody algebra $\widehat{\mfrak{g}}_\kappa$ of level $\kappa$ associated to $\mfrak{g}$ is the $1$-dimensional universal central extension $\widehat{\mfrak{g}}_\kappa = \mfrak{g}(\!(t)\!) \oplus \C c$ of the formal loop algebra $\mfrak{g}(\!(t)\!)= \mfrak{g} \otimes_\C \C(\!(t)\!)$ with the commutation relations
\begin{align}
  [a \otimes f(t), b \otimes g(t)] = [a,b] \otimes f(t)g(t) - \kappa(a,b)\Res_{t=0} (f(t)dg(t))c \label{eq:commutation relation}
\end{align}
for $a, b \in \mfrak{g}$ and $f(t), g(t) \in \C(\!(t)\!)$, where $c$ is the central element of $\widehat{\mfrak{g}}_\kappa$. By introducing the notation $a_n = a \otimes t^n$ for $a \in \mfrak{g}$ and $n \in \Z$, the commutation relations \eqref{eq:commutation relation} can be simplified into the form
\begin{align}
  [a_m,b_n]=[a,b]_{m+n}+m \kappa(a,b) \delta_{m,-n} c \label{eq:commutation relation modes}
\end{align}
for $m,n \in \Z$ and $a,b \in \mfrak{g}$. As $\mfrak{h}$ is a Cartan subalgebra of $\mfrak{g}$, we introduce a Cartan subalgebra $\smash{\widehat{\mfrak{h}}}$ of $\widehat{\mfrak{g}}_\kappa$ by
\begin{align*}
  \widehat{\mfrak{h}} = \mfrak{h} \otimes_\C \C 1 \oplus \C c.
\end{align*}
For $\lambda \in \widehat{\mfrak{h}}^*$, we define its integral root system $\widehat{\Delta}(\lambda)$ by
\begin{align*}
 \widehat{\Delta}(\lambda) =\{\alpha\in \widehat{\Delta}^{\rm re};\, \bra \lambda+\widehat{\rho}, \alpha^\vee \ket \in \Z\},
\end{align*}
where $\widehat{\rho}=\rho+h^\vee \Lambda_0$. Further, let $\widehat{\Delta}(\lambda)_+ = \widehat{\Delta}(\lambda) \cap \widehat{\Delta}^{\rm re}_+$ be the set of positive roots of $\widehat{\Delta}(\lambda)$ and $\smash{\widehat{\Pi}}(\lambda) \subset \smash{\widehat{\Delta}}(\lambda)_+$ be the set of simple roots. Then we say that a weight $\lambda \in \smash{\widehat{\mfrak{h}}^*}$ is \emph{admissible} (\cite{Kac-Wakimoto1989}) provided
\begin{enumerate}[topsep=0pt,itemsep=0pt,parsep=0pt]
 \item[i)] $\lambda$ is \emph{regular dominant}, that is $\bra \lambda + \widehat{\rho}, \alpha^\vee \ket \notin -\N_0$ for all $\alpha\in \widehat{\Delta}^{\rm re}_+$;
 \item[ii)] the $\Q$-span of $\smash{\widehat{\Delta}}(\lambda)$ contains $\smash{\widehat{\Delta}^{\rm re}}$.
\end{enumerate}
In particular, if $\lambda=k\Lambda_0$ with $k \in \C$ is admissible, then $k$ is usually called an \emph{admissible number}. Admissible numbers were described in \cite{Kac-Wakimoto1989,Kac-Wakimoto2008} as follows. A number $k \in \C$ is admissible if and only if
\begin{align*}
  k+h^\vee ={p \over q} \ \text{with } p,q\in \N,\ (p,q)=1,\
  p\geq
    \begin{cases}
      h^\vee & \text{if $(r^\vee,q)=1$}, \\
      h & \text{if $(r^\vee,q)=r^\vee$},
    \end{cases}
\end{align*}
where $r^\vee$ is the lacing number of $\mfrak{g}$, i.e.\ the maximal number of edges in the Dynkin diagram of the Lie algebra $\mfrak{g}$. In this case we have $\smash{\widehat{\Pi}}(k\Lambda_0) = \{\dot{\alpha}_0,\alpha_1,\dots,\alpha_\ell\}$, where
\begin{align*}
  \dot{\alpha}_0 = \begin{cases}
    -\theta+q\delta & \text{if $(r^\vee,q)=1$}, \\
    -\theta_s+{q \over r^\vee}\delta & \text{if $(r^\vee,q)=r^\vee$}
  \end{cases}
\end{align*}
and $\ell = \rank \mfrak{g}$. Since the admissibility of a number $k \in \C$  depends only on $\mfrak{g}$, we shall say that $k$ is an admissible number for $\mfrak{g}$.
\medskip

Let $k \in \Q$ be an admissible number for $\mfrak{g}$ and let $\kappa=k\kappa_0$. We denote by $\mcal{V}_\kappa(\mfrak{g})$ the \emph{universal affine vertex algebra} associated to the affine Kac--Moody algebra $\widehat{\mfrak{g}}_\kappa$, we have
\begin{align*}
 \mcal{V}_\kappa(\mfrak{g})=U(\widehat{\mfrak{g}}_\kappa)_{U(\mfrak{g}[[t]]\oplus \C c)}\C v_\kappa,
\end{align*}
where $\mfrak{g}[[t]]v_\kappa=0$ and $cv_\kappa=v_\kappa$. The \emph{simple affine vertex algebra} $\mcal{L}_\kappa(\mfrak{g})$ associated to $\widehat{\mfrak{g}}_\kappa$ is the unique simple graded quotient of $\mcal{V}_\kappa(\mfrak{g})$. If we denote by $A(\mcal{L}_\kappa(\mfrak{g}))$ the \emph{Zhu's algebra} of $\mcal{L}_\kappa(\mfrak{g})$, then there is a one-to-one correspondence between simple positive energy $\mcal{L}_\kappa(\mfrak{g})$-modules and simple $A(\mcal{L}_\kappa(\mfrak{g}))$-modules, see \cite{Zhu1996}.

For a $\mfrak{g}$-module $E$, let us consider the induced $\widehat{\mfrak{g}}_\kappa$-module
\begin{align*}
   \mathbb{M}_{\kappa,\mfrak{g}}(E)=U(\widehat{\mfrak{g}}_\kappa)_{U(\mfrak{g}[[t]] \oplus \C c)} E,
\end{align*}
where $E$ is considered as the $\mfrak{g}[[t]] \oplus \C c$-module on which $\mfrak{g} \otimes_\C t\C[[t]]$ acts trivially and $c$ acts as the identity. Since $\mathbb{M}_{\kappa,\mfrak{g}}(E)$ has a unique maximal $\widehat{\mfrak{g}}_\kappa$-submodule having zero intersection with $E$, we denote by $\mathbb{L}_{\kappa,\mfrak{g}}(E)$ the corresponding quotient. We say that a $\mfrak{g}$-module $E$ is \emph{admissible of level $k$} if $\mathbb{L}_{\kappa,\mfrak{g}}(E)$ is an $\mcal{L}_\kappa(\mfrak{g})$-module, or equivalently if $E$ is an $A(\mcal{L}_\kappa(\mfrak{g}))$-module. As we have $A(\mcal{L}_\kappa(\mfrak{g})) \simeq U(\mfrak{g})/I_k$, where $I_k$ is a two-sided ideal of $U(\mfrak{g})$, we obtain that a $\mfrak{g}$-module $E$ is admissible of level $k$ if and only if the ideal $I_k$ is contained in the annihilator $\Ann_{U(\mfrak{g})}\!E$.

Admissible simple highest weight $\mfrak{g}$-modules of level $k$ were classified in \cite{Arakawa2016} as follows. Let us denote by ${\rm Pr}_k$ the set of admissible weights $\lambda \in \smash{\widehat{\mfrak{h}}^*}$ of level $k$ such that there is an element $y \in \smash{\widetilde{W}}$ of the extended affine Weyl group $\smash{\widetilde{W}}$ of $\mfrak{g}$ satisfying $\smash{\widehat{\Delta}(\lambda)} = \smash{y(\widehat{\Delta}(k\Lambda_0))}$. Besides, let us introduce the subset
\begin{align*}
  \widebar{{\rm Pr}}_k = \{\lambda \in \mfrak{h}^*;\, \lambda+k\Lambda_0 \in {\rm Pr}_k\}
\end{align*}
of $\mfrak{h}^*$, which is the canonical projection of ${\rm Pr}_k$ to $\mfrak{h}^*$. Then we have the following statement.
\medskip

\theorem{\cite{Arakawa2016}\label{thm-hweight}
Let $k \in \Q$ be an admissible number for $\mfrak{g}$. Then the simple highest weight $\mfrak{g}$-module $L^\mfrak{g}_\mfrak{b}(\lambda)$ with highest weight $\lambda \in \mfrak{h}^*$ is admissible of level $k$ if and only if $\lambda \in \widebar{{\rm Pr}}_k$.}

Hence, the main difficulty in the classification of admissible simple highest weight $\mfrak{g}$-modules consists in an explicit description of the set $\smash{\widebar{{\rm Pr}}_k}$ for an admissible number $k$ of $\mfrak{g}$ which will be our next goal.

%%%%%%%%%%%%%%%%%%%%%%%%%%%%%%%%%%%%%%%%%%%%%%%%%%%%%%%%%%%%%%%%%%%%%%%%%%%%%%%%%%%%%%%%%%%

\subsection{Primitive ideals and nilpotent orbits}

%%%%%%%%%%%%%%%%%%%%%%%%%%%%%%%%%%%%%%%%%%%%%%%%%%%%%%%%%%%%%%%%%%%%%%%%%%%%%%%%%%%%%%%%%%%

\subsubsection{Joseph's associated variety}
\label{subsec:Primitive ideals}

Let us consider a semisimple Lie algebra $\mfrak{g}$. We denote by
\begin{align*}
  U(\mfrak{g})^\mfrak{h} = \{a \in U(\mfrak{g});\, \ad(h)(a)=0\ \text{for all $h \in \mfrak{h}$}\}
\end{align*}
the centralizer of $\mfrak{h}$ in $U(\mfrak{g})$ and by
\begin{align*}
  \Upsilon \colon U(\mfrak{g})^\mfrak{h} \rarr U(\mfrak{h})
\end{align*}
the restriction of the projection $U(\mfrak{g}) = U(\mfrak{h}) \oplus (\widebar{\mfrak{n}}U(\mfrak{g}) + U(\mfrak{g})\mfrak{n}) \rarr U(\mfrak{h})$ to $U(\mfrak{g})^\mfrak{h}$. It is known that $\Upsilon$ is an algebra homomorphism. Let us note that the algebra homomorphism $\chi_\lambda \colon Z(\mfrak{g}) \rarr \C$ with $\lambda \in \mfrak{h}^*$ defined by
\begin{align*}
  \chi_\lambda(z) = \Upsilon(z)(\lambda)
\end{align*}
for $z \in Z(\mfrak{g})$ is the central character of $M^\mfrak{g}_\mfrak{b}(\lambda)$ and $L^\mfrak{g}_\mfrak{b}(\lambda)$.

Further, let us consider the PBW filtration on the universal enveloping algebra $U(\mfrak{g})$ of $\mfrak{g}$ and the associated graded algebra $\gr U(\mfrak{g}) \simeq S(\mfrak{g}) \simeq \C[\mfrak{g}^*]$. The \emph{associated variety} $\mcal{V}(I)$ of a left ideal $I$ of $U(\mfrak{g})$ is defined as the zero locus in $\mfrak{g}^*$ of the associated graded ideal $\gr I$ of $S(\mfrak{g})$.
We have
\begin{align}
  \mcal{V}(I) = \Specm (S(\mfrak{g})/\!\gr I) = \Specm(S(\mfrak{g})/\sqrt{\gr I}),
\end{align}
where $\sqrt{\gr I}$ denotes the radical of $\gr I$. For $\lambda \in \mfrak{h}^*$, we denote  by
\begin{align*}
  J_\lambda = \Ann_{U(\mfrak{g})}\! L^\mfrak{g}_\mfrak{b}(\lambda)
\end{align*}
a primitive ideal of $U(\mfrak{g})$. A theorem of Duflo \cite{Duflo1977} states that for any primitive ideal $I$ of $U(\mfrak{g})$ there exists $\lambda \in \mfrak{h}^*$ such that $I=J_\lambda$. This implies that a simple $\mfrak{g}$-module $E$ is admissible of level $k$ if and only if $\Ann_{U(\mfrak{g})}\!E=J_\lambda$ for some $\lambda \in \widebar{{\rm Pr}}_k$.
\medskip

\proposition{\cite[Section 8.5.8]{Dixmier1977-book}\label{prop:max ideal} Let $\lambda \in \mfrak{h}^*$ be a dominant weight. Then $J_\lambda$ is the unique maximal two-sided ideal of $U(\mfrak{g})$ containing $U(\mfrak{g})\ker \chi_\lambda$.}

Let us note that if $J_\mu \subset J_\lambda$ for $\lambda,\mu \in \mfrak{h}^*$, then necessarily $\mu \in W\!\cdot \lambda$. The opposite implication is given by the following statement, which is an easy consequence of Proposition \ref{prop:max ideal}.
\medskip

\proposition{\label{prop:primitive ideal} Let $\lambda,\mu \in \mfrak{h}^*$ be dominant weights. Then $J_\lambda = J_\mu$ if and only if $\mu \in W\! \cdot \lambda$.}

We have the following statement.

\lemma{\label{lem:admissibility} Let $k \in \Q$ be an admissible number for $\mfrak{g}$. Further, let $\lambda \in \smash{\widebar{{\rm Pr}}_k}$ and $\mu \in W\! \cdot \lambda$. Then we have $\mu \in \smash{\widebar{{\rm Pr}}_k}$ if and only if $\mu$ is a dominant weight.}

\proof{If $\mu \in \smash{\widebar{{\rm Pr}}_k}$, then by definition the weight $\mu$ is dominant. On the other hand, if $\mu$ is dominant and $\mu \in W\!\cdot \lambda$, then by Proposition \ref{prop:primitive ideal} we obtain that $J_\lambda = J_\mu$, since $\lambda$ is also dominant. This gives us $\mu \in \smash{\widebar{{\rm Pr}}_k}$.}

For a weight $\lambda \in \mfrak{h}^*$, we introduce subsets $W^\lambda$ and $W_\lambda$ of $W$ by
\begin{align*}
    W^\lambda = \{w \in W;\, \Delta_+ \cap w(-\Delta_+) \subset \Delta_+ \setminus \Delta(\lambda)\} \quad \text{and} \quad
    W_\lambda = \langle s_\alpha;\, \alpha \in \Delta(\lambda) \rangle,
\end{align*}
where $\Delta(\lambda) = \{\alpha\in\Delta;\, \bra \lambda+\rho,\alpha^\vee \ket \in \Z\}$ is the integral root system of $\lambda$. Obviously, the set $W_\lambda$ is the Weyl group associated to the integral root system $\Delta(\lambda)$. Let us note that $W^\lambda$ and $W_\lambda$ depend only on $\Delta(\lambda)$ not on $\lambda$ itself. We will be mainly interested in the case when the integral root system $\Delta(\lambda)$ is generated by a subset of $\Pi$.

\medskip

\lemma{Let $k\in \Q$ be an admissible number for $\mfrak{g}$. Further, let $\lambda \in \smash{\widebar{{\rm Pr}}_k}$ and $\mu \in W\! \cdot \lambda$. Then we have $\mu \in \smash{\widebar{{\rm Pr}}_k}$ if and only if $\mu = w^{-1} \cdot \lambda$ for some element $w \in W^\lambda$.}

\proof{By Lemma \ref{lem:admissibility} we know that $\mu \in \smash{\widebar{{\rm Pr}}_k}$ if and only if $\mu$ is a dominant weight. Hence, we need to show that $\mu$ is dominant if and only if $\mu = w^{-1} \cdot \lambda$ for some element $w \in W^\lambda$. By the assumption we have that $\mu = w^{-1} \cdot \lambda$ for some $w \in W$. We may write
\begin{align*}
  \bra \mu+\rho,\alpha^\vee \ket = \bra w^{-1}(\lambda+\rho),\alpha^\vee \ket = \bra \lambda+\rho, w(\alpha)^\vee \ket
\end{align*}
for $\alpha \in \Delta_+$. Therefore, the condition that $\mu$ is dominant is equivalent to
\begin{align*}
  \bra \lambda+\rho,\beta^\vee \ket \notin -\N\ \text{for}\ \beta \in \Delta_+ \cap w(\Delta_+) \quad \text{and} \quad \bra \lambda+\rho,\beta^\vee \ket \notin \N\ \text{for}\ \beta \in \Delta_+ \cap w(-\Delta_+).
\end{align*}
However, the first condition is always satisfied since $\lambda$ is dominant. On the other hand, the second condition is satisfied if and only if $\Delta_+ \cap w(-\Delta_+) \subset \Delta_+ \setminus \Delta(\lambda)$, or in other words if and only if $w \in W^\lambda$.}

Let $\mfrak{p}$ be a standard parabolic subalgebra of $\mfrak{g}$ associated to a subset $\Sigma$ of $\Pi$. We set
\begin{align}
  \widebar{{\rm Pr}}_k(\mfrak{p}) = \{\lambda \in \widebar{{\rm Pr}}_k;\, \Delta(\lambda) = \Delta_\Sigma\}, \label{eq:Pr_k(p) definition}
\end{align}
where $\Delta_\Sigma$ is the root subsystem of $\Delta$ generated by $\Sigma$. Then by \cite{Jantzen1977} we have that $L^\mfrak{g}_\mfrak{p}(\lambda) \simeq M^\mfrak{g}_\mfrak{p}(\lambda)$ for $\lambda \in \smash{\widebar{{\rm Pr}}}_k(\mfrak{p})$. Further, since the subsets $W^\lambda$ and $W_\lambda$ of $W$ for $\lambda \in \mfrak{h}^*$ satisfying $\Delta(\lambda) = \Delta_\Sigma$ depend only on $\Sigma$, we shall denote them by $W^\mfrak{p}$ and $W_\mfrak{p}$, respectively. The subgroup $W_\mfrak{p}$ of $W$ we may identify with the Weyl group of the Levi subalgebra $\mfrak{l}$ of $\mfrak{p}$. The set $W^\mfrak{p}$ is the set of minimal length coset representatives of $W_\mfrak{p}\backslash W$. In the literature, the set $W^\mfrak{p}$ is often characterized in different equivalent ways. In fact, for $w \in W$, we have that the following statements are equivalent:
\begin{enumerate}[topsep=3pt,itemsep=0pt]
  \item[i)] $\Delta_+ \cap w(-\Delta_+) \subset \Delta_+ \setminus \Delta_\Sigma$;
  \item[ii)] $w^{-1}(\Delta_+ \cap \Delta_\Sigma) \subset \Delta_+$;
  \item[iii)] $\langle w(\rho), \alpha^\vee \rangle \in \N$ for all $\alpha \in \Delta_+ \cap \Delta_\Sigma$;
  \item[iv)] $\ell(s_\alpha w) = \ell(w)+1$ for all $\alpha \in \Sigma$.
\end{enumerate}

\medskip

The elements of $W^\mfrak{p}$ can be determined by using the following observation, see \cite[Proposition 3.2.16]{Cap-Slovak2009-book}.

\medskip

\proposition{The mapping $w \mapsto w^{-1}(\rho^\mfrak{p})$ restricts to a bijection between $W^\mfrak{p}$ and the orbit of $\rho^\mfrak{p}$ under $W$, where $\rho^\mfrak{p} = \sum_{\alpha \in \Pi \setminus \Sigma} \omega_\alpha$.}

The following statement is clear.

\lemma{\label{lem:W^p condition}Let $\mfrak{p}$ be a standard parabolic subalgebra of $\mfrak{g}$ and $k\in \Q$ be an admissible number for $\mfrak{g}$. Furthermore, let $\lambda \in \smash{\widebar{{\rm Pr}}_k}(\mfrak{p})$ and $\mu \in W\! \cdot \lambda$. Then $\mu \in \smash{\widebar{{\rm Pr}}_k}$ if and only if $\mu = w^{-1} \cdot \lambda$ for some element $w \in W^\mfrak{p}$.}

\vspace{-2mm}

%%%%%%%%%%%%%%%%%%%%%%%%%%%%%%%%%%%%%%%%%%%%%%%%%%%%%%%%%%%%%%%%%%%%%%%%%%%%%%%%%%%%%%%%%%

\subsubsection{Nilpotent orbits}

Let $G$ be a complex connected semisimple algebraic group with its Lie algebra $\mfrak{g}$. We denote by $\mcal{N}(\mfrak{g})$ the nilpotent cone of $\mfrak{g}$, i.e.\ the set of nilpotent elements of $\mfrak{g}$. It is an irreducible closed algebraic subvariety of $\mfrak{g}$ and a finite union of $G$-orbits. There is a unique nilpotent orbit of $\mfrak{g}$, denoted by $\mcal{O}_{\rm prin}$ and called the \emph{principal nilpotent orbit} of $\mfrak{g}$, which is a dense open subset of $\mcal{N}(\mfrak{g})$. Next, since $\mfrak{g}$ is simple, there exists a unique nilpotent orbit of $\mfrak{g}$ that is a dense open subset of $\mcal{N}(\mfrak{g}) \setminus \mcal{O}_{\rm prin}$, denoted by $\mcal{O}_{\rm subreg}$ and called the \emph{subregular nilpotent orbit} of $\mfrak{g}$. Besides, there is a unique nonzero nilpotent orbit of $\mfrak{g}$ of minimal dimension, denoted by $\mcal{O}_{\rm min}$ and called the \emph{minimal nilpotent orbit} of $\mfrak{g}$, such that it is contained in the closure of all nonzero nilpotent orbits of $\mfrak{g}$. By $\mcal{O}_{\rm zero}$ we denote the \emph{zero nilpotent orbit} of $\mfrak{g}$.
For the dimension of these distinguished nilpotent orbits of $\mfrak{g}$ see Figure \ref{fig:hasse diagram general}.

\begin{figure}[ht]
\centering
{\begin{tikzpicture}
[yscale=1.1,xscale=1.7,vector/.style={circle,draw=white,fill=black,ultra thick, inner sep=0.8mm},vector2/.style={circle,draw=white,fill=white,ultra thick, inner sep=1mm}]
\begin{scope}
  \node (A) at (0,0)  {$\mcal{O}_{{\rm zero}}$};
  \node (B) at (0,1)  {$\mcal{O}_{{\rm min}}$};
  \node (C1) at (0,1.7)  {};
  \node (C) at (0,2)  {$\,\dots$};
  \node (C2) at (0,2.3) {};
  \node (D) at (0,3)  {$\mcal{O}_{{\rm subreg}}$};
  \node (E) at (0,4)  {$\mcal{O}_{{\rm prin}}$};
  \node at (2,0) {$0$};
  \node at (2,1) {$2h^\vee - 2$};
  \node at (2,2) {$\dots$};
  \node at (2,3) {$\dim \mfrak{g} - \rank \mfrak{g} - 2$};
  \node at (2,4) {$\dim \mfrak{g} - \rank \mfrak{g}$};
  \draw [thin, -] (B) -- (C1);
  \draw [thin, -] (D) -- (C2);
  \draw [thin, -] (A) -- (B);
  \draw [thin, -] (D) -- (E);
  \node at (0,4.7) {nilpotent orbit};
  \node at (2,4.7) {dimension};
\end{scope}
\end{tikzpicture}}
\caption{Hasse diagram of nilpotent orbits of $\mfrak{g}$}
\label{fig:hasse diagram general}
\vspace{-2mm}
\end{figure}

If $I$ is a two-sided ideal of $\mfrak{g}$, then $I$ and $\gr I$ are invariant under the adjoint action of $G$. Consequently, the associated variety is a union of $G$-orbits of $\mfrak{g}^*$.  As the Cartan--Killing form $\kappa_\mfrak{g}$ is a $\mfrak{g}$-invariant symmetric bilinear form on $\mfrak{g}$, it provides a one-to-one correspondence between adjoint orbits of $\mfrak{g}$ and coadjoint orbits of $\mfrak{g}^*$. For an adjoint orbit $\mcal{O}$ of $\mfrak{g}$ we denote by $\mcal{O}^*$ the corresponding coadjoint orbit of $\mfrak{g}^*$. In addition, for a primitive ideal $I$ of $U(\mfrak{g})$ the associated variety $\mcal{V}(I)$ is the closure of $\mcal{O}^*$ for some nilpotent orbit $\mcal{O}$ of $\mfrak{g}$, see \cite{Joseph1985}.
\medskip

\definition{Let $E$ be a $\mfrak{g}$-module. We say that $E$ belongs to a nilpotent orbit $\mcal{O}$ of $\mfrak{g}$ if the associated variety $\mcal{V}(\Ann_{U(\mfrak{g})}\!E)$ is the closure of $\mcal{O}^*$.}

\theorem{\cite[Theorem 9.3]{Arakawa2015b}\label{thm:nilpotent orbits}
Let $k \in \Q$ be an admissible number for $\mfrak{g}$ with denominator $q \in \N$. Then there exists a nilpotent orbit $\mcal{O}_q$ of $\mfrak{g}$ such that
\begin{align*}
  \mcal{V}(I_k) = \widebar{{\mcal{O}}_q^*}
\end{align*}
and we have
\begin{align*}
  \widebar{{\mcal{O}}_q} = \begin{cases}
    \{x\in\mfrak{g};\, (\ad x)^{2q}=0\}  & \text{if $(r^\vee,q)=1$}, \\
    \{x\in\mfrak{g};\, \pi_{\theta_s}\!(x)^{2q/r^\vee}\!=0\}  & \text{if $(r^\vee,q)=r^\vee$},
  \end{cases}
\end{align*}
where $\pi_{\theta_s} \colon \mfrak{g} \rarr \End L^\mfrak{g}_\mfrak{b}(\theta_s)$ is the simple finite-dimensional $\mfrak{g}$-module with highest weight $\theta_s$.}

Let $\mcal{O}$ be a nilpotent orbit of $\mfrak{g}$ and let $k \in \Q$ be an admissible number for $\mfrak{g}$ with denominator $q \in \N$. We define the subset
\begin{align}
  \widebar{{\rm Pr}}_k^\mcal{O} =\{\lambda\in \widebar{{\rm Pr}}_k;\, \mcal{V}(J_\lambda) = \widebar{\mcal{O}^*}\}
\end{align}
of $\mfrak{h}^*$. It easily follows that a simple $\mfrak{g}$-module $E$ belonging to the nilpotent orbit $\mcal{O}$ is admissible of level $k$ if and only if $\Ann_{U(\mfrak{g})}\!E=J_\lambda$ for some $\lambda \in \widebar{{\rm Pr}}_k^\mcal{O}$. Moreover, as $I_k \subset J_\lambda$ for $\lambda \in \smash{\widebar{{\rm Pr}}_k}$, we have $\mcal{V}(J_\lambda) \subset \mcal{V}(I_k) = \smash{\widebar{{\mcal{O}}_q^*}}$ by Theorem \ref{thm:nilpotent orbits}, which gives us a decomposition
\begin{align}
  \widebar{{\rm Pr}}_k = \bigsqcup_{\mcal{O} \subset \widebar{{\mcal{O}}_q}}  \widebar{{\rm Pr}}_k^\mcal{O}. \label{eq:Pr_k decomposition}
\end{align}
Hence, we need to describe the subset $\widebar{{\rm Pr}}_k^\mcal{O}$ of $\widebar{{\rm Pr}}_k$ for a nilpotent orbit $\mcal{O}$ of $\mfrak{g}$.

In addition, from the decomposition \eqref{eq:Pr_k decomposition} it follows immediately that the set $\widebar{{\rm Pr}}_k^{\smash{\mcal{O}_{\rm prin}}}$ is non-empty if and only if $\mcal{O}_{\rm prin} = \mcal{O}_q$, or equivalently if and only if
\begin{align*}
    q \geq \begin{cases}
             h & \text{if $(r^\vee,q)=1$},\\
             {}^L h^\vee r^\vee & \text{if $(r^\vee,q)=r^\vee$},
           \end{cases}
\end{align*}
where $h$ is the Coxeter number of $\mfrak{g}$ and ${}^L h^\vee$ is the dual Coxeter number of the Langlands dual Lie algebra ${}^L\mfrak{g}$ of $\mfrak{g}$ (see \cite{Arakawa-Futorny-Ramirez2017}).

%%%%%%%%%%%%%%%%%%%%%%%%%%%%%%%%%%%%%%%%%%%%%%%%%%%%%%%%%%%%%%%%%%%%%%%%%%%%%%%%%%%%%%%%%%

\subsubsection{Richardson orbits}
\label{subsec:Richardson orbits}

For $x \in \mfrak{g}$, we denote by $\mfrak{g}^x$ the centralizer of $x$ in $\mfrak{g}$. Besides, for a subset $X$ of $\mfrak{g}$ we define the set
\begin{align*}
  X^{\rm reg} = \{x \in X;\, \dim \mfrak{g}^x = \min\nolimits_{y \in X} \dim \mfrak{g}^y \}
\end{align*}
and call it the \emph{set of regular elements} in $X$.
\medskip

\theorem{\cite[Corollary 4.7]{Borho-Brylinski1982}\label{thm:annihilator verma module} Let $\mfrak{p}$ be a standard parabolic subalgebra of $\mfrak{g}$. Then for a weight $\lambda \in \Lambda^+(\mfrak{p})$ the associated variety $\mcal{V}(\Ann_{U(\mfrak{g})}\!M^\mfrak{g}_\mfrak{p}(\lambda))$ is the closure of the nilpotent orbit $\mcal{O}^*_\mfrak{p}$ of $\mfrak{g}^*$, where
\begin{align*}
  \mcal{O}_\mfrak{p} = (G.\mfrak{p}^\perp)^{\rm reg}
\end{align*}
and $\mfrak{p}^\perp$ is the orthogonal complement of $\mfrak{p}$ with respect to the Cartan--Killing form $\kappa_\mfrak{g}$. In particular, we have
\begin{align*}
  \widebar{\mcal{O}_{\mfrak{p}}} = G.\mfrak{p}^\perp
\end{align*}
and the associated variety $\mcal{V}(\Ann_{U(\mfrak{g})}\!M^\mfrak{g}_\mfrak{p}(\lambda))$ is irreducible. The nilpotent orbit $\mcal{O}_\mfrak{p}$ is the \emph{Richardson orbit} attached to $\mfrak{p}$.}

Let us recall that when two parabolic subalgebras $\mfrak{p}_1$ and $\mfrak{p}_2$ of $\mfrak{g}$ have conjugate Levi subalgebras, the corresponding Richardson orbits $\mcal{O}_{\mfrak{p}_1}$ and $\mcal{O}_{\mfrak{p}_2}$ are exactly the same. Moreover, if $\mfrak{p}_1$ and $\mfrak{p}_2$ are the standard parabolic subalgebras associated to subsets $\Sigma_1$ and $\Sigma_2$ of $\Pi$, then
the Richardson orbits $\mcal{O}_{\mfrak{p}_1}$ and $\mcal{O}_{\mfrak{p}_2}$ coincide if and only if the root subsystems $\Delta_{\Sigma_1}$ and $\Delta_{\Sigma_2}$ of $\Delta$ are conjugate by an element of the Weyl group $W$ of $\mfrak{g}$. Therefore, we may introduce an equivalence relation on the set of all standard parabolic subalgebras of $\mfrak{g}$ by
\begin{align*}
  \mfrak{p}_1 \sim \mfrak{p}_2  \Longleftrightarrow  \text{there exits $w \in W$ such that $\Delta_{\Sigma_1} = w(\Delta_{\Sigma_2})$}.
\end{align*}
For a standard parabolic subalgebra $\mfrak{p}$ of $\mfrak{g}$, we shall denote by $[\mfrak{p}]$ the corresponding equivalence class.

%%%%%%%%%%%%%%%%%%%%%%%%%%%%%%%%%%%%%%%%%%%%%%%%%%%%%%%%%%%%%%%%%%%%%%%%%%%%%%%%%%%%%%%%%%%
%%%%%%%%%%%%%%%%%%%%%%%%%%%%%%%%%%%%%%%%%%%%%%%%%%%%%%%%%%%%%%%%%%%%%%%%%%%%%%%%%%%%%%%%%%%

\section{Admissible highest weight modules for $\mfrak{sl}_{n+1}$}
\label{sec:adm highest weight modules}

Let us consider the simple Lie algebra $\mfrak{g}=\mfrak{sl}_{n+1}$ for $n \in \N$ and let $k \in \Q$ be an admissible number for $\mfrak{g}$ with denominator $q \in \N$. In this section we refine the  classification of highest weight modules over the simple affine vertex algebra $\mcal{L}_\kappa(\mfrak{g})$ of admissible level $\kappa=k\kappa_0$.

%%%%%%%%%%%%%%%%%%%%%%%%%%%%%%%%%%%%%%%%%%%%%%%%%%%%%%%%%%%%%%%%%%%%%%%%%%%%%%%%%%%%%%%%%%%

\subsection{Explicit description of admissible highest weights}

The nilpotent orbits of $\mfrak{g}$ are parameterized by the set $\mcal{P}_{n+1}$ of partitions of $n+1$. We denote by $\mcal{O}_\lambda$ the nilpotent orbit of $\mfrak{g}$ corresponding to a partition $\lambda \in \mcal{P}_{n+1}$. By \cite[Corollary 7.2.4]{Collingwood-McGovern1993-book} the dimension of $\mcal{O}_\lambda$ is given by
\begin{align*}
  \dim \mcal{O}_\lambda = (n+1)^2 - \sum_{i=1}^r (\lambda_i^t)^2,
\end{align*}
where $\lambda^t$ is the transpose partition of $\lambda$ and $r=\ell(\lambda^t)$ is the length of $\lambda^t$. Besides, $\sum_{i=1}^r (\lambda_i^t)^2 - 1$ is the dimension of $\mfrak{g}^e$ for $e \in \mcal{O}_\lambda$ and so we have $\dim \mcal{O}_\lambda = \dim \mfrak{g} - \dim \mfrak{g}^e$.

For a partition $\lambda = [\lambda_1,\lambda_2,\dots,\lambda_r] \in \mcal{P}_{n+1}$, we denote by $\mfrak{p}_\lambda$ the standard parabolic subalgebra of $\mfrak{g}$ associated to the subset
\begin{align*}
  \Sigma_\lambda = \{\alpha_1,\dots,\alpha_{\lambda_1-1},\alpha_{\lambda_1+1},\dots, \alpha_{\lambda_1+\lambda_2-1},\dots,\alpha_{\lambda_1+\lambda_2+\dots+\lambda_{r-1}+1},\dots, \alpha_{\lambda_1+\lambda_2+\dots+\lambda_r-1} \}
\end{align*}
of $\Pi=\{\alpha_1,\alpha_2,\dots,\alpha_{n}\}$. In fact, we can define $\mfrak{p}_\lambda$ in a more elegant way as follows. Let $V=\C^{n+1}$ be the standard representation of $\mfrak{g}$ and let $\{e_1,e_2,\dots,e_{n+1}\}$ be the canonical basis of $V$. We define a partial flag $(V_0,V_1,V_2,\dots,V_r)$ in $V$ by setting
\begin{align*}
  V_i = \bigoplus_{k=1}^{s_i} \C e_k \quad \text{with} \quad s_i = \sum_{j=1}^i \lambda_j
\end{align*}
for $i=1,2,\dots,r$ and $V_0 = 0$. Then one checks easily that
\begin{align*}
  \mfrak{p}_\lambda = \{a \in \mfrak{g};\, a(V_i) \subset V_i\ \text{for}\ i=1,2,\dots, r\}.
\end{align*}
The Richardson orbit $\mcal{O}_{\mfrak{p}_\lambda}\!$ attached to $\mfrak{p}_\lambda$ is the nilpotent orbit $\mcal{O}_{\lambda^t}$, see \cite[Theorem 7.2.3]{Collingwood-McGovern1993-book}. This has an important consequence that every nilpotent orbit of $\mfrak{g}$ is a Richardson orbit.
\medskip

The nilpotent orbit $\mcal{O}_q$ from Theorem \ref{thm:nilpotent orbits} is given by
\begin{align*}
  \mcal{O}_q = \mcal{O}_{\lambda_q},
\end{align*}
where $\lambda_q \in \mcal{P}_{n+1}$ is the partition defined through $\lambda_q = [q^r,s]$, where $n+1=qr+s$ with $r,s\in \N_0$ and $0 \leq s \leq q-1$. Hence, we immediately get that
\begin{align*}
  \mcal{O}_q = \begin{cases}
    \mcal{O}_{\rm prin} = \mcal{O}_{[n+1]} &  \text{if $q \geq n+1$}, \\
    \mcal{O}_{\rm subreg} = \mcal{O}_{[n,1]} & \text{if $q=n$}, \\
    \mcal{O}_{\rm zero} = \mcal{O}_{[1^{n+1}]} & \text{if $q=1$}
  \end{cases}
\end{align*}
for $n \geq 1$. Let us introduce the subset
\begin{align*}
  {\rm Pr}_{k,\Z} = \{\lambda \in {\rm Pr}_k;\, \bra \lambda, \alpha^\vee \ket \in \Z \text{ for $\alpha \in \Pi$}\}
\end{align*}
of admissible weights of level $k$. Then by \cite{Kac-Wakimoto1989} we have
\begin{align*}
  {\rm Pr}_{k,\Z} = \{\lambda \in \widehat{\mfrak{h}}^*;\, \lambda(c)=k,\, \bra \lambda, \alpha^\vee \ket \in \N_0\ \text{for $\alpha \in \Pi$},\, \bra \lambda, \theta^\vee \ket \leq p-n-1\},
\end{align*}
where $p = (k+n+1)q$. Further, from \cite{Kac-Wakimoto1989} we get that
\begin{align*}
  {\rm Pr}_k = \bigcup_{\substack{y \in \widetilde{W}, \\ y(\widehat{\Delta}(k\Lambda_0)_+) \subset \widehat{\Delta}^{\rm re}_+}} {\rm Pr}_{k,y}, \quad {\rm Pr}_{k,y} = \{y\cdot \lambda;\, \lambda \in {\rm Pr}_{k,\Z}\},
\end{align*}
where the extended affine Weyl group $\smash{\widetilde{W}}$ of $\mfrak{g}$ is defined as $\smash{\widetilde{W}} = W \ltimes P^\vee$. For $\eta \in P^\vee$, we denote by $t_\eta$ the corresponding element of the group $\smash{\widetilde{W}}$.
The action of $t_\eta$ on $\smash{\widetilde{\mfrak{h}}^*}$ is given by
\begin{align*}
  t_\eta(\lambda) = \lambda + (\lambda,\delta)\eta - \bigg({(\eta,\eta) \over 2}\,(\lambda,\delta) + (\lambda,\eta)\!\!\bigg)\delta
\end{align*}
for $\lambda \in \smash{\widetilde{\mfrak{h}}^*}$. For $y,y' \in \smash{\widetilde{W}}$ satisfying $y(\widehat{\Delta}(k\Lambda_0)_+) \subset \widehat{\Delta}_+^{\rm re}$, $y'(\widehat{\Delta}(k\Lambda_0)_+) \subset \widehat{\Delta}_+^{\rm re}$ and $y' \neq y$ we have
\begin{align}
  {\rm Pr}_{k,y} \cap {\rm Pr}_{k,y'} \neq \emptyset \quad \Longleftrightarrow \quad  {\rm Pr}_{k,y} = {\rm Pr}_{k,y'} \quad \Longleftrightarrow \quad y' = yt_{q\omega_j} w_j \label{eq:Pr equality}
\end{align}
for some $j \in \{1,2,\dots,n\}$, where $w_j$ is the unique element of the Weyl group $W$ of $\mfrak{g}$ preserving the set $\{\alpha_1,\alpha_2,\dots,\alpha_n,-\theta\}$ and satisfying $w_j(-\theta) = \alpha_j$. Let us note that $w_j$ for $j=1,2,\dots,n$ gives rise to an automorphism of the extended Dynkin diagram of $\smash{\mfrak{g}}$. We denote by $W_+$ the subgroup of $W$ generated by the elements $w_1,w_2,\dots,w_n$.
\medskip

\lemma{\label{lem:w_j properties} We have
\begin{enumerate}[topsep=3pt,itemsep=1pt]
  \item[i)] $w_1 = s_1s_2\cdots s_n$;
  \item[ii)] $w_1(\alpha_i) = \alpha_{i+1}$ for $i=1,2,\dots,n-1$, $w_1(\alpha_n) = -\theta$, and $w_1(-\theta) = \alpha_1$;
  \item[iii)] $w_j = w_1^j$ for $j=1,2,\dots,n$;
  \item[iv)] $w_1^{n+1} = e$;
  \item[v)] $W_+ = \{e,w_1,w_2,\dots,w_n\}$;
  \item[vi)] $w_j \big(\!\sum_{i=1}^n\! \lambda_i \omega_i\big) = \sum_{i=1}^{j-1}\! \lambda_{n+1+i-j}\omega_i -\big(\!\sum_{i=1}^n\!\lambda_i \big) \omega_j + \sum_{i=j+1}^n\! \lambda_{i-j}\omega_i$ for $j=1,2,\dots,n$, where $\lambda_1,\lambda_2,\dots,\lambda_n \in \C$.
  %$w_j \cdot \lambda = \sum_{i=1}^{j-1} \bra\lambda,\alpha_{n+1+i-j}^\vee\ket\omega_i -(\bra\lambda,\theta^\vee\ket+h^\vee)\omega_j + \sum_{i=j+1}^n \bra\lambda,\alpha_{i-j}^\vee\ket\omega_i$
\end{enumerate}
}

\proof{Let us denote $w=s_1s_2\cdots s_n$. To prove (i) and (iii) we only need to verify that the element $w^j$ for $j \in \{1,2,\dots,n\}$ satisfies $w^j(-\theta)=\alpha_j$ and preserves the set $\{\alpha_1,\alpha_2,\dots,\alpha_n,-\theta\}$, since $w_j$ is uniquely determined by these two conditions. The rest of the statements is a straightforward computation. We left the details to the reader.}

Let us recall that an element $w \in W$ is called a \emph{Coxeter element} if a reduced expression of $w$ contains every simple reflection exactly once. Hence, by Lemma \ref{lem:w_j properties} we see that $w_1$ and $w_n$ are Coxeter elements. The order of a Coxeter element is the Coxeter number $h$ of $\mfrak{g}$, i.e.\ we have $h=n+1$.
\medskip

On the set of all standard parabolic subalgebras of $\mfrak{g}$ we introduce another equivalence relation given by
\begin{align*}
  \mfrak{p}_1 \sim_+ \mfrak{p}_2  \Longleftrightarrow  \text{there exits $w \in W_+$ such that $\Sigma_{\mfrak{p}_1} = w(\Sigma_{\mfrak{p}_2})$}.
\end{align*}
For a standard parabolic subalgebra $\mfrak{p}$ of $\mfrak{g}$, we shall denote by $[\mfrak{p}]_+$ the corresponding equivalence class. Obviously, we have that $[\mfrak{q}]_+ \subset [\mfrak{p}]$ for any $\mfrak{q} \in [\mfrak{p}]$. Besides, we assign to $\mfrak{p}$ a subgroup $W_+^\mfrak{p}$ of $W_+$ defined by
\begin{align*}
  W_+^\mfrak{p} = \{w \in W_+;\, w(\Sigma_\mfrak{p}) = \Sigma_\mfrak{p}\}.
\end{align*}
Let us note that we have $W_+^\mfrak{b}=W_+$ and $W_+^\mfrak{g}=\{e\}$.
\medskip

Since the dot action of the Weyl group $W$ on $\mfrak{h}^*$ gives rise to an equivalence relation on $\mfrak{h}^*$, we may define an equivalence relation on $\smash{\widebar{{\rm Pr}}}_k$ for an admissible number $k$ of $\mfrak{g}$ by
\begin{align*}
  \lambda \sim \mu \Longleftrightarrow \text{there exits $w \in W$ such that $\mu = w \cdot \lambda$}
\end{align*}
and set
\begin{align*}
  [\widebar{{\rm Pr}}_k] = \widebar{{\rm Pr}}_k /\!\sim\!.
\end{align*}
Moreover, if $\lambda, \mu \in \widebar{{\rm Pr}}_k$ then $J_\lambda=J_\mu$ if and only if there exists $w\in W$ such that $ \mu=w\cdot \lambda$, which is a consequence of Proposition \ref{prop:primitive ideal}. By using this fact together with \eqref{eq:Pr_k decomposition} we get a decomposition
\begin{align}
  [\widebar{{\rm Pr}}_k] = \bigsqcup_{\mcal{O} \subset \widebar{{\mcal{O}}_q}}  [\widebar{{\rm Pr}}_k^\mcal{O}], \label{eq:Pr_k decomposition class}
\end{align}
where $[\widebar{{\rm Pr}}_k^\mcal{O}]$ is the image of $\widebar{{\rm Pr}}_k^\mcal{O}$ in $[\widebar{{\rm Pr}}_k]$. Let us recall that by \cite[Proposition 2.8]{Arakawa-Futorny-Ramirez2017}
we have
\begin{align*}
  [\widebar{{\rm Pr}}_k] = \bigcup_{\substack{\eta \in P_+^\vee, \\ (\eta,\theta) \leq q-1}} [\widebar{{\rm Pr}}_{k,t_{-\eta}}],
\end{align*}
where $P_+^\vee$ is the set of dominant coweights of $\mfrak{g}$. Further, for a standard parabolic subalgebra $\mfrak{p}$ of $\mfrak{g}$, we define a subset $\Lambda_k(\mfrak{p})$ of $\smash{\widebar{{\rm Pr}}}_k$ by
\begin{align*}
  \Lambda_k(\mfrak{p}) = \bigcup_{\substack{\eta \in P_+^\vee, \\ (\eta,\theta)\leq q-1,\, \Pi^\eta = \Sigma_\mfrak{p}}} \widebar{{\rm Pr}}_{k,t_{-\eta}},
\end{align*}
where $\Pi^\eta = \{\alpha \in \Pi;\, (\eta,\alpha)=0\}$, which enables us to write
\begin{align}
  [\widebar{{\rm Pr}}_k] = \bigcup_{\Sigma \subset \Pi}\ [\Lambda_k(\mfrak{p}_\Sigma)]. \label{eq:[Pr_k] deconposition}
\end{align}
As we have $\smash{\widebar{{\rm Pr}}}_{k,t_{-\eta}}=\{\lambda-(k+n+1)\eta;\, \lambda \in \widebar{{\rm Pr}}_{k,\Z}\}$ for $\eta \in P_+^\vee$ satisfying $(\eta,\theta) \leq q-1$, we get
\begin{align}
  \Lambda_k(\mfrak{p}) = \{\lambda-(k+n+1)\eta;\, \lambda \in \widebar{{\rm Pr}}_{k,\Z},\, \eta \in P_+^\vee,\, (\eta,\theta) \leq q-1,\, \Pi^\eta = \Sigma_\mfrak{p}\}.
\end{align}
Moreover, for a weight $\lambda \in \Lambda_k(\mfrak{p})$ we have $\Delta(\lambda) = \Delta_{\Sigma_\mfrak{p}}\!$, which by \eqref{eq:Pr_k(p) definition} gives us $\Lambda_k(\mfrak{p}) \subset \smash{\widebar{{\rm Pr}}}_k(\mfrak{p})$. By \cite{Jantzen1977} we get that $M^\mfrak{g}_\mfrak{p}(\lambda) \simeq L^\mfrak{g}_\mfrak{b}(\lambda)$ for $\lambda \in \smash{\widebar{{\rm Pr}}}_k(\mfrak{p})$, which together with Theorem \ref{thm:annihilator verma module} implies an inclusion $\smash{\widebar{{\rm Pr}}}_k(\mfrak{p}) \subset \smash{\widebar{{\rm Pr}}}_k^{\smash{\mcal{O}_\mfrak{p}}}\!$. Hence, by using these facts and by the decomposition \eqref{eq:[Pr_k] deconposition} we may write
\begin{align}
  [\widebar{{\rm Pr}}_k^{\mcal{O}}] = \bigcup_{\substack{\Sigma \subset \Pi,\, \mcal{O}=\mcal{O}_{\mfrak{p}_\Sigma}}} [\Lambda_k(\mfrak{p}_\Sigma)]. \label{eq:[Pr_k^O] decomposition}
\end{align}
Let us note that nilpotent orbits of $\mfrak{g}$ are usually parameterize by partitions of $n+1$. We have the following statement.
\medskip

\proposition{\label{prop:[Pr_k^O] decomposition} Let $\mcal{O}_\lambda$ be the nilpotent orbit of $\mfrak{g}$ given by a partition $\lambda \in \mcal{P}_{n+1}$. Then we have
\begin{align*}
 [\widebar{{\rm Pr}}_k^{\smash{\mcal{O}_\lambda}}] = \bigcup_{\mfrak{p} \in [\mfrak{p}_{\lambda^t}]} [\Lambda_k(\mfrak{p})],
\end{align*}
where $\lambda^t$ is the transpose partition of $\lambda$.}

\proof{By results of Section \ref{subsec:Richardson orbits}, the nilpotent orbit $\mcal{O}_\lambda$ is equal to $\smash{\mcal{O}_{\mfrak{p}_{\lambda^t}}}\!$. Besides, we know that $\mcal{O}_{\mfrak{p}_1} = \mcal{O}_{\mfrak{p}_2}$ for standard parabolic subalgebras $\mfrak{p}_1$ and $\mfrak{p}_2$ of $\mfrak{g}$ if and only if $\mfrak{p}_1 \sim \mfrak{p}_2$. The statement then follows easily from \eqref{eq:[Pr_k^O] decomposition}.}

\lemma{\label{lem:condition}Let $y \in W$ and $\eta \in P^\vee$. Then the condition $yt_{-\eta}(\smash{\widehat{\Delta}}(k\Lambda_0)_+) \subset \smash{\widehat{\Delta}}_+^{\rm re}$ is equivalent to
\begin{align*}
  0 \leq (\eta,\alpha) \leq q-1 \ \text{if $y(\alpha) \in \Delta_+$} \qquad \text{and} \qquad 1 \leq (\eta,\alpha) \leq q \ \text{if $y(\alpha) \in \Delta_-$}
\end{align*}
for all $\alpha \in \Delta_+$. Moreover, if $\eta \in P_+^\vee$ satisfies $(\eta,\theta) \leq q-1$, then $yt_{-\eta}(\smash{\widehat{\Delta}}(k\Lambda_0)_+) \subset \smash{\widehat{\Delta}}_+^{\rm re}$ for $y \in W$ if and only if $y^{-1} \in W^{\mfrak{p}_\Sigma}$ with $\Sigma = \{\alpha \in \Pi;\, (\eta,\alpha)=0\}$.}

\proof{The first part of the statement follows from the fact that $\smash{\widehat{\Pi}}(k\Lambda_0) = \{-\theta+q\delta,\alpha_1,\dots,\alpha_n\}$ and that
\begin{align*}
  yt_{-\eta}(\alpha) = y(\alpha) + (\eta,\alpha)\delta \qquad \text{and} \qquad yt_{-\eta}(-\alpha+q\delta) = -y(\alpha) + (q-(\eta,\alpha))\delta
\end{align*}
for $\alpha \in \Delta_+$.

Now, let $y \in W$ and let $\eta \in P_+^\vee$ satisfy $(\eta,\theta) \leq q-1$. Then we have that  $yt_{-\eta}(\smash{\widehat{\Delta}}(k\Lambda_0)_+) \subset \smash{\widehat{\Delta}}_+^{\rm re}$ is equivalent to $(\eta,\alpha) \geq 1$ for all $\alpha \in \Delta_+ \cap y^{-1}(\Delta_-)$ by the first part of the statement. However, this condition is satisfied if and only if $\Delta_+ \cap y^{-1}(\Delta_-) \subset \Delta_+ \setminus \Delta_\Sigma$, where $\Sigma = \{\alpha \in \Pi;\, (\eta,\alpha)=0\}$, which in other words means that $y^{-1} \in W^{\mfrak{p}_\Sigma}$.}

\lemma{\label{lem:W^p symmetries cond}Let $\mfrak{p}$ be a standard parabolic subalgebra of $\mfrak{g}$ and let $j \in \{1,2,\dots,n\}$. Then we have $w_j \in W^\mfrak{p}$ if and only if $\alpha_j \notin \Sigma_\mfrak{p}$. Moreover, the elements of $W_+^\mfrak{p}$ preserve the set $W^\mfrak{p}$ with respect to the left multiplication.}

\proof{Let us assume that $\alpha_j \notin \Sigma_\mfrak{p}$. Let $\alpha \in \Delta_+ \cap w_j(\Delta_-)$. Then there exists $\beta \in \Delta_+$ such that $\alpha= -w_j(\beta)$. Moreover, we have $\beta = \sum_{i=1}^n\! a_i \alpha_i$ with $a_i \in \{0,1\}$ for $i=1,2,\dots,n$.
Further, by Lemma \ref{lem:w_j properties} we may write
\begin{align*}
  \alpha = -w_j(\beta) &= - \sum_{i=1}^{n-j} a_iw_j(\alpha_i) -a_{n+1-j}w_j(\alpha_{n+1-j}) - \sum_{i=n+2-j}^n a_iw_j(\alpha_i) \\
  &= a_{n+1-j}\theta - \sum_{i=1}^{n-j} a_i\alpha_{i+j} - \sum_{i=1}^{j-1} a_{n+1-j+i}\alpha_i,
\end{align*}
which gives us $a_{n+1-j}=1$. As $\alpha_j \notin \Sigma_\mfrak{p}$, we get $\alpha \in \Delta_+\!\setminus \Delta_{\Sigma_\mfrak{p}}$. Therefore, we have $\Delta_+ \cap w_j(\Delta_-) \subset \Delta_+\! \setminus \Delta_{\Sigma_\mfrak{p}}$, which by definition means that $w_j \in W^\mfrak{p}$.

On the other hand, let us assume that $w_j \in W^\mfrak{p}$. Then have $\Delta_+ \cap w_j(\Delta_-) \subset \Delta_+\! \setminus \Delta_{\Sigma_\mfrak{p}}$. It is enough to show that $\alpha_j \in \Delta_+ \cap w_j(\Delta_-)$. But we have $\alpha_j = w_j(-\theta)$ which implies the statement.

Let $w_j \in W_+^\mfrak{p}$ for $j \in \{1,2,\dots,n\}$. Since $w_j(\Sigma_\mfrak{p}) = \Sigma_\mfrak{p}$, we get that $\alpha_j, \alpha_{n+1-j} \notin \Sigma_\mfrak{p}$. Let us assume that $w \in W^\mfrak{p}$. We need to show that $\Delta_+ \cap w_jw(\Delta_-) \subset \Delta_+ \setminus \Delta_{\Sigma_\mfrak{p}}$. If $\alpha \in \Delta_+ \cap w_jw(\Delta_-)$, then there exists $\beta \in \Delta_+$ such that $\alpha = -w_jw(\beta)$. Hence, by using Lemma \ref{lem:w_j properties} we may write
\begin{align*}
  -w(\beta) = w_{n+1-j}(\alpha) =  \sum_{i=1}^{j-1} a_i\alpha_{n+1-j+i} - a_j\theta + \sum_{i=j+1}^n a_i\alpha_{i-j},
\end{align*}
where $\alpha = \sum_{i=1}^n\! a_i\alpha_i$ with $a_i \in \{0,1\}$ for $i=1,2,\dots,n$. If $a_j = 1$ or $a_{n+1-j}=1$, then we have $\alpha \in \Delta_+ \setminus \Delta_{\Sigma_\mfrak{p}}$ since $\alpha_j, \alpha_{n+1-j} \notin \Sigma_\mfrak{p}$. On the other hand, if $a_j=a_{n+1-j}=0$, then $-w(\beta) \in \Delta_+$, which implies $-w(\beta) \in \Delta_+ \setminus \Delta_{\Sigma_\mfrak{p}}$ since $w \in W^\mfrak{p}$. As we have $w_j(\Sigma_\mfrak{p})= \Sigma_\mfrak{p}$ and $\alpha=-w_jw(\beta)$, we get $\alpha \in \Delta_+ \setminus \Delta_{\Sigma_\mfrak{p}}$. This finishes the proof.}

By the previous lemma we may introduce an equivalence relation on the set $W^\mfrak{p}$ in such a way that the equivalence classes coincide with the orbits of $W_+^\mfrak{p}$ on $W^\mfrak{p}$.
\medskip

\lemma{\label{lem:Pr_k Lambda_k(p) intersec}Let $\mfrak{p}$ be a standard parabolic subalgebra of $\mfrak{g}$. Then $\smash{\widebar{{\rm Pr}}}_k \cap w^{-1} \cdot \Lambda_k(\mfrak{p}) \neq \emptyset$ if and only if $w \in W^\mfrak{p}$. Moreover, if $w \in W^\mfrak{p}$ then $w^{-1} \cdot \Lambda_k(\mfrak{p}) \subset \smash{\widebar{{\rm Pr}}}_k$.}

\proof{Since we have $\Lambda_k(\mfrak{p}) \subset \smash{\widebar{{\rm Pr}}}_k$, by Lemma \ref{lem:W^p condition} we get that $w^{-1} \cdot \lambda \in \smash{\widebar{{\rm Pr}}}_k$ for $\lambda \in \Lambda_k(\mfrak{p})$ if and only if $w \in W^\mfrak{p}$, which implies the required statement.}

\lemma{\label{lem:Lambda_k(p)} Let $\mfrak{p}_1, \mfrak{p}_2$ be standard parabolic subalgebras of $\mfrak{g}$. Further, let $y_1,y_2 \in W$ satisfy $y_1 \cdot \Lambda_k(\mfrak{p}_1) \subset \smash{\widebar{{\rm Pr}}}_k$, $y_2 \cdot \Lambda_k(\mfrak{p}_2) \subset \smash{\widebar{{\rm Pr}}}_k$ and $y_1 \neq y_2$. Then the following statements
\begin{enumerate}[topsep=3pt,itemsep=0pt]
  \item[i)] $y_1 \cdot \Lambda_k(\mfrak{p}_1) \cap y_2 \cdot \Lambda_k(\mfrak{p}_2) \neq \emptyset$;
  \item[ii)] $y_1 \cdot \Lambda_k(\mfrak{p}_1) = y_2 \cdot \Lambda_k(\mfrak{p}_2)$;
  \item[iii)] there exists $j \in \{1,2,\dots,n\}$ such that $y_2=y_1w_j$ and $\Sigma_{\mfrak{p}_1} = w_j(\Sigma_{\mfrak{p}_2})$
\end{enumerate}
are equivalent. In addition, if $\Lambda_k(\mfrak{p}_1) \cap \Lambda_k(\mfrak{p}_2) \neq \emptyset$ then $\mfrak{p}_1 = \mfrak{p}_2$.}

\proof{We will prove the implications (iii) $\Longrightarrow$ (ii) $\Longrightarrow$ (i) $\Longrightarrow$ (iii).

Let us assume that there exists $j \in \{1,2,\dots,n\}$ such that $y_2=y_1w_j$ and $\Sigma_{\mfrak{p}_1} = w_j(\Sigma_{\mfrak{p}_2})$. By Lemma \ref{lem:W^p symmetries cond} we have $w_j^{-1} \in W^{\mfrak{p}_2}$, which means that $w_j \cdot \Lambda_k(\mfrak{p}_2) \subset \smash{\widebar{{\rm Pr}}}_k$. Let $\eta \in P_+^\vee$ satisfy $(\eta,\theta) \leq q-1$ and $\Pi^\eta = \Sigma_{\mfrak{p}_1}\!$. By definition of $\Lambda_k(\mfrak{p}_1)$ we have $\smash{\widebar{{\rm Pr}}_{k,t_{-\eta}}} \subset \Lambda_k(\mfrak{p}_1)$. Further, we set $\eta'=w_j^{-1}(\eta-q\omega_j)$. By using Lemma \ref{lem:w_j properties} we may write
\begin{align*}
  \eta' = w_{n+1-j}(\eta-q\omega_j) = \sum_{i=1}^{n-j} a_{i+j}\omega_i + \bigg(q -\sum_{i=1}^n a_i\!\bigg)\omega_{n+1-j} + \sum_{i=n+2-j}^n a_{i+j-n-1}\omega_i,
\end{align*}
where $\eta = \sum_{i=1}^n\! a_i \omega_i$. Since $\eta \in P_+^\vee$ and $(\eta,\theta)\leq q-1$ by the assumption, we have also $\eta' \in P_+^\vee$. Besides, as $w_j(-\theta)=\alpha_j$ and $\Sigma_{\mfrak{p}_1} = w_j(\Sigma_{\mfrak{p}_2})$, we get $\alpha_j \notin \Sigma_{\mfrak{p}_1}\!$ which gives us $a_j \neq 0$. Thus, we have $(\eta',\theta) = q - a_j \leq q-1$. Further, we see that $\smash{\Pi^{\eta'}}=\Sigma_{\mfrak{p}_2}$, which implies $\smash{\widebar{{\rm Pr}}}_{k,t_{\smash{-\eta'}}} \subset \Lambda_k(\mfrak{p}_2)$. By Lemma \ref{lem:condition} we have $t_{-\eta}(\smash{\widehat{\Delta}}(k\Lambda_0)_+) \subset \smash{\widehat{\Delta}}_+^{\rm re}$ and $w_jt_{-\smash{\eta'}} (\smash{\widehat{\Delta}}(k\Lambda_0)_+) \subset \smash{\widehat{\Delta}}_+^{\rm re}$, hence we
may write
\begin{align*}
  \widebar{{\rm Pr}}_{k,t_{-\eta}} = \widebar{{\rm Pr}}_{k,w_jt_{\smash{-\eta'}}} = w_j \cdot \widebar{{\rm Pr}}_{k,t_{\smash{-\eta'}}} \subset w_j \cdot \Lambda_k(\mfrak{p}_2),
\end{align*}
where we used \eqref{eq:Pr equality} in the first equality, and thus $\Lambda_k(\mfrak{p}_1) \subset w_j \cdot \Lambda_k(\mfrak{p}_2)$. By the same argument we obtain the opposite inclusion. Therefore, we have $\Lambda_k(\mfrak{p}_1) = w_j \cdot \Lambda_k(\mfrak{p}_2)$, which immediately implies $y_1 \cdot \Lambda_k(\mfrak{p}_1) = y_2 \cdot \Lambda_k(\mfrak{p}_2)$. This proves the first implication (iii) $\Longrightarrow$ (ii).

Let us assume that $y_1 \cdot \Lambda_k(\mfrak{p}_1) \cap y_2 \cdot \Lambda_k(\mfrak{p}_2) \neq \emptyset$, which is the same as $\Lambda_k(\mfrak{p}_1) \cap w \cdot \Lambda_k(\mfrak{p}_2) \neq \emptyset$ with $w=y_1^{-1}y_2$. As we have $\smash{\widebar{{\rm Pr}}}_k \cap w \cdot \Lambda_k(\mfrak{p}_2) \neq \emptyset$, by Lemma \ref{lem:Pr_k Lambda_k(p) intersec} we obtain that $w^{-1} \in W^{\mfrak{p}_2}$. Hence, there exist $\eta, \eta' \in P_+^\vee$ satisfying $\Pi^\eta=\Sigma_{\mfrak{p}_1}$, $\smash{\Pi^{\eta'}}=\Sigma_{\mfrak{p}_2}$, $(\eta,\theta)\leq q-1$, $(\eta',\theta)\leq q-1$ such that $\smash{\widebar{{\rm Pr}}_{k,t_{-\eta}}}\!=\smash{\widebar{{\rm Pr}}_{kwt_{\smash{-\eta'}}}}\!$, which gives us $wt_{-\eta'} = t_{-\eta}t_{q\omega_j}w_j$ for some $j \in \{1,2,\dots,n\}$. Therefore, we have $w=w_j$ and $\eta'=w_j^{-1}(\eta-q\omega_j)$. Further, by using Lemma \ref{lem:w_j properties} we may write
\begin{align*}
  \eta' = w_{n+1-j}(\eta-q\omega_j) = \sum_{i=1}^{n-j} a_{i+j}\omega_i + \bigg(q -\sum_{i=1}^n a_i\!\bigg)\omega_{n+1-j} + \sum_{i=n+2-j}^n a_{i+j-n-1}\omega_i,
\end{align*}
where $\eta = \sum_{i=1}^n\! a_i \omega_i$. As we have $(\eta',\theta)=q-a_j \leq q-1$, we obtain that $a_j\neq 0$. Then it follows immediately that $\Sigma_{\mfrak{p}_1}= \Pi^\eta = w_j(\smash{\Pi^{\eta'}}) = w_j(\Sigma_{\mfrak{p}_2})$. This proves the last implication (i) $\Longrightarrow$ (iii). The implication (ii) $\Longrightarrow$ (i) is obvious.}

\lemma{\label{lem:[Lambda_k(p)]}Let $\mfrak{p}_1, \mfrak{p}_2$ be standard parabolic subalgebras of $\mfrak{g}$. Then the statements
\begin{enumerate}[topsep=3pt,itemsep=0pt]
  \item[i)] $[\Lambda_k(\mfrak{p}_1)] \cap [\Lambda_k(\mfrak{p}_2)] \neq \emptyset$;
  \item[ii)] $[\Lambda_k(\mfrak{p}_1)] = [\Lambda_k(\mfrak{p}_2)]$;
  \item[iii)] $\mfrak{p}_1 \sim_+ \mfrak{p}_2$
\end{enumerate}
are equivalent.}

\proof{We may assume that the parabolic subalgebras $\mfrak{p}_1$ and $\mfrak{p}_2$ are distinct, otherwise the statement is trivial. We will prove the implications (iii) $\Longrightarrow$ (ii) $\Longrightarrow$ (i) $\Longrightarrow$ (iii).

If we have $\mfrak{p}_1 \sim_+ \mfrak{p}_2$, then there exists $j\in \{1,2,\dots,n\}$ such that $\Sigma_{\mfrak{p}_1} = w_j(\Sigma_{\mfrak{p}_2})$. Hence, by Lemma \ref{lem:Lambda_k(p)} we have $\Lambda_k(\mfrak{p}_1) = w_j\cdot \Lambda_k(\mfrak{p}_2)$, which gives us $[\Lambda_k(\mfrak{p}_1)] = [\Lambda_k(\mfrak{p}_2)]$. This proves the first implication (iii) $\Longrightarrow$ (ii).

Further, let us assume that $[\Lambda_k(\mfrak{p}_1)] \cap [\Lambda_k(\mfrak{p}_2)] \neq \emptyset$. Then there exist $y_1, y_2 \in W$ satisfying $y_1 \cdot \Lambda_k(\mfrak{p}_1) \subset \smash{\widebar{{\rm Pr}}}_k$, $y_2 \cdot \Lambda_k(\mfrak{p}_2) \subset \smash{\widebar{{\rm Pr}}}_k$ and $y_1 \cdot \Lambda_k(\mfrak{p}_1) \cap y_2 \cdot \Lambda_k(\mfrak{p}_2) \neq \emptyset$. By Lemma \ref{lem:Lambda_k(p)} we get that $\mfrak{p}_1 \sim_+ \mfrak{p}_2$. This proves the implication (i) $\Longrightarrow$ (iii). We are done since the implication (ii) $\Longrightarrow$ (i) is obvious.}

The following statement is the main theorem of this section.
\medskip

\theorem{\label{thm:Pr_k parametrization} Let $\mcal{O}_\lambda$ be the nilpotent orbit of $\mfrak{g}$ given by a partition $\lambda \in \mcal{P}_{n+1}$. Then we have
\begin{align*}
\widebar{{\rm Pr}}_k^{\smash{\mcal{O}_\lambda}} = \bigsqcup_{[\mfrak{p}]_+ \in [\mfrak{p}_{\lambda^t}]/\sim_+} \, \bigsqcup_{[w] \in W_+^\mfrak{p}\backslash W^\mfrak{p}} w^{-1} \cdot \Lambda_k(\mfrak{p}),
\end{align*}
where $W_+^\mfrak{p}\backslash W^\mfrak{p}$ denotes the set of orbits of $W_+^\mfrak{p}$ on $W^\mfrak{p}$.}

\proof{By using Proposition \ref{prop:[Pr_k^O] decomposition} and Lemma \ref{lem:[Lambda_k(p)]} we may write
\begin{align*}
  [\widebar{{\rm Pr}}_k^{\smash{\mcal{O}_\lambda}}] = \bigcup_{\mfrak{p} \in [\mfrak{p}_{\lambda^t}]} [\Lambda_k(\mfrak{p})] = \bigsqcup_{[\mfrak{p}]_+ \in [\mfrak{p}_{\lambda^t}]/\sim_+} [\Lambda_k(\mfrak{p})].
\end{align*}
Further, by Lemma \ref{lem:Pr_k Lambda_k(p) intersec} we have
\begin{align*}
  \widebar{{\rm Pr}}_k^{\smash{\mcal{O}_\lambda}} = \bigsqcup_{[\mfrak{p}]_+ \in [\mfrak{p}_{\lambda^t}]/\sim_+}\, \bigcup_{w \in W^\mfrak{p}} w^{-1} \cdot \Lambda_k(\mfrak{p}).
\end{align*}
As by Lemma \ref{lem:Lambda_k(p)} we have $w_1^{-1} \cdot \Lambda_k(\mfrak{p}) = w_2^{-1} \cdot \Lambda_k(\mfrak{p})$ for $w_1,w_2 \in W^\mfrak{p}$ if only if there exists $w \in W^\mfrak{p}_+$ such that $w_2 = ww_1$, we immediately get the required statement.}

On two concrete examples we show how to determine the set $\smash{\widebar{{\rm Pr}}_k^\mcal{O}}$ for a nilpotent orbit $\mcal{O}$ of $\mfrak{g}$ by using Theorem \ref{thm:Pr_k parametrization}.
\smallskip

\noindent
{\bf Example.} Let us start with $\mfrak{g}=\mfrak{sl}_6$ and $\mcal{O}=\mcal{O}_{[4,2]}$. The nilpotent orbit $\mcal{O}$ is the Richardson orbit attached to the standard parabolic subalgebra $\mfrak{p}_{[2,2,1,1]} = \dynkin[x/.style={thin}]{A}{oxoxx}$ with the corresponding equivalence class given by
\begin{align*}
  [\mfrak{p}_{[2,2,1,1]}] = \{\dynkin[x/.style={thin}]{A}{oxoxx}, \dynkin[x/.style={thin}]{A}{oxxox}, \dynkin[x/.style={thin}]{A}{oxxxo}, \dynkin[x/.style={thin}]{A}{xoxox}, \dynkin[x/.style={thin}]{A}{xoxxo}, \dynkin[x/.style={thin}]{A}{xxoxo}\}.
\end{align*}
Further, we have
\begin{align*}
  \dynkin[x/.style={thin}]{A}{oxoxx} \xrightarrow{w_1} \dynkin[x/.style={thin}]{A}{xoxox} \xrightarrow{w_1} \dynkin[x/.style={thin}]{A}{xxoxo} \xrightarrow{w_2} \dynkin[x/.style={thin}]{A}{oxxxo} \xrightarrow{w_2} \dynkin[x/.style={thin}]{A}{oxoxx}
\end{align*}
and
\begin{align*}
  \dynkin[x/.style={thin}]{A}{oxxox} \xrightarrow{w_1} \dynkin[x/.style={thin}]{A}{xoxxo} \xrightarrow{w_2} \dynkin[x/.style={thin}]{A}{oxxox}
\end{align*}
which gives us that $[\mfrak{p}_{[2,2,1,1]}]$ decomposes into two equivalence classes with respect to the equivalence relation $\sim_+$, i.e.\ we obtain
\begin{align*}
  [\mfrak{p}_1]_+ &= \{\dynkin[x/.style={thin}]{A}{oxoxx}, \dynkin[x/.style={thin}]{A}{xoxox}, \dynkin[x/.style={thin}]{A}{xxoxo}, \dynkin[x/.style={thin}]{A}{oxxxo}\}, \\
  [\mfrak{p}_2]_+ &= \{\dynkin[x/.style={thin}]{A}{oxxox}, \dynkin[x/.style={thin}]{A}{xoxxo}\},
\end{align*}
where $\mfrak{p}_1 = \dynkin[x/.style={thin}]{A}{oxoxx}$ and $\mfrak{p}_2 = \dynkin[x/.style={thin}]{A}{oxxox}$. Besides, we see that
\begin{align*}
  W_+^{\mfrak{p}_1} = \{e\} \qquad \text{and} \qquad W_+^{\mfrak{p}_2} = \{e,w_3\}.
\end{align*}
Therefore, by Theorem \ref{thm:Pr_k parametrization} we may write
\begin{align*}
  \widebar{{\rm Pr}}_k^\mcal{O} = \bigsqcup_{w \in W^{\mfrak{p}_1}} w^{-1} \cdot \Lambda_k(\mfrak{p}_1) \sqcup \bigsqcup_{[w] \in W_+^{\mfrak{p}_2}\backslash W^{\mfrak{p}_2}} w^{-1} \cdot \Lambda_k(\mfrak{p}_2),
\end{align*}
where
\begin{align*}
  \Lambda_k(\mfrak{p}_1) & = \bigg\{\lambda- {p\over q}(a\omega_2 + b\omega_4 + c\omega_5);\, \lambda \in \widebar{{\rm Pr}}_{k,\Z},\, a,b,c \in \N,\, a+b+c \leq q-1\bigg\}, \\
  \Lambda_k(\mfrak{p}_2) & = \bigg\{\lambda- {p\over q}(a\omega_2 + b\omega_3 + c\omega_5);\, \lambda \in \widebar{{\rm Pr}}_{k,\Z},\, a,b,c \in \N,\, a+b+c \leq q-1\bigg\}.
\end{align*}
In addition, we can easily see that $\smash{\widebar{{\rm Pr}}_k^\mcal{O}}$ is non-empty if $q \geq 4$. The last step is to determine the sets $W^{\mfrak{p}_1}$ and $W^{\mfrak{p}_2}$, which we usually describe by the corresponding Hasse diagrams. In fact, we do not need to know the structure of the Bruhat ordering on $W^{\mfrak{p}_1}$ and $W^{\mfrak{p}_2}$.
\smallskip

\noindent
{\bf Example.} Let now $\mfrak{g}=\mfrak{sl}_3$ and $\mcal{O}=\mcal{O}_{[3]}=\mcal{O}_{\rm prin}$. The nilpotent orbit $\mcal{O}$ is attached to the standard parabolic subalgebra $\mfrak{p}_{[1,1,1]} = \dynkin[x/.style={thin}]{A}{xx} = \mfrak{b}$ with the corresponding equivalence class given by
\begin{align*}
  [\mfrak{p}_{[1,1,1]}] = \{\dynkin[x/.style={thin}]{A}{xx}\}.
\end{align*}
Hence, we have
\begin{align*}
  \dynkin[x/.style={thin}]{A}{xx} \xrightarrow{w_1} \dynkin[x/.style={thin}]{A}{xx},
\end{align*}
which gives us $[\mfrak{b}]_+ = \{\dynkin[x/.style={thin}]{A}{xx}\}$ and $W_+^\mfrak{b}=\{e,w_1,w_2\} =W_+$. Therefore, by Theorem \ref{thm:Pr_k parametrization} we may write
\begin{align*}
  \widebar{{\rm Pr}}_k^\mcal{O} = \bigsqcup_{[w] \in W_+^\mfrak{b} \backslash W^\mfrak{b}} w^{-1} \cdot \Lambda_k(\mfrak{b}),
\end{align*}
where
\begin{align*}
  \Lambda_k(\mfrak{b}) & = \bigg\{\lambda- {p\over q}(a\omega_1 + b\omega_2);\, \lambda \in \widebar{{\rm Pr}}_{k,\Z},\, a,b \in \N,\, a+b \leq q-1\bigg\}.
\end{align*}
Moreover, we can see that $\smash{\widebar{{\rm Pr}}_k^\mcal{O}}$ is non-empty if $q \geq 3$. The elements of $W^\mfrak{b}$ we can read off directly from the corresponding Hasse diagram
\medskip

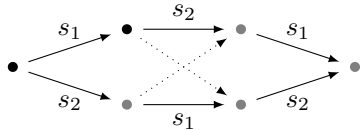
\hfil
\begin{minipage}[h]{0.4\textwidth}
\begin{tikzpicture}
[>=latex, yscale=1, xscale=1.5,point/.style={circle,draw=white,fill=black,line width=1.5mm, inner sep=1mm}, point2/.style={circle,draw=white,fill=gray,line width=1.5mm, inner sep=1mm}]
  \node (A) at (0,0) [point] {};
  \node (B) at (1,0.5) [point] {};
  \node (C) at (1,-0.5) [point2] {};
  \node (D) at (2,0.5) [point2] {};
  \node (E) at (2,-0.5) [point2] {};
  \node (F) at (3,0) [point2] {};
  \draw [thin, ->] (A) to node[above] {$s_1$} (B);
  \draw [thin, ->] (A) to node[below] {$s_2$} (C);
  \draw [thin, ->] (B) to node[above] {$s_2$} (D);
  \draw [thin, ->] (C) to node[below] {$s_1$} (E);
  \draw [thin, ->] (D) to node[above] {$s_1$} (F);
  \draw [thin, ->] (E) to node[below] {$s_2$} (F);
  \draw [thin,dotted, ->] (B) -- (E);
  \draw [thin,dotted, ->] (C) -- (D);
\end{tikzpicture}
\end{minipage}
\hfil

\medskip

\noindent
as follows. The vertices of the graph represent the elements of $W^\mfrak{b}$. The least element is the unit $e$ of $W$ and if there is an arrow
\medskip

\hfil
\begin{minipage}[h]{0.2\textwidth}
\begin{tikzpicture}
[>=latex, yscale=1, xscale=1.5,point/.style={circle,draw=white,fill=black,line width=1.5mm, inner sep=1mm}, point2/.style={circle,draw=white,fill=gray,line width=1.5mm, inner sep=1mm}]
  \node (A) at (0,0) [] {$w$};
  \node (B) at (1,0) [] {$w\smash{'}$};
  \draw [thin, ->] (A) to node[above] {$s_i$} (B);
\end{tikzpicture}
\end{minipage}
\hfil

\medskip

\noindent
for $w, w' \in W$, then we have $w' = ws_i$. Hence, we easily get $W^\mfrak{b}=\{e, s_1, s_2, s_1s_2, s_2s_1, s_1s_2s_1\}$. In addition, the black vertices of the graph stand for representatives of the equivalence classes in $W_+^\mfrak{b} \backslash W^\mfrak{b}$.

%%%%%%%%%%%%%%%%%%%%%%%%%%%%%%%%%%%%%%%%%%%%%%%%%%%%%%%%%%%%%%%%%%%%%%%%%%%%%%%%%%%%%%%%%%

\subsubsection{Minimal and subregular nilpotent orbit}

Let us consider the simple Lie algebra $\mfrak{g}=\mfrak{sl}_{n+1}$ with $n \geq 2$ and let $k\in \Q$ be an admissible number for $\mfrak{g}$ with denominator $q \in \N$.
\medskip

The minimal nilpotent orbit $\mcal{O}_{\rm min}$ of $\mfrak{g}$ corresponds to the partition $\lambda=[2,1^{n-1}]$ with the transpose partition $\lambda^t=[n,1]$. Therefore, the nilpotent orbit $\mcal{O}_{\rm min}$ is attached to the standard parabolic subalgebra

\begin{align*}
  \mfrak{p}_{\lambda^t} = \begin{dynkinDiagram}[edge length=6mm, x/.style={thin}]{A}{oo..ox}
  \dynkinBrace[n-1]{1}{3}
  \end{dynkinDiagram}
\end{align*}

with the corresponding equivalence class given by
\begin{align*}
  [\mfrak{p}_{\lambda^t}] = \{\mfrak{p}_{\alpha_1}^{\rm max}, \mfrak{p}_{\alpha_n}^{\rm max}\},
\end{align*}
where $\mfrak{p}_{\alpha}^{\rm max}$ for $\alpha \in \Pi$ is the standard parabolic subalgebra of $\mfrak{g}$ associated to the subset $\Sigma = \Pi \setminus \{\alpha\}$ of $\Pi$. Besides, we see that $[\mfrak{p}_{\lambda^t}]$ decomposes into one equivalence class with respect to the equivalence relation $\sim_+$, i.e.\ we have
\begin{align*}
  [\mfrak{p}_{\alpha_1}^{\rm max}]_+ = [\mfrak{p}_{\lambda^t}]
\end{align*}
together with $W_+^{\smash{\mfrak{p}_{\alpha_1}^{\rm max}}} = \{e\}$. By using Theorem \ref{thm:Pr_k parametrization} we get the following statement.
\medskip

\theorem{Let $k\in\Q$ be an admissible number for $\mfrak{g}$ with denominator $q \in \N$. Then we have
\begin{align*}
  \widebar{{\rm Pr}}_k^{\smash{\mcal{O}_{\rm min}}} = \bigsqcup_{w \in W^{\smash{\mfrak{p}_{\alpha_1}^{\rm max}}}} w^{-1} \cdot \Lambda_k(\mfrak{p}_{\alpha_1}^{\rm max}),
\end{align*}
where
\begin{align*}
  \Lambda_k(\mfrak{p}_{\alpha_1}^{\rm max}) = \bigg\{\lambda - {p \over q}\,a\omega_1;\, \lambda \in \widebar{{\rm Pr}}_{k,\Z},\ a \in \N,\, a \leq q-1\bigg\}
\end{align*}
and
\begin{align*}
  W^{\smash{\mfrak{p}_{\alpha_1}^{\rm max}}} = \{e, s_1, s_1s_2,\dots, s_1s_2\cdots s_n\}.
\end{align*}
Moreover, the set $\widebar{{\rm Pr}}_k^{\smash{\mcal{O}_{\rm min}}}$ is non-empty if $q\geq 2$.}

The subregular nilpotent orbit $\mcal{O}_{\rm subreg}$ of $\mfrak{g}$ corresponds to the partition $\lambda=[n,1]$ with the transpose partition $\lambda^t=[2,1^{n-1}]$. Hence, the subregular orbit $\mcal{O}_{\rm subreg}$ is attached to the standard parabolic subalgebra
\begin{align*}
  \mfrak{p}_{\lambda^t} = \begin{dynkinDiagram}[edge length=6mm, x/.style={thin}]{A}{ox..xx}
  \dynkinBrace[n-1]{2}{4}
  \end{dynkinDiagram}
\end{align*}
with the corresponding equivalence class given by
\begin{align*}
  [\mfrak{p}_{\lambda^t}] = \{\mfrak{p}_{\alpha_1}^{\rm min}, \mfrak{p}_{\alpha_2}^{\rm min}, \dots, \mfrak{p}_{\alpha_n}^{\rm min}\},
\end{align*}
where $\mfrak{p}_{\alpha}^{\rm min}$ for $\alpha \in \Pi$ is the standard parabolic subalgebra of $\mfrak{g}$ associated to the subset $\Sigma = \{\alpha\}$ of $\Pi$. In addition, we see that $[\mfrak{p}_{\lambda^t}]$ decomposes into one equivalence class with respect to the equivalence relation $\sim_+$, i.e.\ we have
\begin{align*}
  [\mfrak{p}_{\alpha_n}^{\rm min}]_+ = [\mfrak{p}_{\lambda^t}]
\end{align*}
together with $W_+^{\smash{\mfrak{p}_{\alpha_n}^{\rm min}}} = \{e\}$. By using Theorem \ref{thm:Pr_k parametrization} we get the following statement.
\medskip

\theorem{\label{thm:subregular orbit}Let $k\in\Q$ be an admissible number for $\mfrak{g}$ with denominator $q \in \N$. Then we have
\begin{align*}
  \widebar{{\rm Pr}}_k^{\smash{\mcal{O}_{\rm subreg}}} = \bigsqcup_{w \in W^{\smash{\mfrak{p}_{\alpha_n}^{\rm min}}}} w^{-1} \cdot \Lambda_k(\mfrak{p}_{\alpha_n}^{\rm min}),
\end{align*}
where
\begin{align*}
  \Lambda_k(\mfrak{p}_{\alpha_n}^{\rm min}) = \bigg\{\lambda - {p \over q} \sum_{i=1}^{n-1} a_i\omega_i;\, \lambda \in \widebar{{\rm Pr}}_{k,\Z},\ a_1,a_2,\dots, a_{n-1} \in \N,\, \sum_{i=1}^{n-1} a_i \leq q-1\bigg\}.
\end{align*}
and
\begin{align*}
  W^{\smash{\mfrak{p}_{\alpha_n}^{\rm min}}} = \{w\in W;\, \ell(s_nw)=\ell(w)+1\}.
\end{align*}
Moreover, the set $\widebar{{\rm Pr}}_k^{\smash{\mcal{O}_{\rm subreg}}}$ is non-empty if $q\geq n$.}

\vspace{-2mm}

%%%%%%%%%%%%%%%%%%%%%%%%%%%%%%%%%%%%%%%%%%%%%%%%%%%%%%%%%%%%%%%%%%%%%%%%%%%%%%%%%%%%%%%%%%

\subsubsection{Zero and principal nilpotent orbit}

Let us consider the simple Lie algebra $\mfrak{g}=\mfrak{sl}_{n+1}$ with $n \geq 1$ and let $k\in \Q$ be an admissible number for $\mfrak{g}$ with denominator $q \in \N$.
\medskip

The zero nilpotent orbit $\mcal{O}_{\rm zero}$ of $\mfrak{g}$ corresponds to the partition $\lambda=[1^{n+1}]$ with the transpose partition $\lambda^t=[n+1]$. Therefore, the nilpotent orbit $\mcal{O}_{\rm zero}$ is attached to the standard parabolic subalgebra
\begin{align*}
  \mfrak{p}_{\lambda^t} = \begin{dynkinDiagram}[edge length=6mm, x/.style={thin}]{A}{oo..oo}
  \dynkinBrace[n]{1}{4}
  \end{dynkinDiagram}
\end{align*}
with the corresponding equivalence class $[\mfrak{p}_{\lambda^t}] = \{\mfrak{g}\}$. Hence, we have $[\mfrak{g}]_+ = [\mfrak{p}_{\lambda^t}]$, $W^\mfrak{g}=\{e\}$ and $W_+^\mfrak{g} = \{e\}$. By using Theorem \ref{thm:Pr_k parametrization} we get the following statement, cf.\ \cite{Malikov-Frenkel1999}.
\medskip

\theorem{\label{thm:zero orbit}Let $k\in\Q$ be an admissible number for $\mfrak{g}$ with denominator $q \in \N$. Then we have
\begin{align*}
  \widebar{{\rm Pr}}_k^{\smash{\mcal{O}_{\rm zero}}} = \Lambda_k(\mfrak{g}) = \widebar{{\rm Pr}}_{k,\Z}.
\end{align*}
Moreover, the set $\widebar{{\rm Pr}}_k^{\smash{\mcal{O}_{\rm zero}}}$ is always non-empty.}

The principal nilpotent orbit $\mcal{O}_{\rm prin}$ of $\mfrak{g}$ corresponds to the partition $\lambda=[n+1]$ with the transpose partition $\lambda^t=[1^{n+1}]$. Therefore, the nilpotent orbit $\mcal{O}_{\rm prin}$ is attached to the standard parabolic subalgebra
\begin{align*}
  \mfrak{p}_{\lambda^t} = \begin{dynkinDiagram}[edge length=6mm, x/.style={thin}]{A}{xx..xx}
  \dynkinBrace[n]{1}{4}
  \end{dynkinDiagram}
\end{align*}
with the corresponding equivalence class $[\mfrak{p}_{\lambda^t}] = \{\mfrak{b}\}$. Hence, we have $[\mfrak{b}]_+ = [\mfrak{p}_{\lambda^t}]$, $W^\mfrak{b}=W$ and $W_+^\mfrak{b}=W_+$. Again, by using Theorem \ref{thm:Pr_k parametrization} we get the following statement.
\medskip

\theorem{\label{thm-prin} Let $k\in\Q$ be an admissible number for $\mfrak{g}$ with denominator $q \in \N$. Then we have
\begin{align*}
  \widebar{{\rm Pr}}_k^{\smash{\mcal{O}_{\rm prin}}} = \bigsqcup_{[w] \in W_+\! \backslash W} w^{-1} \cdot \Lambda_k(\mfrak{b}),
\end{align*}
where
\begin{align*}
  \Lambda_k(\mfrak{b}) = \bigg\{\lambda - {p \over q} \sum_{i=1}^n a_i\omega_i;\, \lambda \in \widebar{{\rm Pr}}_{k,\Z},\ a_1,a_2,\dots,a_n \in \N,\, \sum_{i=1}^n a_i \leq q-1\bigg\}.
\end{align*}
Moreover, the set $\widebar{{\rm Pr}}_k^{\smash{\mcal{O}_{\rm prin}}}$ is non-empty if $q\geq n+1$.}

\vspace{-2mm}

%%%%%%%%%%%%%%%%%%%%%%%%%%%%%%%%%%%%%%%%%%%%%%%%%%%%%%%%%%%%%%%%%%%%%%%%%%%%%%%%%%%%%%%%%%
%%%%%%%%%%%%%%%%%%%%%%%%%%%%%%%%%%%%%%%%%%%%%%%%%%%%%%%%%%%%%%%%%%%%%%%%%%%%%%%%%%%%%%%%%%

\subsection{Examples}

In this section we give some examples of admissible highest weights. Let us consider the simple Lie algebra $\mfrak{g} = \mfrak{sl}_{n+1}$ with $n \geq 1$ and let $k\in \Q$ be an admissible number for $\mfrak{g}$ with denominator $q\in \N$, i.e.\ we have
\begin{align*}
  k+n+1 = {p \over q}\ \text{with}\ p,q \in \N,\, (p,q)=1,\, p \geq n+1.
\end{align*}
%Based on the discussion in Section \ref{sec:associated varieties},
We have a decomposition
\begin{align*}
  \widebar{{\rm Pr}}_k = \bigsqcup_{\mcal{O} \subset \mcal{N}(\mfrak{g})} \widebar{{\rm Pr}}_k^\mcal{O},
\end{align*}
where the set $\smash{\widebar{{\rm Pr}}_k^\mcal{O}}$ of admissible weights of level $k$ belonging to a nilpotent orbit $\mcal{O}$ of $\mfrak{g}$ may be empty set. In fact, we have that $\smash{\widebar{{\rm Pr}}_k^\mcal{O}} \neq \emptyset$ if and only if $\mcal{O} \subset \widebar{{\mcal{O}}_q}$, where the nilpotent orbit $\mcal{O}_q$ of $\mfrak{g}$ is determined by Theorem \ref{thm:nilpotent orbits}.

We focus on two possible directions for which a complete description of admissible weights is relatively easy. Either we can consider the Lie algebra $\mfrak{g}$ of low rank, i.e.\ $n=1,2,3$, or take into consideration the level $k$ with small denominator, i.e.\ $q=1,2$.

\begin{figure}[ht]
\centering
\subcaptionbox{$\mfrak{g}=\mfrak{sl}_2$\label{fig:sl2}}[0.25\textwidth]
{\begin{tikzpicture}
[yscale=1.2, xscale=1]
  \node (A) at (0,0)  {$\mcal{O}_{[1^2]}$};
  \node (B) at (0,1)  {$\mcal{O}_{[2]}$};
  \node at (1,0) {$0$};
  \node at (1,1) {$2$};
  \node at (2,0) {\dynkin[x/.style={thin}]{A}{o}};
  \node at (2,1) {\dynkin[x/.style={thin}]{A}{x}};
 \draw [thin, -] (A) -- (B);
\end{tikzpicture}}
\hfill
\subcaptionbox{$\mfrak{g}=\mfrak{sl}_3$\label{fig:sl3}}[0.25\textwidth]
{\begin{tikzpicture}
[yscale=1.2, xscale=1]
  \node (A) at (0,0)  {$\mcal{O}_{[1^3]}$};
  \node (B) at (0,1)  {$\mcal{O}_{[2,1]}$};
  \node (C) at (0,2)  {$\mcal{O}_{[3]}$};
  \node at (1,0) {$0$};
  \node at (1,1) {$4$};
  \node at (1,2) {$6$};
  \node at (2,0) {\dynkin[x/.style={thin}]{A}{oo}};
  \node at (2,0.9) {\dynkin[x/.style={thin}]{A}{ox}};
  \node at (2,1.1) {\dynkin[x/.style={thin}]{A}{xo}};
  \node at (2,2) {\dynkin[x/.style={thin}]{A}{xx}};
  \node at (0,4.3) {$\mcal{O}$};
  \node at (1,4.3) {$\dim$};
  \node at (2,4.3) {$\mfrak{p}$};
  \draw [thin, -] (A) -- (B);
  \draw [thin, -] (B) -- (C);
\end{tikzpicture}}
\hfill
\subcaptionbox{$\mfrak{g}=\mfrak{sl}_4$\label{fig:sl4}}[0.25\textwidth]
{\begin{tikzpicture}
[yscale=1.2, xscale=1]
  \node (A) at (0,0)  {$\mcal{O}_{[1^4]}$};
  \node (B) at (0,1)  {$\mcal{O}_{[2,1^2]}$};
  \node (C) at (0,2)  {$\mcal{O}_{[2^2]}$};
  \node (D) at (0,3)  {$\mcal{O}_{[3,1]}$};
  \node (E) at (0,4)  {$\mcal{O}_{[4]}$};
  \node at (1,0) {$0$};
  \node at (1,1) {$6$};
  \node at (1,2) {$8$};
  \node at (1,3) {$10$};
  \node at (1,4) {$12$};
  \node at (2,0) {\dynkin[x/.style={thin}]{A}{ooo}};
  \node at (2,0.9) {\dynkin[x/.style={thin}]{A}{oox}};
  \node at (2,1.1) {\dynkin[x/.style={thin}]{A}{xoo}};
  \node at (2,2) {\dynkin[x/.style={thin}]{A}{oxo}};
  \node at (2,2.8) {\dynkin[x/.style={thin}]{A}{oxx}};
  \node at (2,3) {\dynkin[x/.style={thin}]{A}{xox}};
  \node at (2,3.2) {\dynkin[x/.style={thin}]{A}{xxo}};
  \node at (2,4) {\dynkin[x/.style={thin}]{A}{xxx}};
  \draw [thin, -] (A) -- (B);
  \draw [thin, -] (B) -- (C);
  \draw [thin, -] (C) -- (D);
  \draw [thin, -] (D) -- (E);
\end{tikzpicture}}
\hspace{\the\parindent}
\caption{Richardson orbits for $\mfrak{sl}_n$}
\label{fig:nilpotent orbits sl}
\vspace{-2mm}
\end{figure}
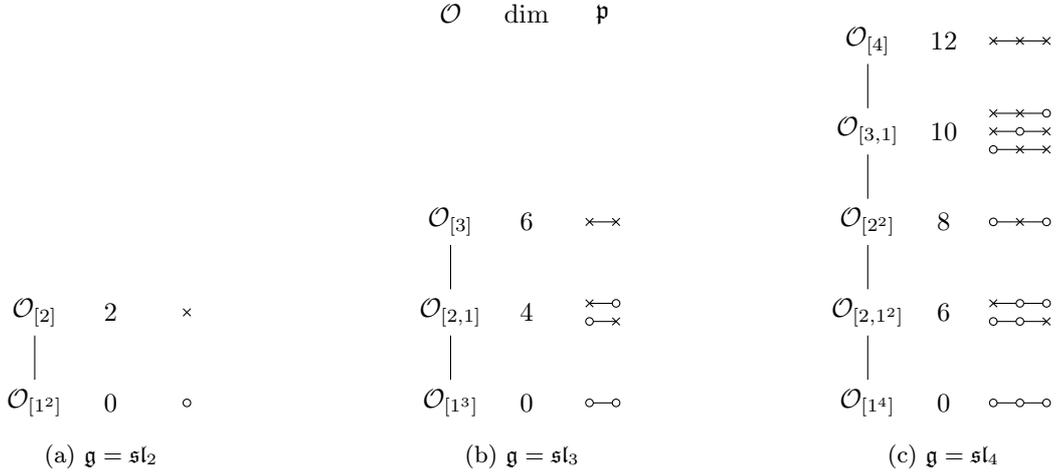

%%%%%%%%%%%%%%%%%%%%%%%%%%%%%%%%%%%%%%%%%%%%%%%%%%%%%%%%%%%%%%%%%%%%%%%%%%%%%%%%%%%%%%%%%%

\subsubsection{Lie algebra $\mfrak{sl}_2$}

We start with the simplest and well-known case $\mfrak{g}=\mfrak{sl}_2$ of admissible level $k \in \Q$, i.e.\ $n=1$ and
\begin{align*}
  k+2 = {p \over q}\ \text{with}\ p,q \in \N,\, (p,q)=1,\, p \geq 2.
\end{align*}
The Hasse diagram describing the structure of nilpotent orbits for $\mfrak{g}$ is given in Figure \ref{fig:sl2}. Let us note that $\mcal{O}_{[2]} = \mcal{O}_{\rm prin} = \mcal{O}_{\rm min}$ and $\mcal{O}_{[1^2]} = \mcal{O}_{\rm zero} = \mcal{O}_{\rm subreg}$. Therefore, we have a decomposition
\begin{align*}
  \widebar{{\rm Pr}}_k = \widebar{{\rm Pr}}_k^{\smash{\mcal{O}_{\rm zero}}} \sqcup \widebar{{\rm Pr}}_k^{\smash{\mcal{O}_{\rm prin}}},
\end{align*}
where
\begin{align*}
  \widebar{{\rm Pr}}_k^{\smash{\mcal{O}_{\rm zero}}} &= \widebar{{\rm Pr}}_{k,\Z} = \{\lambda_1\omega_1;\, \lambda_1 \in \N_0,\, \lambda_1 \leq p-2\}, \\
  \widebar{{\rm Pr}}_k^{\smash{\mcal{O}_{\rm prin}}} &= \Lambda_k(\mfrak{b}) = \bigg\{\lambda - {p\over q}\,a\omega_1;\, \lambda \in \widebar{{\rm Pr}}_{k,\Z},\, a \in \N,\, a \leq q-1\bigg\}.
\end{align*}
The Hasse diagram of $W^\mfrak{p}$ for the corresponding standard parabolic subalgebra $\mfrak{p}$ of $\mfrak{g}$ is given in Figure \ref{fig:hasse diagram Wp}. Let us note that $\smash{\widebar{{\rm Pr}}}_k^{\smash{\mcal{O}_{\rm prin}}}$ is non-empty if $q \geq 2$.

%%%%%%%%%%%%%%%%%%%%%%%%%%%%%%%%%%%%%%%%%%%%%%%%%%%%%%%%%%%%%%%%%%%%%%%%%%%%%%%%%%%%%%%%%%

\subsubsection{Lie algebra $\mfrak{sl}_3$}

The next case is $\mfrak{g}=\mfrak{sl}_3$ of admissible level $k\in \Q$, i.e.\ $n=2$ and
\begin{align*}
  k+3 = {p \over q}\ \text{with}\ p,q \in \N,\, (p,q)=1,\, p \geq 3.
\end{align*}
The Hasse diagram describing the structure of nilpotent orbits for $\mfrak{g}$ is given in Figure \ref{fig:sl3}. Let us note that $\mcal{O}_{[3]} = \mcal{O}_{\rm prin}$, $\mcal{O}_{[2,1]} = \mcal{O}_{\rm min} = \mcal{O}_{\rm subreg}$ and $\mcal{O}_{[1^3]} = \mcal{O}_{\rm zero}$. Therefore, we have a decomposition
\begin{align*}
  \widebar{{\rm Pr}}_k = \widebar{{\rm Pr}}_k^{\smash{\mcal{O}_{\rm zero}}} \sqcup \widebar{{\rm Pr}}_k^{\smash{\mcal{O}_{\rm min}}} \sqcup \widebar{{\rm Pr}}_k^{\smash{\mcal{O}_{\rm prin}}},
\end{align*}
where
\begin{align*}
  \widebar{{\rm Pr}}_k^{\smash{\mcal{O}_{\rm zero}}} &= \widebar{{\rm Pr}}_{k,\Z} = \{\lambda_1\omega_1+\lambda_2\omega_2;\, \lambda_1,\lambda_2 \in \N_0,\, \lambda_1+\lambda_2 \leq p-3\}, \\
  \widebar{{\rm Pr}}_k^{\smash{\mcal{O}_{\rm min}}} &= \Lambda_k(\mfrak{p}_{\alpha_1}^{\rm max}) \cup s_1 \cdot \Lambda_k(\mfrak{p}_{\alpha_1}^{\rm max}) \cup s_2s_1 \cdot \Lambda_k(\mfrak{p}_{\alpha_1}^{\rm max}),\\
  \widebar{{\rm Pr}}_k^{\smash{\mcal{O}_{\rm prin}}} &= \Lambda_k(\mfrak{b}) \cup s_1 \cdot \Lambda_k(\mfrak{b})
\end{align*}
and
\begin{align*}
  \Lambda_k(\mfrak{p}_{\alpha_1}^{\rm max}) &= \bigg\{\lambda-{p\over q}\,a\omega_1;\, \lambda \in \widebar{{\rm Pr}}_{k,\Z},\, a \in \N,\, a \leq q-1\bigg\}, \\
  \Lambda_k(\mfrak{b}) &= \bigg\{\lambda -{p\over q}(a\omega_1 + b\omega_2);\, \lambda \in \widebar{{\rm Pr}}_{k,\Z},\, a,b \in \N,\, a+b \leq q-1\bigg\}.
\end{align*}
The Hasse diagram of $W^\mfrak{p}$ for the corresponding standard parabolic subalgebra $\mfrak{p}$ of $\mfrak{g}$ is given in Figure \ref{fig:hasse diagram Wp}. Let us note that $\smash{\widebar{{\rm Pr}}}_k^{\smash{\mcal{O}_{\rm min}}}$ is non-empty if $q \geq 2$ and $\smash{\widebar{{\rm Pr}}}_k^{\smash{\mcal{O}_{\rm prin}}}$ is non-empty for $q\geq 3$.

\afterpage{\clearpage}
\begin{figure}[p]
\centering
\hfill
\subcaptionbox{$\smash{W^\mfrak{b}}$, $\mfrak{g}=\mfrak{sl}_2$\label{fig:sl2 W_prin}}[0.26\textwidth]
{\begin{tikzpicture}
[>=latex, yscale=1, xscale=1.5,point/.style={circle,draw=white,fill=black,line width=1.5mm, inner sep=1mm}, point2/.style={circle,draw=white,fill=gray,line width=1.5mm, inner sep=1mm}]
  \node (A) at (0,0) [point] {};
  \node (B) at (1,0) [point2] {};
 \draw [thin, ->] (A) to node[above] {$s_1$} (B);
\end{tikzpicture}}\hfill
\vspace{6mm}

\hfill
\subcaptionbox{$W^{\smash{\mfrak{p}_{\alpha_1}^{\rm max}}}$, $\mfrak{g}=\mfrak{sl}_3$\label{fig:sl3 W_min}}[0.26\textwidth]
{\begin{tikzpicture}
[>=latex, yscale=1, xscale=1.5,point/.style={circle,draw=white,fill=black,line width=1.5mm, inner sep=1mm}, point2/.style={circle,draw=white,fill=gray,line width=1.5mm, inner sep=1mm}]
  \node (A) at (0,0) [point] {};
  \node (B) at (1,0) [point] {};
  \node (C) at (2,0) [point] {};
  \draw [thin, ->] (A) to node[above] {$s_1$} (B);
  \draw [thin, ->] (B) to node[above] {$s_2$} (C);
\end{tikzpicture}}\hfill
\vspace{6mm}

\hfill
\subcaptionbox{$\smash{W^\mfrak{b}}$, $\mfrak{g}=\mfrak{sl}_3$\label{fig:sl3 W_prin}}[0.39\textwidth]
{\begin{tikzpicture}
[>=latex, yscale=1, xscale=1.5,point/.style={circle,draw=white,fill=black,line width=1.5mm, inner sep=1mm}, point2/.style={circle,draw=white,fill=gray,line width=1.5mm, inner sep=1mm}]
  \node (A) at (0,0) [point] {};
  \node (B) at (1,0.5) [point] {};
  \node (C) at (1,-0.5) [point2] {};
  \node (D) at (2,0.5) [point2] {};
  \node (E) at (2,-0.5) [point2] {};
  \node (F) at (3,0) [point2] {};
  \draw [thin, ->] (A) to node[above] {$s_1$} (B);
  \draw [thin, ->] (A) to node[below] {$s_2$} (C);
  \draw [thin, ->] (B) to node[above] {$s_2$} (D);
  \draw [thin, ->] (C) to node[below] {$s_1$} (E);
  \draw [thin, ->] (D) to node[above] {$s_1$} (F);
  \draw [thin, ->] (E) to node[below] {$s_2$} (F);
  \draw [thin,dotted, ->] (B) -- (E);
  \draw [thin,dotted, ->] (C) -- (D);
\end{tikzpicture}}\hfill
\vspace{6mm}

\hfill
\subcaptionbox{$W^{\smash{\mfrak{p}_{\alpha_1}^{\rm max}}}$, $\mfrak{g}=\mfrak{sl}_4$\label{fig:sl4 W_min}}[0.39\textwidth]
{\begin{tikzpicture}
[>=latex, yscale=1, xscale=1.5,point/.style={circle,draw=white,fill=black,line width=1.5mm, inner sep=1mm}, point2/.style={circle,draw=white,fill=gray,line width=1.5mm, inner sep=1mm}]
  \node (A) at (0,0) [point] {};
  \node (B) at (1,0) [point] {};
  \node (C) at (2,0) [point] {};
  \node (D) at (3,0) [point] {};
  \draw [thin, ->] (A) to node[above] {$s_1$} (B);
  \draw [thin, ->] (B) to node[above] {$s_2$} (C);
  \draw [thin, ->] (C) to node[above] {$s_3$} (D);
\end{tikzpicture}} \hfill
\vspace{6mm}

\hfill
\subcaptionbox{$W^{\smash{\mfrak{p}_{\alpha_2}^{\rm max}}}$, $\mfrak{g}=\mfrak{sl}_4$\label{fig:sl4 W_rect}}[0.52\textwidth]
{\begin{tikzpicture}
[>=latex, yscale=1, xscale=1.5,point/.style={circle,draw=white,fill=black,line width=1.5mm, inner sep=1mm}, point2/.style={circle,draw=white,fill=gray,line width=1.5mm, inner sep=1mm}]
  \node (A) at (0,0) [point] {};
  \node (B) at (1,0) [point] {};
  \node (C) at (2,0.5) [point] {};
  \node (D) at (2,-0.5) [point2] {};
  \node (E) at (3,0) [point2] {};
  \node (F) at (4,0) [point2] {};
  \draw [thin, ->] (A) to node[above] {$s_2$} (B);
  \draw [thin, ->] (B) to node[above] {$s_1$} (C);
  \draw [thin, ->] (B) to node[below] {$s_3$} (D);
  \draw [thin, ->] (C) to node[above] {$s_3$} (E);
  \draw [thin, ->] (D) to node[below] {$s_1$} (E);
  \draw [thin, ->] (E) to node[above] {$s_2$} (F);
\end{tikzpicture}} \hfill
\vspace{6mm}

\hfill
\subcaptionbox{$W^{\smash{\mfrak{p}_{\alpha_3}^{\rm min}}}$, $\mfrak{g}=\mfrak{sl}_4$\label{fig:sl4 W_subreg}}[0.65\textwidth]
{\begin{tikzpicture}
[>=latex, yscale=1, xscale=1.5,point/.style={circle,draw=white,fill=black,line width=1.5mm, inner sep=1mm}, point2/.style={circle,draw=white,fill=gray,line width=1.5mm, inner sep=1mm}]
  \node (A) at (0,0) [point] {};
  \node (B) at (1,0.5) [point] {};
  \node (C) at (1,-0.5) [point] {};
  \node (D) at (2,0.5) [point] {};
  \node (E) at (2,-0.5) [point] {};
  \node (F) at (2,-1.5) [point] {};
  \node (G) at (3,0.5) [point] {};
  \node (H) at (3,-0.5) [point] {};
  \node (I) at (3,-1.5) [point] {};
  \node (J) at (4,0.5) [point] {};
  \node (K) at (4,-0.5) [point] {};
  \node (L) at (5,0) [point] {};
  \draw [thin, ->] (A) to node[above] {$s_1$} (B);
  \draw [thin, ->] (A) to node[below] {$s_2$} (C);
  \draw [thin, ->] (B) to node[pos=0.4,above] {$s_2$} (E);
  \draw [thin, ->] (C) to node[pos=0.4,below] {$s_1$} (D);
  \draw [thin, ->] (C) to node[below] {$s_3$} (F);
  \draw [thin, ->] (D) to node[above] {$s_2$} (G);
  \draw [thin, ->] (D) to node[pos=0.4,above] {$s_3$} (H);
  \draw [thin, ->] (E) to node[pos=0.4,below] {$s_1$} (G);
  \draw [thin, ->] (E) to node[pos=0.4,above] {$s_3$} (I);
  \draw [thin, ->] (F) to node[pos=0.4,below] {$s_1$} (H);
  \draw [thin, ->] (G) to node[pos=0.4,above] {$s_3$} (K);
  \draw [thin, ->] (H) to node[pos=0.4,below] {$s_3$} (J);
  \draw [thin, ->] (I) to node[below] {$s_1$} (K);
  \draw [thin, ->] (J) to node[above] {$s_3$} (L);
  \draw [thin, ->] (K) to node[below] {$s_2$} (L);
  \draw [thin,dotted, ->] (B) -- (D);
  \draw [thin,dotted, ->] (C) -- (E);
  \draw [thin,dotted, ->] (F) -- (I);
  \draw [thin,dotted, ->] (G) -- (J);
  \draw [thin,dotted, ->] (H) -- (K);
\end{tikzpicture}} \hfill
\vspace{6mm}

\hfill
\subcaptionbox{$\smash{W^\mfrak{b}}$, $\mfrak{g}=\mfrak{sl}_4$\label{fig:sl4 W_prin}}[0.78\textwidth]
{\begin{tikzpicture}
[>=latex, yscale=1, xscale=1.5,point/.style={circle,draw=white,fill=black,line width=1.5mm, inner sep=1mm}, point2/.style={circle,draw=white,fill=gray,line width=1.5mm, inner sep=1mm}]
  \node (A) at (0,0) [point] {};
  \node (B) at (1,1) [point] {};
  \node (C) at (1,0) [point] {};
  \node (D) at (1,-1) [point] {};
  \node (E) at (2,2) [point] {};
  \node (F) at (2,1) [point2] {};
  \node (G) at (2,0) [point] {};
  \node (H) at (2,-1) [point2] {};
  \node (I) at (2,-2) [point2] {};
  \node (J) at (3,2) [point2] {};
  \node (K) at (3,1) [point2] {};
  \node (L) at (3,0) [point2] {};
  \node (M) at (3,-1) [point2] {};
  \node (N) at (3,-2) [point2] {};
  \node (O) at (3,-3) [point2] {};
  \node (P) at (4,1) [point2] {};
  \node (Q) at (4,0) [point2] {};
  \node (R) at (4,-1) [point2] {};
  \node (S) at (4,-2) [point2] {};
  \node (T) at (4,-3) [point2] {};
  \node (U) at (5,0) [point2] {};
  \node (V) at (5,-1) [point2] {};
  \node (W) at (5,-2) [point2] {};
  \node (X) at (6,-1) [point2] {};
  \draw [thin, ->] (A) to node[above] {$s_1$} (B);
  \draw [thin, ->] (A) to node[above] {$s_2$} (C);
  \draw [thin, ->] (A) to node[below] {$s_3$} (D);
  \draw [thin, ->] (B) to node[pos=0.3,yshift=-2,above] {$s_2$} (F);
  \draw [thin, ->] (B) to node[pos=0.6,below] {$s_3$} (G);
  \draw [thin, ->] (C) to node[pos=0.9,yshift=-2,below] {$s_1$} (E);
  \draw [thin, ->] (C) to node[pos=0.9,yshift=2,above] {$s_3$} (I);
  \draw [thin, ->] (D) to node[pos=0.6,above] {$s_1$} (G);
  \draw [thin, ->] (D) to node[pos=0.3,yshift=2,below] {$s_2$} (H);
  \draw [thin, ->] (E) to node[pos=0.4,above] {$s_2$} (K);
  \draw [thin, ->] (E) to node[pos=0.9,yshift=2,above] {$s_3$} (L);
  \draw [thin, ->] (F) to node[pos=0.15,yshift=2,below] {$s_1$} (K);
  \draw [thin, ->] (F) to node[pos=0.7,yshift=4,xshift=2,above] {$s_3$} (O);
  \draw [thin, ->] (G) to node[pos=0.55,yshift=-1,above] {$s_2$} (M);
  \draw [thin, ->] (H) to node[pos=0.55,yshift=-3,xshift=1,below] {$s_1$} (J);
  \draw [thin, ->] (H) to node[pos=0.,yshift=1,below] {$s_3$} (N);
  \draw [thin, ->] (I) to node[pos=0.2,yshift=-2,below] {$s_1$} (L);
  \draw [thin, ->] (I) to node[pos=0.3,yshift=2,below] {$s_2$} (N);
  \draw [thin, ->] (J) to node[above] {$s_2$} (P);
  \draw [thin, ->] (J) to node[pos=0.1,yshift=-2,below] {$s_3$} (Q);
  \draw [thin, ->] (K) to node[pos=0.05,xshift=-1,yshift=-4,below] {$s_3$} (S);
  \draw [thin, ->] (L) to node[pos=0.4,yshift=3,xshift=-1,below] {$s_2$} (R);
  \draw [thin, ->] (M) to node[pos=0.55,yshift=2,above] {$s_1$} (P);
  \draw [thin, ->] (M) to node[pos=0.6,yshift=2,above] {$s_3$} (T);
  \draw [thin, ->] (N) to node[pos=0.5,yshift=-2,below] {$s_1$} (Q);
  \draw [thin, ->] (O) to node[pos=0.4,below] {$s_1$} (S);
  \draw [thin, ->] (O) to node[below] {$s_2$} (T);
  \draw [thin, ->] (P) to node[pos=0.1,yshift=-2,below] {$s_3$} (V);
  \draw [thin, ->] (Q) to node[pos=0.3,yshift=2,below] {$s_2$} (U);
  \draw [thin, ->] (R) to node[pos=0.4,below] {$s_1$} (U);
  \draw [thin, ->] (R) to node[pos=0.4,above] {$s_3$} (W);
  \draw [thin, ->] (S) to node[pos=0.3,yshift=-2,above] {$s_2$} (W);
  \draw [thin, ->] (T) to node[pos=0.1,yshift=2,above] {$s_1$} (V);
  \draw [thin, ->] (U) to node[above] {$s_3$} (X);
  \draw [thin, ->] (V) to node[above] {$s_2$} (X);
  \draw [thin, ->] (W) to node[below] {$s_1$} (X);
  \draw [thin,dotted, ->] (B) -- (E);
  \draw [thin,dotted, ->] (C) -- (F);
  \draw [thin,dotted, ->] (C) -- (H);
  \draw [thin,dotted, ->] (D) -- (I);
  \draw [thin,dotted, ->] (E) -- (J);
  \draw [thin,dotted, ->] (F) -- (M);
  \draw [thin,dotted, ->] (G) -- (J);
  \draw [thin,dotted, ->] (G) -- (L);
  \draw [thin,dotted, ->] (G) -- (O);
  \draw [thin,dotted, ->] (H) -- (M);
  \draw [thin,dotted, ->] (I) -- (O);
  \draw [thin,dotted, ->] (K) -- (P);
  \draw [thin,dotted, ->] (K) -- (R);
  \draw [thin,dotted, ->] (L) -- (Q);
  \draw [thin,dotted, ->] (L) -- (S);
  \draw [thin,dotted, ->] (M) -- (R);
  \draw [thin,dotted, ->] (N) -- (R);
  \draw [thin,dotted, ->] (N) -- (T);
  \draw [thin,dotted, ->] (P) -- (U);
  \draw [thin,dotted, ->] (Q) -- (V);
  \draw [thin,dotted, ->] (S) -- (V);
  \draw [thin,dotted, ->] (T) -- (W);
\end{tikzpicture}} \hfill
\hspace{\the\parindent}
\caption{Hasse diagrams of $W^\mfrak{p}$}
\label{fig:hasse diagram Wp}
%\vspace{-2mm}
\end{figure}

%%%%%%%%%%%%%%%%%%%%%%%%%%%%%%%%%%%%%%%%%%%%%%%%%%%%%%%%%%%%%%%%%%%%%%%%%%%%%%%%%%%%%%%%%%

\subsubsection{Lie algebra $\mfrak{sl}_4$}
\label{subsec:lie algebra sl4}

The last case is $\mfrak{g}=\mfrak{sl}_4$ of admissible level $k\in \Q$, i.e.\ $n=3$ and
\begin{align*}
  k+4 = {p \over q}\ \text{with}\ p,q \in \N,\, (p,q)=1,\, p \geq 4.
\end{align*}
The Hasse diagram describing the structure of nilpotent orbits for $\mfrak{g}$ is given in Figure \ref{fig:sl4}. Let us note that $\mcal{O}_{[4]} = \mcal{O}_{\rm prin}$, $\mcal{O}_{[3,1]} = \mcal{O}_{\rm subreg}$, $\mcal{O}_{[2^2]} = \mcal{O}_{\rm rect}$, $\mcal{O}_{[2,1^2]} = \mcal{O}_{\rm min}$ and $\mcal{O}_{[1^4]} = \mcal{O}_{\rm zero}$. Therefore, we have a decomposition
\begin{align*}
  \widebar{{\rm Pr}}_k = \widebar{{\rm Pr}}_k^{\smash{\mcal{O}_{\rm zero}}} \sqcup \widebar{{\rm Pr}}_k^{\smash{\mcal{O}_{\rm min}}} \sqcup \widebar{{\rm Pr}}_k^{\smash{\mcal{O}_{\rm rect}}} \sqcup \widebar{{\rm Pr}}_k^{\smash{\mcal{O}_{\rm subreg}}} \sqcup \widebar{{\rm Pr}}_k^{\smash{\mcal{O}_{\rm prin}}},
\end{align*}
where
\begin{align*}
  \widebar{{\rm Pr}}_k^{\smash{\mcal{O}_{\rm zero}}} &= \widebar{{\rm Pr}}_{k,\Z} = \{\lambda_1\omega_1+\lambda_2\omega_2 +\lambda_3\omega_3;\, \lambda_1,\lambda_2, \lambda_3 \in \N_0,\, \lambda_1+\lambda_2+\lambda_3 \leq p-4\}, \\
  \widebar{{\rm Pr}}_k^{\smash{\mcal{O}_{\rm min}}} &= \Lambda_k(\mfrak{p}_{\alpha_1}^{\rm max}) \cup s_1 \cdot \Lambda_k(\mfrak{p}_{\alpha_1}^{\rm max}) \cup s_2s_1 \cdot \Lambda_k(\mfrak{p}_{\alpha_1}^{\rm max}) \cup s_3s_2s_1 \cdot \Lambda_k(\mfrak{p}_{\alpha_1}^{\rm max}),\\
  \widebar{{\rm Pr}}_k^{\smash{\mcal{O}_{\rm rect}}} &= \Lambda_k(\mfrak{p}_{\alpha_2}^{\rm max}) \cup s_2 \cdot \Lambda_k(\mfrak{p}_{\alpha_2}^{\rm max}) \cup s_1s_2 \cdot \Lambda_k(\mfrak{p}_{\alpha_2}^{\rm max}),\\
  \widebar{{\rm Pr}}_k^{\smash{\mcal{O}_{\rm subreg}}} &= \Lambda_k(\mfrak{p}_{\alpha_3}^{\rm min}) \cup s_1 \cdot \Lambda_k(\mfrak{p}_{\alpha_3}^{\rm min}) \cup s_2 \cdot \Lambda_k(\mfrak{p}_{\alpha_3}^{\rm min}) \cup s_1s_2 \cdot \Lambda_k(\mfrak{p}_{\alpha_3}^{\rm min}) \cup s_2s_1 \cdot \Lambda_k(\mfrak{p}_{\alpha_3}^{\rm min}) \\
  & \quad \cup s_3s_2 \cdot \Lambda_k(\mfrak{p}_{\alpha_3}^{\rm min}) \cup s_1s_2s_1 \cdot \Lambda_k(\mfrak{p}_{\alpha_3}^{\rm min}) \cup s_3s_1s_2 \cdot \Lambda_k(\mfrak{p}_{\alpha_3}^{\rm min}) \cup s_3s_2s_1 \cdot \Lambda_k(\mfrak{p}_{\alpha_3}^{\rm min}) \\
  & \quad \cup s_2s_3s_1s_2 \cdot \Lambda_k(\mfrak{p}_{\alpha_3}^{\rm min}) \cup s_3s_1s_2s_1 \cdot \Lambda_k(\mfrak{p}_{\alpha_3}^{\rm min}) \cup s_2s_3s_1s_2s_1 \cdot \Lambda_k(\mfrak{p}_{\alpha_3}^{\rm min}),\\
  \widebar{{\rm Pr}}_k^{\smash{\mcal{O}_{\rm prin}}} &= \Lambda_k(\mfrak{b}) \cup s_1 \cdot \Lambda_k(\mfrak{b}) \cup s_2 \cdot \Lambda_k(\mfrak{b}) \cup s_3 \cdot \Lambda_k(\mfrak{b}) \cup s_1s_2 \cdot \Lambda_k(\mfrak{b}) \cup s_1s_3 \cdot \Lambda_k(\mfrak{b})
\end{align*}
and
\begin{align*}
  \Lambda_k(\mfrak{p}_{\alpha_1}^{\rm max}) &= \bigg\{\lambda-{p\over q}\,a\omega_1;\, \lambda \in \widebar{{\rm Pr}}_{k,\Z},\, a \in \N,\, a \leq q-1\bigg\}, \\
  \Lambda_k(\mfrak{p}_{\alpha_2}^{\rm max}) &= \bigg\{\lambda -{p\over q}\,a\omega_2;\, \lambda \in \widebar{{\rm Pr}}_{k,\Z},\, a \in \N,\, a \leq q-1\bigg\}, \\
  \Lambda_k(\mfrak{p}_{\alpha_3}^{\rm min}) &= \bigg\{\lambda -{p\over q}(a\omega_1 + b\omega_2);\, \lambda \in \widebar{{\rm Pr}}_{k,\Z},\, a,b \in \N,\, a+b \leq q-1\bigg\}, \\
  \Lambda_k(\mfrak{b}) &= \bigg\{\lambda -{p\over q}(a\omega_1 + b\omega_2 + c\omega_3);\, \lambda \in \widebar{{\rm Pr}}_{k,\Z},\, a,b,c \in \N,\, a+b+c \leq q-1\bigg\}.
\end{align*}
The Hasse diagram of $W^\mfrak{p}$ for the corresponding standard parabolic subalgebra $\mfrak{p}$ of $\mfrak{g}$ is given in Figure \ref{fig:hasse diagram Wp}. Let us note that $\smash{\widebar{{\rm Pr}}}_k^{\smash{\mcal{O}_{\rm min}}}$ and $\smash{\widebar{{\rm Pr}}}_k^{\smash{\mcal{O}_{\rm rect}}}$ are non-empty if $q \geq 2$, $\smash{\widebar{{\rm Pr}}}_k^{\smash{\mcal{O}_{\rm subreg}}}$ is non-empty for $q\geq 3$ and $\smash{\widebar{{\rm Pr}}}_k^{\smash{\mcal{O}_{\rm prin}}}$ is non-empty provided $q\geq 4$.

%%%%%%%%%%%%%%%%%%%%%%%%%%%%%%%%%%%%%%%%%%%%%%%%%%%%%%%%%%%%%%%%%%%%%%%%%%%%%%%%%%%%%%%%%%

\subsubsection{Admissible levels with small denominator}
\label{subsec:small denominator}

The simplest possible and very well-known case is $\mfrak{g}$ of admissible level $k \in \Q$ with denominator $q=1$, i.e.\ we have
\begin{align*}
   k+n+1 = p\ \text{with}\ p \in \N,\, p \geq n+1.
\end{align*}
By Theorem \ref{thm:nilpotent orbits} we have $\mcal{O}_1 = \mcal{O}_{[1^{n+1}]} = \mcal{O}_{\rm zero}$. Moreover, since $\smash{\widebar{{\mcal{O}}_1}} = \mcal{O}_{\rm zero}$, we easily obtain that
\begin{align}
  \widebar{{\rm Pr}}_k = \widebar{{\rm Pr}}_k^{\smash{\mcal{O}_{\rm zero}}} = \widebar{{\rm Pr}}_{k,\Z}
\end{align}
by Theorem \ref{thm:zero orbit}.
\smallskip

The more interesting case is $\mfrak{g}$ of admissible level $k \in \Q$ with denominator $q=2$, i.e.\ we have
\begin{align*}
   k+n+1 = {p \over 2}\ \text{with}\ p \in \N,\, (p,2) =1,\, p \geq n+1.
\end{align*}
By Theorem \ref{thm:nilpotent orbits} we have
\begin{align*}
  \mcal{O}_2 = \begin{cases}
    \mcal{O}_{[2^m]} & \text{if $n=2m-1$}, \\
    \mcal{O}_{[2^m,1]} & \text{if $n=2m$},
  \end{cases}
\end{align*}
which implies that $\mcal{O}_\lambda \subset \smash{\widebar{{\mcal{O}}_2}}$ provided $\lambda =[2^r,1^{n+1-2r}]$ for $r=0,1,\dots,m$. Therefore, by \eqref{eq:Pr_k decomposition} we have a decomposition
\begin{align*}
  \widebar{{\rm Pr}}_k = \widebar{{\rm Pr}}_k^{\smash{\mcal{O}_{\rm zero}}} \sqcup \bigsqcup_{r=1}^m\, \widebar{{\rm Pr}}_k^{\smash{\mcal{O}_{[2^r,1^{n+1-2r}]}}}.
\end{align*}
For $r=1,2,\dots,m$, we denote by $\lambda_r$ the partition $[2^r,1^{n+1-2r}]$ of $n+1$. As $\lambda_r^t = [n+1-r,r]$, we have
\begin{align*}
  [\mfrak{p}_{\lambda_r^t}] = \begin{cases}
    \{\mfrak{p}_{\alpha_r}^{\rm max}, \mfrak{p}_{\alpha_{n+1-r}}^{\rm max}\} & \text{if $n+1 \neq 2r$}, \\
    \{\mfrak{p}_{\alpha_r}^{\rm max}\} & \text{if $n+1 = 2r$}
  \end{cases}
\end{align*}
for $r=1,2,\dots,m$. Besides, we see that $[\mfrak{p}_{\alpha_r}^{\rm max}]_+= [\mfrak{p}_{\lambda_r^t}]$ together with
\begin{align*}
  W_+^{\smash{\mfrak{p}_{\alpha_r}^{\rm max}}} = \begin{cases}
    \{e\} & \text{if $n+1 \neq 2r$}, \\
    \{e, w_r\} & \text{if $n+1 = 2r$}.
  \end{cases}
\end{align*}
Hence, by using Theorem \ref{thm:Pr_k parametrization} we may write
\begin{align}
  \widebar{{\rm Pr}}_k = \widebar{{\rm Pr}}_{k,\Z} \sqcup  \bigsqcup_{r=1}^m \bigsqcup_{w \in W^{\smash{\mfrak{p}_{\alpha_r}^{\rm max}}}} w^{-1} \cdot \Lambda_k(\mfrak{p}_{\alpha_r}^{\rm max})
\end{align}
if $n$ is even and
\begin{align}
  \widebar{{\rm Pr}}_k = \widebar{{\rm Pr}}_{k,\Z} \sqcup  \bigsqcup_{r=1}^{m-1} \bigsqcup_{w \in W^{\smash{\mfrak{p}_{\alpha_r}^{\rm max}}}} w^{-1} \cdot \Lambda_k(\mfrak{p}_{\alpha_r}^{\rm max}) \sqcup \bigsqcup_{[w] \in W_+^{\smash{\mfrak{p}_{\alpha_m}^{\rm max}}} \backslash W^{\smash{\mfrak{p}_{\alpha_m}^{\rm max}}}} w^{-1} \cdot \Lambda_k(\mfrak{p}_{\alpha_m}^{\rm max})
\end{align}
if $n$ is odd, where $m = \lfloor {n+1 \over 2} \rfloor$ and
\begin{align*}
  \Lambda_k(\mfrak{p}_{\alpha_r}^{\rm max}) = \bigg\{\lambda - {p\over 2}\,\omega_r;\, \lambda \in \widebar{{\rm Pr}}_{k,\Z} \bigg\}
\end{align*}
for $r=1,2,\dots,n$.

The Hasse diagram $W^\mfrak{p}$ for the standard parabolic subalgebra $\mfrak{p}=\mfrak{p}_\alpha^{\rm max}$ with $\alpha \in \Pi$ we can determine by the methods described at the end of Section \ref{subsec:Primitive ideals}. However, since $\mfrak{p}$ is a maximal parabolic subalgebra of $\mfrak{g}$, which in other words means that the corresponding nilradical $\mfrak{u}$ is commutative, there is a relatively easy way how to compute the Hasse diagram $W^\mfrak{p}$, see \cite{Enright-Hunziker-Pruett2014} for more details.

%%%%%%%%%%%%%%%%%%%%%%%%%%%%%%%%%%%%%%%%%%%%%%%%%%%%%%%%%%%%%%%%%%%%%%%%%%%%%%%%%%%%%%%%%%%
%%%%%%%%%%%%%%%%%%%%%%%%%%%%%%%%%%%%%%%%%%%%%%%%%%%%%%%%%%%%%%%%%%%%%%%%%%%%%%%%%%%%%%%%%%%

\section{Tame and strongly tame Gelfand--Tsetlin modules}

%%%%%%%%%%%%%%%%%%%%%%%%%%%%%%%%%%%%%%%%%%%%%%%%%%%%%%%%%%%%%%%%%%%%%%%%%%%%%%%%%%%%%%%%%%%

\subsection{Tame Gelfand--Tsetlin modules}

Let us consider the simple Lie algebra $\mfrak{g}=\mfrak{sl}_{n+1}$ for $n \in \N$. A Cartan subalgebra $\mfrak{h}$ of $\mfrak{g}$ is given by diagonal matrices
\begin{align*}
  \mfrak{h} = \{\diag(a_1,a_2,\dots,a_{n+1});\, a_1,a_2,\dots,a_{n+1} \in \C,\, {\textstyle \sum_{i=1}^{n+1}} a_i=0\}.
\end{align*}
For $i=1,2,\dots,{n+1}$ we define $\veps_i \in \mfrak{h}^*$ through $\veps_i(\diag(a_1,a_2,\dots,a_{n+1}))=a_i$. The root system of $\mfrak{g}$ with respect to $\mfrak{h}$ is then given by $\Delta=\{\veps_i-\veps_j;\, 1\leq i \neq j \leq n+1\}$. A positive root system in $\Delta$ is $\Delta_+=\{\veps_i - \veps_j;\, 1 \leq i < j \leq n+1\}$ with the set of simple roots $\Pi_{\rm st}=\{\alpha_1,\alpha_2,\dots,\alpha_n\}$, $\alpha_i= \veps_i - \veps_{i+1}$ for $i=1,2,\dots,n$. We will also use the notation
\begin{align*}
  \alpha_{i,j} = \sum_{k=i}^j \alpha_k = \veps_i - \veps_{j+1}
\end{align*}
for $1 \leq i < j \leq n$. The fundamental weights are
\begin{align*}
\omega_i = \sum_{j=1}^i \veps_j - {i \over n+1} \sum_{j=1}^{n+1} \veps_j
\end{align*}
for $i=1,2,\dots,n$. Further, we define root vectors of $\mfrak{g}$ by
\begin{align*}
  e_{\veps_i-\veps_j} = e_{i,j} = E_{i,j} \quad \text{and} \quad f_{\veps_i-\veps_j} = f_{i,j} = E_{j,i}
\end{align*}
together with the corresponding coroot
\begin{align*}
  h_{\veps_i-\veps_j} = E_{i,i} - E_{j,j}
\end{align*}
for $1 \leq i < j \leq n+1$, where $E_{i,j} \in M_{n+1\times n+1}(\C)$ for $1\leq i,j \leq n+1$ is the $(n+1 \times n+1)$-matrix having $1$ at the intersection of the $i$-th row and the $j$-th column and $0$ elsewhere. In particular, we have
\begin{align*}
   e_i = E_{i,i+1}, \qquad h_i = E_{i,i}-E_{i+1,i+1}, \qquad f_i = E_{i+1,i}
\end{align*}
for $i =1,2,\dots,n$, for the Chevalley generators of $\mfrak{g}$. Besides, we denote by $\mfrak{b}_{\rm st}$ the standard Borel subalgebra of $\mfrak{g}$ with respect to the positive root system $\Delta_+$. Let us note that for two Borel subalgebras $\mfrak{b}_1, \mfrak{b}_2$ of $\mfrak{g}$ there exists an automorphism $\varphi \in \Aut(\mfrak{g})$ such that $\varphi(\mfrak{b}_1) = \mfrak{b}_2$.

Further, let us recall that there is the Tits extension $W_{\rm Tits}$ of the Weyl group $W$ of $\mfrak{g}$ which fits into the short exact sequence
\begin{align*}
  1 \rarr \Z_2^{\dim \mfrak{h}} \rarr W_{\rm Tis} \rarr W \rarr 1
\end{align*}
of finite groups. Although, the Weyl group $W$ does not have an action on $\mfrak{g}$, its Tits extension $W_{\rm Tits}$ has an action on $\mfrak{g}$ given as follows. For $\alpha \in \Delta$ and an $\mfrak{sl}_2$-triple $(e_\alpha,h_\alpha,f_\alpha)$, we define $r_\alpha \in \Aut(\mfrak{g})$ by
\begin{align*}
  r_\alpha = \exp(\ad(f_\alpha))\exp(-\ad(e_\alpha))\exp(\ad(f_\alpha)).
\end{align*}
Then the assignment $\tilde{s}_{\alpha_i} \mapsto r_{\alpha_i}$ for $i=1,2,\dots,n$, where $\tilde{s}_{\alpha_i}\!$ is the generator of $W_{\rm Tits}$ lifting the generator $s_{\alpha_i}\!$ of $W$, gives us a representation of $W_{\rm Tits}$ on $\mfrak{g}$. Moreover, we have $r_\alpha{}_{|\mfrak{h}}=s_\alpha$ and $r_\alpha(\mfrak{g}_\beta) = \mfrak{g}_{s_\alpha(\beta)}$ for $\alpha,\beta \in \Delta$. By a representative $\dot{w} \in \Aut(\mfrak{g})$ of $w \in W$ we will mean an element $\dot{w} = r_{\alpha_1}r_{\alpha_2}\cdots r_{\alpha_k}$ if $w = s_{\alpha_1}s_{\alpha_2}\cdots s_{\alpha_k}$, where $k\in\N_0$ and $\alpha_i \in \Pi_{\rm st}$ for $i=1,2,\dots,k$.

We will also need an automorphism $\sigma \in \Aut(\mfrak{g})$ order $2$ induced from a symmetric transformation
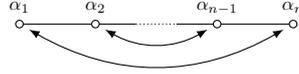
\vspace{-1mm}
\begin{align*}
    \begin{dynkinDiagram}[involutions={14;23}, labels*={\alpha_1,\alpha_2,\alpha_{n-1},\alpha_n}, edge length=10mm, x/.style={thin}]{A}{oo..oo}
    \end{dynkinDiagram}
\end{align*}
of order $2$ of the Dynkin diagram of $\mfrak{g}$ and determined by
\begin{align*}
    e_i \mapsto e_{n+1-i}, \qquad h_i \mapsto -h_{n+1-i}, \qquad f_i \mapsto f_{n+1-i}
\end{align*}
for $i=1,2,\dots,n$. It follows immediately that $\sigma(\alpha_i) = \alpha_{n+1-i}$ for $i=1,2,\dots,n$ and $\sigma(\mfrak{b}_{\rm st}) = \mfrak{b}_{\rm st}$. We denote by $\widebar{W}$ the subgroup of $\Aut(\mfrak{g})$ generated by $\sigma$ and by $r_{\alpha_i}\!$ for $i=1,2,\dots,n$.
\medskip

For a commutative algebra $\Gamma$, we denote by $\Hom(\Gamma,\C)$ the set of all characters of $\Gamma$, i.e.\ algebra homomorphisms from $\Gamma$ to $\C$. Let $M$ be a $\Gamma$-module. For each character $\chi \in \Hom(\Gamma,\C)$ we set
\begin{align*}
  M_\chi = \{v \in M;\, (\exists k \in \N)\,(\forall a \in \Gamma)\, (a-\chi(a))^kv=0\}.
\end{align*}
We say that $M$ is a \emph{$\Gamma$-weight module} if
\begin{align*}
  M = \bigoplus_{\chi \in \Hom(\Gamma,\C)} M_\chi.
\end{align*}
If $M_\chi \neq\{0\}$ then $\chi$ is a $\Gamma$-weight of $M$. The set of all $\Gamma$-weights of $M$ is the \emph{$\Gamma$-support} of $M$.  The dimension of  $M_\chi$ is the \emph{$\Gamma$-multiplicity} of $\chi$ in $M$.

We will be mainly interested in commutative subalgebras of $U(\mfrak{g})$ constructed in the following way. Let
\begin{align*}
  \mcal{F}\colon \mfrak{g}_1 \subset \mfrak{g}_2 \subset \dots \subset \mfrak{g}_n = \mfrak{g}
\end{align*}
be a chain of Lie subalgebras of $\mfrak{g}$ such that $\mfrak{g}_k \simeq \mfrak{sl}_{k+1}$ and $\mfrak{g}_k \cap \mfrak{h}$ is a Cartan subalgebra of $\mfrak{g}_k$ for $k=1,2,\dots,n$. The \emph{Gelfand--Tsetlin subalgebra} $\Gamma$ of $U(\mfrak{g})$ with respect to $\mcal{F}$ is generated by the Cartan subalgebra $\mfrak{h}$ and by the center of $U(\mfrak{g}_k)$ for $k=1,2,\dots,n$, see \cite{Drozd-Futorny-Ovsienko1994}. It is known that $\Gamma$ is a maximal commutative subalgebra of $U(\mfrak{g})$. A finitely generated $\Gamma$-weight $\mfrak{g}$-module $M$ is called a \emph{$\Gamma$-Gelfand--Tsetlin $\mfrak{g}$-module}.

Let $\Pi=\{\gamma_1,\gamma_2,\dots,\gamma_n\}$ be a set of simple roots in $\Delta$. Let us note that we will understand $\Pi$ as an ordered set. Then the most important example of a Gelfand--Tsetlin subalgebra of $U(\mfrak{g})$ is the \emph{standard} Gelfand--Tsetlin subalgebra $\Gamma_{\rm st}(\Pi)$ with respect to $\Pi$ which is attached to a chain $\mcal{F}$ uniquely determined by a sequence
\begin{align*}
  \{\gamma_1\} \subset \{\gamma_1, \gamma_2\} \subset \dots \subset \{\gamma_1,\gamma_2, \dots, \gamma_n\}
\end{align*}
of sets of simple roots of the corresponding Lie subalgebras in the chain $\mcal{F}$. Moreover, since there exists $w \in \widebar{W}$ satisfying $\Pi=w(\Pi_{\rm st})$, we get immediately that $\Gamma_{\rm st}(\Pi) = w(\Gamma_{\rm st}(\Pi_{\rm st}))$. We will use the notation $\Gamma_{\rm st}$ instead of $\Gamma_{\rm st}(\Pi_{\rm st})$ and $\Pi_w$ rather than $w(\Pi_{\rm st})$ for $w \in \widebar{W}$. We denote by $\mfrak{b}_\Pi$ the Borel subalgebra of $\mfrak{g}$ associated to $\Pi$ and by $\rho_\mfrak{b}$ the corresponding Weyl vector, i.e.\ the half-sum of positive roots (with respect to $\Pi$) of $\mfrak{g}$.
\medskip

\example{If $\mfrak{p} \supset \mfrak{b}_{\rm st}$ is a standard parabolic subalgebra of $\mfrak{g}$ and $E$ is a weight  $\mfrak{p}$-module with finite weight multiplicities and trivial action of the nilradical $\mfrak{u}$ of $\mfrak{p}$, then both $\mfrak{g}$-modules $M^\mfrak{g}_\mfrak{p}(E)$ and $L^\mfrak{g}_\mfrak{p}(E)$ are $\Gamma_{\rm st}(\Pi_w)$-Gelfand--Tsetlin $\mfrak{g}$-modules for all $w\in \widebar{W}$. In particular, any highest weight $\mfrak{g}$-module is a $\Gamma_{\rm st}(\Pi_w)$-Gelfand--Tsetlin $\mfrak{g}$-module for any $w\in \widebar{W}$.}
\medskip

\definition{Let $\Gamma$ be a Gelfand--Tseltin subalgebra of $U(\mfrak{g})$. We say that a $\Gamma$-Gelfand--Tsetlin $\mfrak{g}$-module $M$ is \emph{tame} if $\Gamma$ has a simple spectrum on $M$, i.e.\ all $\Gamma$-multiplicities are equal to $1$. In this case $\Gamma$-weights of $M$ parameterize a basis of $M$.}

Finite-dimensional $\mfrak{g}$-modules are examples of tame $\Gamma_{\rm st}(\Pi_w)$-Gelfand--Tsetlin $\mfrak{g}$-modules for arbitrary $w\in \widebar{W}$. Various classes of infinite-dimensional tame Gelfand--Tsetlin $\mfrak{g}$-modules are known, e.g.\ \emph{generic} modules \cite{Futorny-Ramirez-Zhang2018}, or more generally, \emph{relation} modules \cite{Futorny-Ramirez-Zhang2019}. Moreover, relation $\Gamma_{\rm st}$-Gelfand--Tsetlin $\mfrak{g}$-modules can be characterized as those $\mfrak{g}$-modules having an eigenbasis for $\Gamma_{\rm st}$ formed by a set of tableaux $T(v)$ for $v \in \smash{\C^{{n(n+1) \over 2}}}$, where $T(v)$ is the tableau associated to a vector
\begin{align*}
  v = (v_{n+1,1},\dots,v_{n+1,n+1}|v_{n,1},\dots,v_{n,n}|\dots|v_{2,1},v_{2,2}|v_{1,1}),
\end{align*}
and admiting the action of the Chevalley generators in the form
\begin{align}\label{Gelfand--Tsetlin formulas}
\begin{aligned}
e_k(T(v)) &= -\sum_{i=1}^k\! \left(\frac{\prod_{j=1}^{k+1}( v_{k,i}- v_{k+1,j})}{\prod_{j\neq i}^{k}( v_{k,i}- v_{k,j})}\!\right)T(v+\delta^{k,i}),\\[1mm]
f_k(T(v)) &= \sum_{i=1}^k\! \left(\frac{\prod_{j=1}^{k-1}( v_{k,i}- v_{k-1,j})}{\prod_{j\neq i}^{k}( v_{k,i}- v_{k,j})}\!\right)T(v-\delta^{k,i}),\\[1mm]
h_k(T(v)) &= \left(2\sum_{i=1}^k v_{k,i}-\sum_{i=1}^{k-1} v_{k-1,i}-\sum_{i=1}^{k+1} v_{k+1,i}-1\!\right)T(v),
\end{aligned}
\end{align}
where $\delta^{k,i}$ is the vector having $1$ at the position $(k,i)$ and $0$ elsewhere, for $k=1,2,\dots,n$. If a denominator equals zero, then the summand is assumed to be zero. For such a relation  $\mfrak{g}$-module, we will use the name \emph{strongly tame} $\Gamma_{\rm st}$-Gelfand--Tsetlin $\mfrak{g}$-module.
\medskip

Let $M$ be a $\mfrak{g}$-module with the action of $\mfrak{g}$ given via a homomorphism $\pi \colon \mfrak{g} \rarr \End M$ of Lie algebras. For $w \in \widebar{W}$, we define the twisted $\mfrak{g}$-module $M^w$ in such a way that as a vector space it coincides with $M$, however the action of $\mfrak{g}$ is given through a homomorphism $\pi_w \colon \mfrak{g} \rarr \End M$ of Lie algebras defined by $\pi_w(a) = \pi(w^{-1}(a))$ for $a \in \mfrak{g}$.
\medskip

Let $\Pi$ be a set of simple roots in $\Delta$. We say that $M$ is a strongly tame $\Gamma_{\rm st}(\Pi)$-Gelfand--Tsetlin $\mfrak{g}$-module if there is $w \in \widebar{W}$ such that $\Pi=w(\Pi_{\rm st})$ and $\smash{M^{w^{-1}}}\!$ is a strongly tame $\Gamma_{\rm st}$-Gelfand--Tsetlin $\mfrak{g}$-module.
\medskip

The following useful proposition is an immediate consequence of the definition of strongly tame Gelfand--Tsetlin $\mfrak{g}$-modules.
\medskip

\proposition{\label{prop-tame-twisted} Let $\Pi$ be a set of simple roots in $\Delta$ and let $M$ be a strongly tame $\Gamma_{\rm st}(\Pi)$-Gelfand--Tsetlin $\mfrak{g}$-module. Then $M^w$ is a strongly tame $\Gamma_{\rm st}(w(\Pi))$-Gelfand--Tsetlin $\mfrak{g}$-module for $w \in \widebar{W}$.}

Again, finite-dimensional $\mfrak{g}$-modules are examples of strongly tame $\Gamma_{\rm st}(\Pi_w)$-Gelfand--Tsetlin $\mfrak{g}$-modules for any $w\in \widebar{W}$. Besides, strongly tame highest weight $\Gamma_{\rm st}$-Gelfand--Tsetlin $\mfrak{g}$-modules were described in \cite[Theorem 4.8]{Futorny-Morales-Ramirez2020}.
\medskip

Let $\mfrak{b}=\mfrak{b}_\Pi$ for a set of simple roots $\Pi=\{\gamma_1,\gamma_2,\dots,\gamma_n\}$ in $\Delta$.
\medskip

\theorem{\cite[Theorem 4.8]{Futorny-Morales-Ramirez2020} Let $\lambda \in \mfrak{h}^*$ be a regular dominant weight. Then the simple $\mfrak{g}$-module $L^\mfrak{g}_\mfrak{b}(\lambda)$ is a strongly tame $\Gamma_{\rm st}(\Pi)$-Gelfand--Tsetlin $\mfrak{g}$-module.}

Using the fact that $L^\mfrak{g}_{\smash{w(\mfrak{b})}}(w(\lambda)) \simeq L^\mfrak{g}_\mfrak{b}(\lambda)^w$ for $\lambda \in \mfrak{h}^*$ and $w \in \widebar{W}$, we easily obtain from Proposition \ref{prop-tame-twisted} the following statement.
\medskip

\corollary{\label{thm-h.w.-regdom} Let $\lambda \in \mathfrak{h}^*$ be a regular dominant weight. Then $L^\mfrak{g}_{\smash{w(\mfrak{b})}}(w(\lambda))$ is a strongly tame $\Gamma_{\rm st}(w(\Pi))$-Gelfand--Tsetlin $\mfrak{g}$-module for $w\in \widebar{W}$.}

\corollary{\label{cor-diff-Borel}Let $L^\mfrak{g}_\mfrak{b}(\lambda)$ for $\lambda \in \mfrak{h}^*$ be a strongly tame $\Gamma_{\rm st}(\Pi)$-Gelfand--Tsetlin $\mfrak{g}$-module such that $\langle \lambda + \rho_\mfrak{b}, \gamma_i^\vee \rangle \in \N$ for all $i \in I$, where $I$ is some subset of $\{1,2,\dots,n\}$. If $w \in W$ is the element of minimal length satisfying $\langle w(\lambda+\rho_\mfrak{b}), \gamma_i^\vee \rangle \in -\N$ for all $i \in I$, then $L^\mfrak{g}_\mfrak{b}(\lambda) \simeq L^\mfrak{g}_{\smash{w(\mfrak{b})}}(w(\lambda))$ is a strongly tame $\Gamma_{\rm st}(w(\Pi))$-Gelfand--Tsetlin $\mfrak{g}$-module.}

The following observation is a key lemma.
\medskip

\lemma{\label{lem-h.w.-regdom}Let $\lambda \in \mfrak{h}^*$ be a regular dominant weight with respect to $\mfrak{b}$ and let $r \in \{1,2,\dots,n\}$ be such that $\langle \lambda+\rho_\mfrak{b}, \gamma_r^\vee \rangle  \notin \Z$. Let us suppose also that $\lambda$ satisfies one of the following conditions:
\begin{enumerate}[topsep=3pt,itemsep=0pt]
\item[i)] $\langle \lambda+\rho_\mfrak{b}, \gamma_{i,j}^\vee \rangle  \notin \Z$ for all $1\leq i \leq j \leq r$;
\item[ii)]  $\langle \lambda+\rho_\mfrak{b}, \gamma_i^\vee \rangle  \in \N$ for all $i\in \{1,2,\dots, r-1 \}$;
\item[iii)] $\langle\lambda +\rho_\mfrak{b},\gamma_{r-1}^\vee \rangle \notin \mathbb{Z}$, $\langle\lambda+\rho_\mfrak{b},\gamma_{r-1,r}^\vee\rangle \in\mathbb{N}$, and $\langle\lambda+\rho_\mfrak{b},\gamma_k^\vee \rangle \in \mathbb{N}$ for all $k \in \{1,2,\dots,n\} \setminus \{r-1,r\}$ provided $r\neq 1$;
\item[iv)] $\langle\lambda +\rho_\mfrak{b},\gamma_{r+1}^\vee \rangle \notin \mathbb{Z}$, $\langle\lambda+\rho_\mfrak{b},\gamma_{r,r+1}^\vee\rangle \in\mathbb{N}$, and $\langle\lambda+\rho_\mfrak{b},\gamma_k^\vee \rangle \in \mathbb{N}$ for all $k \in \{1,2,\dots,n\} \setminus \{r,r+1\}$ provided $r\neq n$.
\end{enumerate}
Then there exists $w\in W$ such that  $L^\mfrak{g}_{\smash{w^{-1}(\mfrak{b})}}(w^{-1}(\lambda))$ is a strongly tame $\Gamma_{\rm st}(\Pi)$-Gelfand--Tsetlin $\mfrak{g}$-module. Moreover, we have that $L^\mfrak{g}_\mfrak{b}(\lambda)$ is a strongly tame $\Gamma_{\rm st}(w(\Pi))$-Gelfand--Tsetlin $\mfrak{g}$-module.}

\proof{For $r=1$ the statement follows from \cite[Theorem 4.8]{Futorny-Morales-Ramirez2020} with $w=e$. Let us assume that $r>1$. Then in both cases we set $w=s_{\gamma_{r-1}} s_{\gamma_r} s_{\gamma_{r-2}} s_{\gamma_{r-1}} \cdots s_{\gamma_1}s_{\gamma_2}$. Let us choose $v_j \in \C$ for $j=1,2,\dots, n+1$ is such a way that $ v_1-v_2= \langle \lambda + \rho_\mfrak{b}, \gamma_r^\vee \rangle  $,  $v_1- v_{r+1}=-\langle \lambda + \rho_\mfrak{b}, \gamma_{r-1}^\vee \rangle $, $ v_2-v_{r+2}=\langle \lambda + \rho_\mfrak{b}, \gamma_{r+1}^\vee \rangle$, $v_{j+2}-v_{j+3}=\langle \lambda + \rho_\mfrak{b}, \gamma_j^\vee \rangle$ for $1\leq j \leq r-2$, $ v_j-v_{j+1}=\langle \lambda + \rho_\mfrak{b}, \gamma_j^\vee \rangle$ for $r+2 \leq j \leq n$ and
\begin{align*}
  \sum_{j=1}^{n+1}v_j=-\binom{n+1}{2}.
\end{align*}
Further, let us consider the Gelfand--Tsetlin tableau $T(v)$ for $v \in \smash{\C^{{n(n+1) \over 2}}}$ satisfying
\begin{align*}
v_{i,j}=\begin{cases}
v_1+r-1 & \text{if $i=j=1$},\\
v_1+r-i+1& \text{if $1 < i \leq r$ and $j=i-1$},\\
v_1 &  \text{if $i>r$ and $j=r$},\\
v_2+r-i+1 & \text{if $1 < i \leq r$ and $j=i$},\\
v_2 & \text{if $i>r$ and $j=i$},\\
v_3 & \text{if $i\geq 3$ and $j=1$},\\
v_{j+2} & \text{if $i\geq 4$ and $2\leq j <r$},\\
v_j &  \text{if $i>j\geq r+1$}.
\end{cases}
\end{align*}
Let $\mcal{C}$ be the maximal set of relations satisfied by the tableau $T(v)$. Then $\mcal{C}$ is an admissible set of relations and $V_\mcal{C}(T(v))$ is a simple $\mfrak{g}$-module by \cite[Theorem 5.6]{Futorny-Ramirez-Zhang2019}. Therefore, we get that $V_\mcal{C}(T(v))=U(\mfrak{g})T(v)$ together with relations $e_{\gamma_j}(T(v))=0$ for $j\in\{1,2,\dots,n\}\setminus\{2,r+1\}$ and $e_{\gamma_{2,r+1}}(T(v))=0$, $f_{\gamma_{1,r}}(T(v))=0$, which gives us that $T(v)$ is a highest weight vector with respect to the Borel subalgebra $w^{-1}(\mfrak{b})$ of $\mfrak{g}$. Moreover, we have $h_{\gamma_1}(T(v)) = \langle \lambda,\gamma_r^\vee \rangle\,T(v)$, $h_{-\gamma_{1,r}}(T(v))= \langle \lambda,\gamma_{r-1}^\vee \rangle\,T(v)$, $h_{\gamma_{2,r+1}}(T(v)) = \langle \lambda,\gamma_{r+1}^\vee \rangle\,T(v)$, $h_{\gamma_{j+2}}(T(v)) = \langle \lambda,\gamma_j^\vee\rangle \,  T(v)$ for $1\leq j \leq r-2$ and $h_{\gamma_j}(T(v))=\langle \lambda,\gamma_j^\vee \rangle\,T(v)$ for $r+2 \leq j \leq n$, which implies that $V_\mcal{C}(T(v)) \simeq L^\mfrak{g}_{\smash{w^{-1}(b)}}(w^{-1}(\lambda))$.}

\vspace{-2mm}

%%%%%%%%%%%%%%%%%%%%%%%%%%%%%%%%%%%%%%%%%%%%%%%%%%%%%%%%%%%%%%%%%%%%%%%%%%%%%%%%%%%%%%%%%%%

\subsection{Localization of tame Gelfand--Tsetlin modules}
\label{sec:loc functor}

Let $f \in \mfrak{g}$ be a locally $\ad$-nilpotent regular element in $U(\mfrak{g})$. We denote by $U(\mfrak{g})_{(f)}\!$ the left ring of fractions of $U(\mfrak{g})$ with respect to the multiplicative set $\{f^n;\, n \in \N_0\}$ in $U(\mfrak{g})$. In addition, we introduce a $1$-parameter family of algebra automorphisms $\Theta_f^\nu \colon U(\mfrak{g})_{(f)} \rarr U(\mfrak{g})_{(f)}$ by
\begin{align*}
  \Theta_f^\nu(u) = \sum_{k=0}^\infty \binom{\nu+k-1}{k} f^{-k}\ad(f)^k(u)
\end{align*}
for $\nu \in \C$ and $u\in U(\mfrak{g})_{(f)}$. Hence, we may define the \emph{twisted localization functor} $D^\nu_f$ relative to $f$ and $\nu$ on the category of $\mfrak{g}$-modules as
\begin{align*}
  D^\nu_f(M) = U(\mfrak{g})_{(f)} \otimes_{U(\mfrak{g})}\! M
\end{align*}
for a $\mfrak{g}$-module $M$, where the action of $\mfrak{g}$ on $D^\nu_f(M)$ is twisted through $\Theta_f^\nu$, i.e.\ we have
\begin{align*}
  u(v) = \Theta_f^\nu(u)v
\end{align*}
for $u\in U(\mfrak{g})_{(f)}$ and $v \in D_f^\nu(M)$. We will also use the notation $D_f$ instead of $D_f^0$. Furthermore, we define the \emph{twisting functor} $T_f$ on the category of $\mfrak{g}$-modules as
\begin{align*}
  T_f(M) = D_f(M)/M
\end{align*}
for a $\mfrak{g}$-module $M$, which is well defined since $M$ is a $\mfrak{g}$-submodule of $D_f(M)$.
\medskip

Twisted localization functor preserves the annihilator of highest weight modules with dominant highest weights.
\medskip

\lemma{\label{localization annihilator}\cite[Corollary 3.5]{Futorny-Morales-Ramirez2020} Let $\mfrak{b}=\mfrak{b}_\Pi$ for a set of simple roots $\Pi$ in $\Delta$, $\nu \in \C $ and $\lambda\in \mfrak{h}^*$ be a dominant weight with respect to $\mfrak{b}$. Let us assume that $f \in \mfrak{g}$ is a locally $\ad$-nilpotent regular element in $U(\mfrak{g})$. If $f$ is injective on the $\mfrak{g}$-module $L^\mfrak{g}_\mfrak{b}(\lambda)$, then we have $\Ann_{U(\mfrak{g})}\!L^\mfrak{g}_\mfrak{b}(\lambda)= \Ann_{U(\mfrak{g})}\! D_f^\nu (L^\mfrak{g}_\mfrak{b}(\lambda))= \Ann_{U(\mfrak{g})}\!N$, where $N \neq 0$ is a simple subquotient of $\smash{D_f^\nu(L^\mfrak{g}_\mfrak{b}(\lambda))}$.}

In particular, Lemma \ref{localization annihilator} says that if $L^\mfrak{g}_\mfrak{b}(\lambda)$ belongs to a nilpotent orbit $\mcal{O}$ of $\mfrak{g}$, then  $D_f^\nu(L^\mfrak{g}_\mfrak{b}(\lambda))$ belongs to the same nilpotent orbit $\mcal{O}$ for any $\nu \in \C$.
\medskip

Let $\mfrak{b}=\mfrak{b}_\Pi$ for a set of simple roots $\Pi=\{\gamma_1,\gamma_2,\dots,\gamma_n\}$ in $\Delta$. It was shown in \cite{Futorny-Morales-Ramirez2020} that the twisted localization functor $D_f^\nu$ for $\nu \in \C$ defines a functor on the category of strongly tame $\Gamma_{\rm st}(\Pi)$-Gelfand--Tsetlin $\mfrak{g}$-modules via the restriction.
\medskip

\theorem{\label{twistedE21} Let $M$ be a strongly tame $\Gamma_{\rm st}(\Pi)$-Gelfand--Tsetlin $\mfrak{g}$-module and $f=f_{\gamma_1}$.
\begin{enumerate}[topsep=3pt,itemsep=0pt]
\item[i)] \cite[Lemma 5.2]{Futorny-Morales-Ramirez2020} If $f$ is injective on $M$, then $T_f(M)$ is a strongly tame $\Gamma_{\rm st}(\Pi)$-Gelfand--Tsetlin $\mfrak{g}$-module.
\item[ii)] \cite[Theorem 5.4]{Futorny-Morales-Ramirez2020} If $f$ is injective on $M$, then $D_f^\nu(M)$ is a strongly tame $\Gamma_{\rm st}(\Pi)$-Gelfand--Tsetlin $\mfrak{g}$-module for any $\nu \in \C$.
\item[iii)] \cite[Corolary 5.5]{Futorny-Morales-Ramirez2020} If $f$ is bijective on $M$, then $M\simeq D_f^\nu(N)$ for some simple strongly tame $\Gamma_{\rm st}(\Pi)$-Gelfand--Tsetlin $\mfrak{g}$-module $N$ with injective action of $f$ and some $\nu \in \C\setminus \Z$.
\end{enumerate}}
\vspace{-1mm}

\corollary{\label{twistedalpha} Let $M$ be a strongly tame $\Gamma_{\rm st}(\Pi)$-Gelfand--Tsetlin $\mfrak{g}$-module and $f=\smash{f_{w(\gamma_1)}}$ for some $w\in \widebar{W}$. If $f$ is injective on $M^w$, then $D_f^\nu(M^w)$ is a strongly tame $\Gamma_{\rm st}(w(\Pi))$-Gelfand--Tsetlin $\mfrak{g}$-module for any $\nu \in \C$.}

\proof{By Proposition \ref{prop-tame-twisted} we have that $M^w$ for $w \in \widebar{W}$ is a strongly tame  $\Gamma_{\rm st}(w(\Pi))$-Gelfand--Tsetlin $\mfrak{g}$-module. As $\Gamma_{\rm st}(w(\Pi))$ is the standard Gelfand--Tsetlin subalgebra of $U(\mfrak{g})$ with respect to $w(\Pi)$ having $w(\gamma_1)$ as the first simple root and $f$ is injective on $M^w$, the statement follows from Theorem \ref{twistedE21}.}

\corollary{\label{cor-reg-dom-localization} Let $\lambda \in \mfrak{h}^*$ be a regular dominant weight with respect to $\mfrak{b}$ satisfying the assumptions of Lemma \ref{lem-h.w.-regdom} and $f=f_{\gamma_r}$ for $r \in \{1,2,\dots,n\}$. Then
\begin{enumerate}[topsep=3pt,itemsep=0pt]
    \item[i)] $T_f(L^\mfrak{g}_\mfrak{b}(\lambda)) \simeq L^\mfrak{g}_{\smash{s_{\gamma_r}(\mfrak{b})}}(\lambda+\gamma_r)$ is a strongly tame $\Gamma_{\rm st}(w(\Pi))$-Gelfand--Tsetlin $\mfrak{g}$-module;
    \item[ii)] $D_f^\nu(L^\mfrak{g}_\mfrak{b}(\lambda))$ is a
    strongly tame $\Gamma_{\rm st}(w(\Pi))$-Gelfand--Tsetlin $\mfrak{g}$-module for any $\nu \in \C$.
\end{enumerate}}
\vspace{-1mm}

\proof{The first statement follows from Lemma \ref{lem-h.w.-regdom} and Theorem \ref{twistedE21}(i). As $w(\gamma_r)$ is the first simple root of $w(\Pi)$, the second statement follows from Lemma \ref{lem-h.w.-regdom} and Theorem \ref{twistedE21}(ii).}

\vspace{-2mm}

%%%%%%%%%%%%%%%%%%%%%%%%%%%%%%%%%%%%%%%%%%%%%%%%%%%%%%%%%%%%%%%%%%%%%%%%%%%%%%%%%%%%%%%%%%%
%%%%%%%%%%%%%%%%%%%%%%%%%%%%%%%%%%%%%%%%%%%%%%%%%%%%%%%%%%%%%%%%%%%%%%%%%%%%%%%%%%%%%%%%%%%

\section{Explicit realization of admissible modules}

In this section we give an explicit realization of admissible $\mfrak{sl}_{n+1}$-modules.

%%%%%%%%%%%%%%%%%%%%%%%%%%%%%%%%%%%%%%%%%%%%%%%%%%%%%%%%%%%%%%%%%%%%%%%%%%%%%%%%%%%%%%%%%%%

\subsection{Realization of admissible highest weight modules}

Let us remark that strongly tame highest weight $\Gamma_{\rm st}$-Gelfand--Tsetlin $\mfrak{g}$-modules were classified in \cite[Theorem 4.8]{Futorny-Morales-Ramirez2020}. As a consequence of this classification we get the following statement.
\medskip

Let $k\in \Q$ be an admissible number for $\mfrak{g}$ with denominator $q\in \N$, i.e.\ we have
\begin{align*}
  k+n+1 = {p \over q}\ \text{with}\ p,q \in \N,\, (p,q)=1,\, p \geq n+1.
\end{align*}
In the whole section we will assume that $\mfrak{b}$ is a Borel subalgebra of $\mfrak{g}$ associated to a set of simple roots $\Pi=\{\gamma_1,\gamma_2,\dots,\gamma_n\}$ in $\Delta$. The set of positive roots with respect to $\mfrak{b}$ we denote by $\Delta_+^\mfrak{b}$.
\medskip

\theorem{\label{thm-h.w.-adm} For $\lambda \in \smash{\widebar{{\rm Pr}}}_k$, the simple $\mfrak{g}$-module $L^\mfrak{g}_\mfrak{b}(\lambda)$ is a strongly tame $\Gamma_{\rm st}(\Pi)$-Gelfand--Tsetlin $\mfrak{g}$-module. Moreover, if $w\in \widebar{W}$ then $L^\mfrak{g}_{\smash{w(\mfrak{b})}}(w(\lambda))$ is a strongly tame $\Gamma_{\rm st}(w(\Pi))$-Gelfand--Tsetlin $\mfrak{g}$-module.}

\proof{It follows immediately from Corollary \ref{thm-h.w.-regdom}, given that all $\lambda \in \smash{\widebar{{\rm Pr}}}_k$ are regular dominant.}

The following generalizes \cite[Corollary 6.8]{Futorny-Morales-Ramirez2020} for arbitrary nilpotent orbit of $\mfrak{g}$.
\medskip

\corollary{\label{cor-diff-Borel-adm} Let $\lambda \in \smash{\widebar{{\rm Pr}}_k^\mcal{O}}$ for a nilpotent orbit $\mcal{O}$ of $\mfrak{g}$. Let us assume that $f_\gamma$ is injective on $L^\mfrak{g}_\mfrak{b}(\lambda)$ for a simple root $\gamma \in \Pi$. Then $L^\mfrak{g}_{\smash{s_\gamma(\mfrak{b})}}(\lambda+\gamma)$ is an admissible  strongly tame $\Gamma_{\rm st}(s_\gamma(\Pi))$-Gelfand--Tsetlin $\mfrak{g}$-module which belongs to the nilpotent orbit $\mcal{O}$.}

\proof{As $L^\mfrak{g}_{\smash{s_\gamma(\mfrak{b})}}(\lambda+\gamma)$ is a subquotient of $D_{f_\gamma}\!(L^\mfrak{g}_\mfrak{b}(\lambda))$, we get that $L^\mfrak{g}_{\smash{s_\gamma(\mfrak{b})}}(\lambda+\gamma)$ is an admissible $\mfrak{g}$-module which belongs to the nilpotent orbit $\mcal{O}$ by Lemma \ref{localization annihilator}. Applying Theorem \ref{thm-h.w.-adm} we get that $L^\mfrak{g}_{\smash{s_\gamma(\mfrak{b})}}(\lambda+\gamma)$ is a strongly tame $\Gamma_{\rm st}(s_\gamma(\Pi))$-Gelfand--Tsetlin $\mfrak{g}$-module.}

\corollary{Let $\lambda \in \smash{\widebar{{\rm Pr}}}_k$, $1\leq r\leq n$ and
$w=s_{\gamma_{r-1}} s_{\gamma_r} s_{\gamma_{r-2}} s_{\gamma_{r-1}} \cdots s_{\gamma_1}s_{\gamma_2} \in W$ if $r>1$ and $w=e \in W$ if $r=1$. Then $L^\mfrak{g}_\mfrak{b}(\lambda)$ is a strongly tame $\Gamma_{\rm st}(w(\Pi))$-Gelfand--Tsetlin $\mfrak{g}$-module in the following cases:
\begin{enumerate}[topsep=2pt,itemsep=0pt,parsep=0pt]
\item[i)] $\lambda \in \widebar{{\rm Pr}}_k^{\smash{\mcal{O}_{\rm prin}}}$,
\item[ii)] $\lambda \in \Lambda_k(\mfrak{p}^{\rm max}_{\gamma_r}) \subset \widebar{{\rm Pr}}_k^{\smash{\mcal{O}_{[2^{\min\{r,n+1-r\}},1^{n+1-2\min\{r,n+1-r\}}]}}}$ with
    \begin{align*}
      \Lambda_k(\mfrak{p}_{\gamma_r}^{\rm max}) &= \bigg\{\mu-{p\over q}\,a\omega_r;\, \mu \in \widebar{{\rm Pr}}_{k,\Z},\, a \in \N,\, a \leq q-1\bigg\},
   \end{align*}
\item[iii)] $\lambda \in \Lambda_k(\mfrak{p}_{\gamma_n}^{\rm min}) \subset \widebar{{\rm Pr}}_k^{\smash{\mcal{O}_{\rm subreg}}}$ provided $r<n$ with
    \begin{align*}
  \Lambda_k(\mfrak{p}_{\gamma_n}^{\rm min}) = \bigg\{\mu - {p \over q} \sum_{i=1}^{n-1} a_i\omega_i;\, \mu \in \widebar{{\rm Pr}}_{k,\Z},\ a_1,a_2,\dots, a_{n-1} \in \N,\, \sum_{i=1}^{n-1} a_i \leq q-1\bigg\},
  \end{align*}
\end{enumerate}
where $k+n+1={p \over q}$ with $p,q \in \N$, $(p,q)=1$ and $p \geq n+1$.}

\vspace{-2mm}

%%%%%%%%%%%%%%%%%%%%%%%%%%%%%%%%%%%%%%%%%%%%%%%%%%%%%%%%%%%%%%%%%%%%%%%%%%%%%%%%%%%%%%%%%%%

\subsection{Realization of admissible $\mfrak{sl}_2$-induced modules}

We will discuss tableau realizations of simple admissible $\mfrak{sl}_2$-induced $\mfrak{g}$-modules which belong to the principal nilpotent orbit, the subregular nilpotent orbit and the nilpotent orbits for maximal parabolic subalgebras.
\medskip

Let us consider a parabolic subalgebra $\mfrak{p} = \mfrak{l} \oplus \mfrak{u} \supset
\mfrak{b}$ of $\mfrak{g}$ such that the Levi subalgebra $\mfrak{l}$ is isomorphic to $\mfrak{sl}_2 + \mfrak{h}$. For a simple weight $\mfrak{l}$-module $E$, we may consider the generalized Verma module $M^\mfrak{g}_\mfrak{p}(E)$, where $E$ is understood as a $\mfrak{p}$-module on which $\mfrak{u}$ acts trivially, together with its simple quotient $L^\mfrak{g}_\mfrak{p}(E)$. We will call them \emph{$\mfrak{sl}_2$-induced} $\mfrak{g}$-modules. Since $\mfrak{p}$ is a standard parabolic subalgebra of $\mfrak{g}$, there exists a simple root $\gamma \in \Pi = \{\gamma_1,\gamma_2,\dots,\gamma_n\}$ such that the corresponding $\mfrak{sl}_2$-triple $(e_\gamma,h_\gamma,f_\gamma)$ forms a basis of $[\mfrak{l},\mfrak{l}]$. Based on the well-known classification of simple weight $\mfrak{sl}_2$-modules (see e.g.\ \cite{Mazorchuk2010-book}), we easily obtain a classification of simple weight $\mfrak{l}$-modules which looks as follows.
\begin{enumerate}[topsep=2pt,itemsep=3pt,parsep=0pt]
  \item[i)] The finite-dimensional $\mfrak{l}$-modules $E_\lambda$ with highest weight $\lambda \in \mfrak{h}^*$, where $\lambda_\gamma \in \N_0$.
  \item[ii)] The highest weight $\mfrak{l}$-modules $L_\lambda$ with highest weight $\lambda \in \mfrak{h}^*$, where $\lambda_\gamma \in \C \setminus \N_0$.
  \item[iii)] The lowest weight $\mfrak{l}$-modules $\widebar{L}_\lambda$ with lowest weight $\lambda \in \mfrak{h}^*$, where $\lambda_\gamma \in \C\setminus\! -\N_0$.
  \item[iv)] The cuspidal $\mfrak{l}$-modules $R_{\lambda,\xi}$ with weight $\lambda \in \mfrak{h}^*$ and $\xi \in \C$, where $\lambda_\gamma \in \C$ and $2\xi \neq \mu(\mu+2)$ for all $\mu \in \lambda_\gamma + 2\Z$.
\end{enumerate}
We used the notation $\lambda_\gamma = \lambda(h_\gamma)$ for $\lambda \in \mfrak{h}^*$. Let us note that the parameter $\xi \in \C$ for $R_{\lambda,\xi}$ is the eigenvalue of the $\mfrak{sl}_2$-Casimir element $c_\gamma$ defined by
\begin{align*}
  c_\gamma = e_\gamma f_\gamma + f_\gamma e_\gamma + {\textstyle {1\over 2}}h^2_\gamma.
\end{align*}
Besides, we have that $R_{\lambda,\xi} \simeq R_{\lambda+n\gamma,\xi}$ for $n \in \Z$. Also note that all $\mfrak l$-modules above are strongly tame
$\Gamma_{\rm st}$-Gelfand--Tsetlin modules.

In the next, it will be convenient for us to use the following realization of the simple weight $\mfrak{l}$-modules from the classification above, i.e.\ we have
\begin{enumerate}[topsep=2pt,itemsep=3pt,parsep=0pt]
  \item[i)] $\smash{L^\mfrak{l}_{\mfrak{b} \cap \mfrak{l}}(\lambda)} \simeq E_\lambda$ for $\lambda \in \mfrak{h}^*$ if $\langle \lambda+\rho_\mfrak{b},
      \gamma^\vee \rangle \in \N$;
  \item[ii)] $\smash{M^\mfrak{l}_{\mfrak{b} \cap \mfrak{l}}(\lambda)} \simeq L_\lambda$ for $\lambda \in \mfrak{h}^*$ if $\langle \lambda+\rho_\mfrak{b},
      \gamma^\vee \rangle \notin \N$;
  \item[iii)]$T_{f_\gamma}\!(\smash{M^\mfrak{l}_{\mfrak{b} \cap \mfrak{l}}(\lambda)}) \simeq \widebar{L}_{\lambda+\gamma}$ for $\lambda \in \mfrak{h}^*$ if $\langle \lambda+\rho_\mfrak{b},
      \gamma^\vee \rangle \notin -\N$;
  \item[iv)] $D^\nu_{f_\gamma}\!(\smash{M^\mfrak{l}_{\mfrak{b} \cap \mfrak{l}}(\lambda)}) \simeq R_{\lambda+\nu\gamma,\xi_\lambda}$ for $\lambda \in \mfrak{h}^*$, $\nu \in \C \setminus \Z$ if $\nu + \langle \lambda+\rho_\mfrak{b}, \gamma^\vee \rangle \notin \Z$, where $\xi_\lambda = \smash{{\lambda_\gamma(\lambda_\gamma+2) \over 2}}$.
\end{enumerate}
\medskip

Since $L^\mfrak{g}_\mfrak{p}(E_\lambda) \simeq L^\mfrak{g}_\mfrak{b}(\lambda)$ if $\langle \lambda+\rho_\mfrak{b}, \gamma^\vee \rangle \in \N$, $L^\mfrak{g}_\mfrak{p}(L_\lambda) \simeq L^\mfrak{g}_\mfrak{b}(\lambda)$ if $\langle \lambda+\rho_\mfrak{b}, \gamma^\vee \rangle \notin \N$, and $L^\mfrak{g}_\mfrak{p}(\widebar{L}_\lambda) \simeq L^\mfrak{g}_{\smash{s_\gamma(\mfrak{b})}}(\lambda+\gamma)$ if $\langle \lambda+\rho_\mfrak{b}, \gamma^\vee \rangle \notin -\N$, we see that in cases (i)--(iii) the corresponding simple $\mfrak{sl}_2$-induced $\mfrak{g}$-modules are highest weight $\mfrak{g}$-modules. On the other hand, the last case leads to non-highest weight $\mfrak{g}$-modules that we will focus on below.
Let us note that the action of $e_\gamma$ and $f_\gamma$ on $R_{\lambda+\nu\gamma,\xi_\lambda}$ is injective for $\lambda \in \mfrak{h}^*$ and $\nu \in \C \setminus \Z$ provided $\nu + \langle \lambda+\rho_\mfrak{b}, \gamma^\vee \rangle \notin \Z$.
\medskip

\proposition{\cite[Theorem 5.6]{Futorny-Morales-Ramirez2020}\label{sl2 Induced} Let $n>1$, $\lambda \in \mfrak{h}^*$ and $\nu \in \C \setminus \Z$. Then $L^\mfrak{g}_\mfrak{p}(R_{\lambda+\nu\gamma_1,\xi_\lambda})$ is a strongly tame $\Gamma_{\rm st}(\Pi)$-Gelfand--Tsetlin $\mfrak{g}$-module if and only if $L^\mfrak{g}_\mfrak{b}(\lambda)$ is a strongly tame $\Gamma_{\rm st}(\Pi)$-Gelfand--Tsetlin module provided $\langle \lambda+\rho_\mfrak{b},\gamma^\vee_1 \rangle  \notin \N_0$ and $\nu+\langle \lambda+\rho_\mfrak{b},\gamma^\vee_1 \rangle  \notin\Z$. Moreover, we have that $L^\mfrak{g}_\mfrak{p}(R_{\lambda+\nu\gamma_1,\xi_\lambda}) \simeq \smash{D_{f_{\gamma_1}}^\nu\!(L^\mfrak{g}_\mfrak{b}(\lambda))}$.}

Let $k\in \Q$ be an admissible number for $\mfrak{g}$ with denominator $q\in \N$, i.e.\ we have
\begin{align*}
  k+n+1 = {p \over q}\ \text{with}\ p,q \in \N,\, (p,q)=1,\, p \geq n+1,
\end{align*}
and let $E$ be a simple weight $\mfrak{l}$-module. It follows from \cite[Theorem 2.12]{Arakawa-Futorny-Ramirez2017} that if $L^\mfrak{g}_\mfrak{p}(E)$ is an admissible $\mfrak{g}$-module of level $k$, then $E$ is an admissible $\mfrak{sl}_2$-module of level $k_\gamma = k+n-1$. From the realization of $R_{\lambda+\nu\gamma,\xi_\lambda}$ through the twisted localization functor $\smash{D^\nu_{f_\gamma}}$ and by using Lemma \ref{localization annihilator}, we get the following statement.
\medskip

\proposition{\label{prop-sl2}\cite[Proposition 6.9]{Futorny-Morales-Ramirez2020} Let $\lambda \in \mfrak{h}^*$ and $\nu \in \C \setminus \Z$. The $\mfrak{l}$-module $R_{\lambda+\nu\gamma,\xi_\lambda}$ is an admissible $\mfrak{sl}_2$-module of level $k_\gamma$ if and only if $\langle \lambda+\rho_\mfrak{b}, \gamma^\vee \rangle = b - {p \over q}a$ and $\nu - {p\over q}a \notin \Z$, where $a,b \in \N$ satisfy $a\leq q-1$ and $b \leq p-2$.}

Now by applying Proposition \ref{sl2 Induced} and Proposition \ref{prop-sl2} we obtain the following generalization of \cite[Theorem 6.5]{Arakawa-Futorny-Ramirez2017} and \cite[Corollary 6.11]{Futorny-Morales-Ramirez2020}.
\medskip

\corollary{\label{admindsl2} Let $n>1$, $\lambda \in \mfrak{h}^*$ and $\nu \in \C \setminus \Z$. Then $L^\mfrak{g}_\mfrak{p}(R_{\lambda+\nu\gamma_1,\xi_\lambda})$ is an admissible strongly tame $\Gamma_{\rm st}(\Pi)$-Gelfand--Tsetlin $\mfrak{g}$-module if and only if $\lambda \in \smash{\widebar{{\rm Pr}}_k}$ provided $\langle \lambda+\rho_\mfrak{b},\gamma^\vee_1 \rangle  \notin \Z$ and $\nu+\langle \lambda+\rho_\mfrak{b},\gamma^\vee_1 \rangle  \notin\Z$.}

Let us note that admissible $\mfrak{g}$-modules $L^\mfrak{g}_\mfrak{p}(R_{\lambda+\nu\gamma_1,\xi_\lambda})$ in Corollary \ref{admindsl2} are $\mfrak{sl}_2$-induced modules, where the Lie subalgebra $\mfrak{sl}_2$ of $\mfrak{g}$ is given by the first simple root $\gamma_1$. Hence, we need an extension of Corollary \ref{admindsl2} to other simple roots. As a consequence of Proposition \ref{sl2 Induced} and Corollary \ref{cor-reg-dom-localization}(ii) we obtain the following result.
\medskip

\corollary{\label{cor-local-r} Let $n>1$, $\lambda\in \mfrak{h}^*$ and $\nu \in \C \setminus \Z$. Let us assume that $\lambda$ satisfies the conditions of Lemma \ref{lem-h.w.-regdom} for a simple root $\gamma$. Then $\smash{D^\nu_{f_{\gamma}}\!(L^\mfrak{g}_\mfrak{b}(\lambda))}$ is a simple strongly tame $\Gamma_{\rm st}(w(\Pi))$-Gelfand--Tsetlin $\mfrak{g}$-module provided $\nu + \langle \lambda+\rho, \gamma^\vee \rangle \notin \Z$. Besides, we have $\smash{D^\nu_{f_{\gamma}}\!(L^\mfrak{g}_\mfrak{b}(\lambda))} \simeq L^\mfrak{g}_\mfrak{p}(R_{\lambda+\nu\gamma,\xi_\lambda})$.}

\vspace{-2mm}

%%%%%%%%%%%%%%%%%%%%%%%%%%%%%%%%%%%%%%%%%%%%%%%%%%%%%%%%%%%%%%%%%%%%%%%%%%%%%%%%%%%%%%%%%%%

\subsubsection{Principal nilpotent orbit}

By Theorem \ref{thm-prin} we have that $[\widebar{{\rm Pr}}_k^{\smash{\mcal{O}_{\rm prin}}}]$ is represented by the set $\Lambda_k(\mfrak{b})$, where
\begin{align*}
  \Lambda_k(\mfrak{b}) = \bigg\{\mu - {p \over q} \sum_{i=1}^n a_i\omega_i;\, \mu \in \widebar{{\rm Pr}}_{k,\Z},\ a_1,a_2,\dots,a_n \in \N,\, \sum_{i=1}^n a_i \leq q-1\bigg\}.
\end{align*}
Further, since for $\lambda \in \Lambda_k(\mfrak{b})$, $\langle \lambda+\rho_\mfrak{b},\gamma^\vee \rangle \notin \Z$ for all $\gamma \in \Delta_+^\mfrak{b}$, we obtain immediately that a weight $\lambda \in \widebar{{\rm Pr}}_k^{\smash{\mcal{O}_{\rm prin}}}$ is not only regular dominant but also antidominant. Therefore, we have the following refinement of Theorem \ref{thm-h.w.-adm} for the principal nilpotent orbit.
\medskip

\theorem{\label{prop-diff-hw} If $\lambda \in \widebar{{\rm Pr}}_k^{\smash{\mcal{O}_{\rm prin}}}$, then $L^\mfrak{g}_\mfrak{b}(\lambda)$ is a strongly tame $\Gamma_{\rm st}(w(\Pi))$-Gelfand--Tsetlin $\mfrak{g}$-module for any $w\in W$.}

\proof{If $\lambda \in \widebar{{\rm Pr}}_k^{\smash{\mcal{O}_{\rm prin}}}$, then the weight
$\lambda$ is regular dominant and antidominant, and $L^\mfrak{g}_\mfrak{b}(\lambda)\simeq M^\mfrak{g}_\mfrak{b}(\lambda)$ is a strongly tame $\Gamma_{\rm st}(\Pi)$-Gelfand--Tsetlin $\mfrak{g}$-module. Further, for $w\in W$ the weight $w^{-1}(\lambda)$ is also strongly generic, and hence $L^\mfrak{g}_{\smash{w^{-1}(\mfrak{b})}}(w^{-1}(\lambda))$ is a strongly tame $\Gamma_{\rm st}(\Pi)$-Gelfand--Tsetlin $\mfrak{g}$-module, which gives us that
$L^\mfrak{g}_{\smash{w^{-1}(\mfrak{b})}}(w^{-1}(\lambda))^w\simeq L^\mfrak{g}_\mfrak{b}(\lambda)$ is a strongly tame $\Gamma_{\rm st}(w(\Pi))$-Gelfand--Tsetlin $\mfrak{g}$-module.}

We have the following  refinement and generalization of \cite[Theorem 6.5]{Arakawa-Futorny-Ramirez2017}. It provides a classification and an explicit construction of all simple admissible $\mfrak{sl}_2$-induced modules in the principal nilpotent orbit.
\medskip

\theorem{\label{principaltwisted}Let  $n>1$, $\lambda \in \widebar{{\rm Pr}}_k^{\smash{\mcal{O}_{\rm prin}}}$, $\gamma \in \Delta_+^\mfrak{b}$, $\nu \in \C \setminus \Z$ and $\nu + \langle \lambda+\rho_\mfrak{b},\gamma^\vee \rangle \notin \Z$. Further, let us assume that $w \in \smash{\widebar{W}}$ satisfies $\gamma = w(\gamma_1)$. Then
\begin{enumerate}[topsep=3pt,itemsep=0pt]
 \item[i)] $D^\nu_{f_\gamma}\!(L^\mfrak{g}_\mfrak{b}(\lambda))$ is a simple admissible strongly tame $\Gamma_{\rm st}(w(\Pi))$-Gelfand--Tsetlin $\mfrak{g}$-module which belongs to the principal nilpotent orbit;
  \item[ii)] $\mfrak{g}$-modules $D^\nu_{f_\gamma}\!(L^\mfrak{g}_\mfrak{b}(\lambda))$ for $\gamma \in \Pi$ exhaust all simple admissible $\mfrak{sl}_2$-induced $\mfrak{g}$-modules which belongs to the principal nilpotent orbit. All such $\mfrak{g}$-modules have unbounded finite weight multiplicities.
\end{enumerate}}

\proof{By Theorem \ref{prop-diff-hw} we have that $L^\mfrak{g}_\mfrak{b}(\lambda)$ is a strongly tame $\Gamma_{\rm st}(\Pi)$-Gelfand--Tsetlin $\mfrak{g}$-module. Since $f_\gamma$ is injective on $L^\mfrak{g}_\mfrak{b}(\lambda)$, the first  statement follows from Lemma \ref{localization annihilator} and Theorem \ref{twistedE21}. The second statement follows from Proposition \ref{sl2 Induced}.}

\remark{\begin{enumerate}[topsep=0pt,itemsep=0pt,parsep=0pt]
\item[i)] Let us note that $D^\nu_{f_\gamma}\!(L^\mfrak{g}_\mfrak{b}(\lambda))$ has infinite weight multiplicities if $\gamma \notin \Pi$.
\item[ii)] We see that admissible modules in the principal nilpotent orbit admit various realizations as strongly tame Gelfand--Tsetlin modules analogously to finite-dimensional modules. This is typical for admissible modules in all nilpotent orbits.
\end{enumerate}}

%%%%%%%%%%%%%%%%%%%%%%%%%%%%%%%%%%%%%%%%%%%%%%%%%%%%%%%%%%%%%%%%%%%%%%%%%%%%%%%%%%%%%%%%%%%

\subsubsection{Nilpotent orbits for maximal parabolic subalgebras}

Let us consider $1 \leq r \leq n$. Then by the results of Section \ref{subsec:small denominator} we have that the set of equivalence classes $[\widebar{{\rm Pr}}_k^{\smash{\mcal{O}_{[2^{\min\{r,n+1-r\}},1^{n+1-2\min\{r,n+1-r\}}]}}}]$ is represented by the set $\Lambda_k(\mfrak{p}_{\gamma_r}^{\rm max})$, where
\begin{align*}
      \Lambda_k(\mfrak{p}_{\gamma_r}^{\rm max}) &= \bigg\{\mu-{p\over q}\,a\omega_r;\, \mu \in \widebar{{\rm Pr}}_{k,\Z},\, a \in \N,\, a \leq q-1\bigg\}.
\end{align*}
Let us note that for $\lambda \in \Lambda_k(\mfrak{p}_{\gamma_r}^{\rm max})$ with $r=1$ or $r=n$, the $\mfrak{g}$-module $L^\mfrak{g}_\mfrak{b}(\lambda)$ is a simple admissible highest weight $\mfrak{g}$-module which belongs to the the minimal nilpotent orbit.
\medskip

Next theorem provides a realization for certain simple admissible $\mfrak{sl}_2$-induced  $\mfrak{g}$-modules which belong to the nilpotent orbit attached to a maximal parabolic subalgebra of $\mfrak{g}$.
\medskip

\theorem{\label{Pro:maxparabolic}Let $\lambda \in \Lambda_k(\mfrak{p}^{\rm max}_{\gamma_r})$, $1\leq r\leq n$, $\gamma \in \Delta_+^\mfrak{b}$, $\nu \in \C \setminus \Z$ and $\nu+\langle \lambda+\rho_\mfrak{b},\gamma^\vee \rangle \notin\Z$. Further, let us assume that $\langle \lambda + \rho_\mfrak{b},\gamma^\vee \rangle \notin\Z$. Then $D^\nu_{f_\gamma}\!(L^\mfrak{g}_\mfrak{b}(\lambda))$ is an admissible $\mfrak{g}$-module which belongs to the nilpotent orbit attached to the the maximal parabolic subalgebra $\smash{\mfrak{p}_{\gamma_r}^{\rm max}}$. Moreover, there exists $w\in \widebar{W}$ such that $D^\nu_{f_\gamma}\!(L^\mfrak{g}_\mfrak{b}(\lambda))$ is a strongly tame $\Gamma_{\rm st}(w(\Pi))$-Gelfand--Tsetlin $\mfrak{g}$-module.}

\proof{If $r=1$ or $r=n$ then the statement was proved in \cite[Theorem 6.12]{Futorny-Morales-Ramirez2020}, where a complete classification of simple admissible $\mfrak{sl}_2$-induced $\mfrak{g}$-modules in the minimal orbit was obtained. Let us assume now that $1<r<n$. The admissibility of $D^\nu_{f_\gamma}\!(L^\mfrak{g}_\mfrak{b}(\lambda))$ follows from Lemma \ref{localization annihilator}. If we have $\gamma=\gamma_r$ for $1< r \leq \smash{\lfloor \frac{n+1}{2} \rfloor}$, then the statement follows from Corollary \ref{cor-reg-dom-localization} for $w_r=s_{\gamma_{r-1}} s_{\gamma_r} s_{\gamma_{r-2}} s_{\gamma_{r-1}} \cdots s_{\gamma_1}s_{\gamma_2} \in W$. If $r> \smash{\lfloor \frac{n+1}{2} \rfloor}$ we apply the automorphism $\sigma$ and then the element $w_{\gamma_{n-r+1}}$. This proves the statement for simple roots.

Let us assume now that $\gamma=\gamma_{i,j}$ for $1\leq i \leq r \leq j \leq n$ and $i < j$. We proceed by induction on $j-i$ with $j-i=0$ as the base of induction which corresponds to a simple root considered above. For $i<r$, let us consider the $\mfrak{g}$-module $L^\mfrak{g}_{\smash{s_{\gamma_i}\!(\mfrak{b})}}(s_i(\lambda))\simeq L^\mfrak{g}_\mfrak{b}(\lambda)$. Then $L^\mfrak{g}_{\smash{s_{\gamma_i}\!(\mfrak{b})}} (s_i(\lambda))$ is a $\Gamma_{\rm st}(s_{\gamma_i}\!(\Pi))$-Gelfand--Tsetlin module (cf.\ Corollary \ref{cor-diff-Borel}). The Borel subalgebra $s_{\gamma_i}\!(\mfrak{b})$ corresponds to a new set of simple roots $s_{\gamma_i}(\Pi)=\{\gamma'_1,\gamma'_2,\dots,\gamma'_n\}$ and $\gamma=\gamma'_{i+1,j}$. The result then follows by induction. Further, for $i=r$ we have $j>r$. In this case we consider the $\mfrak{g}$-module $L^\mfrak{g}_{\smash{s_{\gamma_j}\!(\mfrak{b})}}(s_{\gamma_j}(\lambda))\simeq L^\mfrak{g}_\mfrak{b}(\lambda)$ and we apply the same argument. This completes the proof for $\lambda \in \Lambda_k(\mfrak{p}^{\rm max}_{\gamma_r})$.}

\remark{Theorem \ref{Pro:maxparabolic} shows how to explicitly construct simple admissible $\mfrak{sl}_2$-induced $\mfrak{g}$-modules for the representatives of maximal parabolic orbits. These modules do not exhaust all admissible $\mfrak{sl}_2$-induced $\mfrak{g}$-modules. For example, by acting on $\lambda$ by a simple reflection $s_\gamma$ corresponding to a non-integral entry one obtains an admissible highest weight $\mfrak{g}$-module $L^\mfrak{g}_\mfrak{b}(s_\gamma\cdot \lambda)$ in a maximal parabolic orbit which can be localized with respect to any root $\gamma_{i,j}$ such that $\langle s_\gamma \cdot \lambda, \gamma_{i,j}^\vee \rangle \notin \Z$. This follows from Corollary \ref{cor-diff-Borel-adm} and Lemma \ref{lem-h.w.-regdom}.  On the other hand, more admissible $\mfrak{sl}_2$-induced $\mfrak{g}$-modules can be obtained by using other highest weights in maximal parabolic orbits.}

%%%%%%%%%%%%%%%%%%%%%%%%%%%%%%%%%%%%%%%%%%%%%%%%%%%%%%%%%%%%%%%%%%%%%%%%%%%%%%%%%%%%%%%%%%%

\subsubsection{Subregular nilpotent orbit}

By Theorem \ref{thm:subregular orbit} we have that $[\widebar{{\rm Pr}}_k^{\smash{\mcal{O}_{\rm subreg}}}]$ is represented by the set $\Lambda_k(\mfrak{p}_{\gamma_n}^{\rm min})$, where
\begin{align*}
  \Lambda_k(\mfrak{p}_{\gamma_n}^{\rm min}) = \bigg\{\mu - {p \over q} \sum_{i=1}^{n-1} a_i\omega_i;\, \mu \in \widebar{{\rm Pr}}_{k,\Z},\ a_1,a_2,\dots, a_{n-1} \in \N,\, \sum_{i=1}^{n-1} a_i \leq q-1\bigg\}.
\end{align*}
Let us note that for $\lambda\in  \Lambda_k(\mfrak{p}_{\gamma_n}^{\rm min})$ and
all $\gamma \in \Delta_+^\mfrak{b} \setminus \{\gamma_n\}$, we have $\langle \lambda+\rho_\mfrak{b},\gamma^\vee \rangle \notin \Z$.
\medskip

As an immediate consequence of Lemma \ref{lem-h.w.-regdom} and \cite[Theorem 5.4]{Futorny-Morales-Ramirez2020} we obtain the following statement. For $1\leq r\leq n$, we set $w_r=s_{\gamma_{r-1}} s_{\gamma_r} s_{\gamma_{r-2}} s_{\gamma_{r-1}} \cdots s_{\gamma_1} s_{\gamma_2} \in W$ if $r>1$, and $w_r=e$ for $r=1$.
\medskip

\corollary{\label{cor-h.w.-subreg}
Let $\lambda \in \Lambda_k(\mfrak{p}_{\gamma_n}^{\rm min})$, $1\leq r\leq n$ and $w=w_r$. Then
\begin{enumerate}[topsep=3pt,itemsep=0pt]
\item[i)] $L^\mfrak{g}_{\smash{w^{-1}(\mfrak{b})}}(w^{-1}(\lambda))$ is a strongly tame $\Gamma_{\rm st}(\Pi)$-Gelfand--Tsetlin $\mfrak{g}$-module;
\item[ii)] $L^\mfrak{g}_\mfrak{b}(\lambda)$ is a strongly tame $\Gamma_{\rm st}(w(\Pi))$-Gelfand--Tsetlin $\mfrak{g}$-module.
\end{enumerate}
}

The next theorem follows from Corollary \ref{cor-local-r}. It provides a realization for certain
simple admissible $\mfrak{sl}_2$-induced $\mfrak{g}$-modules which belong to the subregular nilpotent orbit.
\medskip

\theorem{\label{thm-subreg}Let $\lambda \in \Lambda_k(\mfrak{p}_{\gamma_n}^{\rm min})$, $\gamma\in \{\gamma_1,\dots,\gamma_{n-1}, \gamma_{n-1}+\gamma_n\}$, $\nu \in \C \setminus \Z$, $\nu+\langle \lambda+\rho_\mfrak{b}, \gamma^\vee \rangle  \notin\Z$. Then $\smash{D^\nu_{f_\gamma}}\!(L^\mfrak{g}_\mfrak{b}(\lambda))$ is a simple admissible $\mfrak{g}$-module which belongs to the subregular nilpotent orbit. Moreover, we have that $\smash{D^\nu_{f_\gamma}}\!(L^\mfrak{g}_\mfrak{b}(\lambda))$ is a strongly tame $\Gamma_{\rm st}(w(\Pi))$-Gelfand--Tsetlin $\mfrak{g}$-module with $w=w_r$ for $\gamma=\gamma_r$, $1\leq r<n$,  and $w=w_{n-1}s_{\gamma_n}$ for $\gamma=\gamma_{n-1}+\gamma_n$.}

Let us note that the case $\gamma=\gamma_1$ in Theorem \ref{thm-subreg} is already covered by Theorem \ref{twistedE21}.

%%%%%%%%%%%%%%%%%%%%%%%%%%%%%%%%%%%%%%%%%%%%%%%%%%%%%%%%%%%%%%%%%%%%%%%%%%%%%%%%%%%%%%%%%%%
%%%%%%%%%%%%%%%%%%%%%%%%%%%%%%%%%%%%%%%%%%%%%%%%%%%%%%%%%%%%%%%%%%%%%%%%%%%%%%%%%%%%%%%%%%%

\section{Tableau realization of admissible $\mfrak{sl}_4$-modules}
\label{sl4-modules}

In this section we give examples of a tableau realization of simple admissible modules over the Lie algebra $\mfrak{sl}_4$.
\medskip

Let $\mfrak{g}=\mfrak{sl}_4$ and let $k\in \Q$ be an admissible number for $\mfrak{g}$ with denominator $q \in \N$, i.e.\ we have
\begin{align*}
  k+4= {p \over q}\ \text{with}\ p,q \in \N,\, (p,q)=1,\, p \geq 4.
\end{align*}
By results of Section \ref{subsec:lie algebra sl4} we know that if $\lambda \in \widebar{{\rm Pr}}_k$, then the simple highest weight $\mfrak{g}$-module $L^\mfrak{g}_\mfrak{b}(\lambda)$ belongs to the one of the following nilpotent orbits \begin{align*}
  \mcal{O}_{\rm prin}=\mcal{O}_\mfrak{b}, \qquad \mcal{O}_{\rm subreg}=\mcal{O}_{\mfrak{p}_{\alpha_3}^{\rm min}}, \qquad \mcal{O}_{\rm rect} = \mcal{O}_{\mfrak{p}_{\alpha_2}^{\rm max}}, \qquad \mcal{O}_{\rm min}=\mcal{O}_{\mfrak{p}_{\alpha_1}^{\rm max}}, \qquad \mcal{O}_{\rm zero}=\mcal{O}_\mfrak{g}
\end{align*}
depending on the level $k$. We will restrict ourselves to the case when $\lambda \in \Lambda_k(\mfrak{p})$ for the corresponding nilpotent orbit $\mcal{O}_\mfrak{p}$. For all nilpotent orbits, except the zero nilpotent orbit, we will describe a tableau realization of $\mfrak{g}$-modules $M^\mfrak{g}_\mfrak{p}(\lambda)$, $T_fM^\mfrak{g}_\mfrak{p}(\lambda)$ and $D^\nu_fM^\mfrak{g}_\mfrak{p}(\lambda)$ for $f = f_{\alpha_{1,2}}$ and $\nu \in \C$.

For that reason, let us consider the Gelfand--Tsetlin subalgebra $\Gamma=\Gamma_{\rm st}(\Pi_w)$ for the element $w=s_2$ of the Weyl group $W$. The corresponding chain $\mcal{F}$ of Lie subalgebras of $\mfrak{g}$ is determined by the sequence of sets of their simple roots
\begin{align*}
\{\alpha_{1,2}\} \subset \{\alpha_{1,2},-\alpha_2\} \subset \{\alpha_{1,2},-\alpha_2, \alpha_{2,3}\}.
\end{align*}
By applying a representative $\dot{w}$ of $w$ to the formulas \eqref{Gelfand--Tsetlin formulas} we obtain the action of $\mfrak{g}$ on strongly tame $\Gamma$-Gelfand--Tsetlin $\mfrak{g}$-modules.

For a description of bases we will use the following notation introduced in \cite{Futorny-Ramirez-Zhang2019}. Let $T(v)$ for $v \in \C^{10}$ be a tableau satisfying a set of relations in $\mcal{C}$.
We denote by
\begin{align*}
  \mcal{B}_\mcal{C}(T(v)) = \{T(v+i(z));\, z \in \Z^6,\, T(v+i(z)) \text{ satisfies relations in } \mcal{C}\},
\end{align*}
where the embedding $i \colon \C^6 \rarr \C^{10}$ is given by
\begin{align*}
  (v_{3,1},v_{3,2},v_{3,3}|v_{2,1},v_{2,2}|v_{1,1}) \mapsto (0,0,0,0|v_{3,1},v_{3,2},v_{3,3}|v_{2,1},v_{2,2}|v_{1,1}).
\end{align*}
Besides, since the condition that a tableau $T(v+i(z))$ satisfies relations in $\mcal{C}$ leads to a system of inequalities on
\begin{align*}
  z = (r,s,t|m,n|\ell)
\end{align*}
with integer coefficients, we shall denote the set $\mcal{B}_\mcal{C}(T(v))$ by these inequalities. Furthermore, we denote by $V_\mcal{C}(T(v))$ the $\Gamma$-Gelfand--Tsetlin $\mfrak{g}$-module spanned by $\mcal{B}_\mcal{C}(T(v))$ with the corresponding action of $\mfrak{g}$.

%%%%%%%%%%%%%%%%%%%%%%%%%%%%%%%%%%%%%%%%%%%%%%%%%%%%%%%%%%%%%%%%%%%%%%%%%%%%%%%%%%%%%%%%%%%

\newpage

\subsection{Principal nilpotent orbit}

For $\lambda \in \Lambda_k(\mfrak{b}) \subset \widebar{{\rm Pr}}_k^{\smash{\mcal{O}_{\rm prin}}}$, we have
\begin{align*}
 \lambda=\big(\lambda_1- {\textstyle {p \over q}}a,\lambda_2-{\textstyle {p \over q}}b,\lambda_3-{\textstyle {p \over q}}a\big)\ \text{with}\ a,b,c \in \N,\ a+b+c \leq q-1.
\end{align*}
We set $v_1-v_3=\langle \lambda+\rho,\alpha_1^\vee \rangle \notin\Z$, $v_3-v_2=\langle \lambda+\rho,\alpha_2^\vee \rangle \notin\Z$,  $v_2-v_4=\langle \lambda+\rho,\alpha_3^\vee \rangle \notin\Z$ such that $v_1+v_2+v_3+v_4=-6$.
\begin{align*}
T(v) = \begin{aligned}
\begin{tikzpicture}[x=15mm, y=7mm, every node/.style={draw, minimum width=14.3mm, minimum height=6.3mm}]
    \node at (0,0) {$v_1$};
    \node at (1,0) {$v_2$};
    \node at (2,0) {$v_3$};
    \node at (3,0) {$v_4$};
    \node at (0.5,-1) {$v_1$};
    \node at (1.5,-1) {$v_2$};
    \node at (2.5,-1) {$v_3$};
    \node at (1,-2) {$v_1$};
    \node at (2,-2) {$v_2+1$};
    \node at (1.5,-3) {$v_1$};
\end{tikzpicture}
\end{aligned}
\end{align*}

\vspace{-5mm}

\begin{table}[h]
\centering
\renewcommand{\arraystretch}{2.5}
\begin{tabular}{|c|c|c|c|}
    \hline
    $V_\mcal{C}(T(w))$ & $T(w)$ & $\mcal{C}$ & $\mcal{B}_\mcal{C}(T(w))$ \\
    \hline
    $M^\mfrak{g}_\mfrak{p}(\lambda)$ & $T(v)$ &
    \begin{minipage}{0.31\textwidth}
    \vspace{2mm}
    \begin{tikzpicture}[x=7mm, y=7mm, >=latex, yscale=1, xscale=2, every node/.style={circle, draw=white, fill=black, line width=1.5mm, inner sep=1mm}]
      \node (A) at (0,0) {};
      \node (B) at (1,0) {};
      \node (C) at (2,0) {};
      \node (D) at (3,0) {};
      \node (E) at (0.5,-1) {};
      \node (F) at (1.5,-1) {};
      \node (G) at (2.5,-1) {};
      \node (H) at (1,-2) {};
      \node (I) at (2,-2) {};
      \node (J) at (1.5,-3) {};
      \draw [blue,thin, ->] (A) to (E);
      \draw [blue,thin, ->] (B) to (F);
      \draw [blue,thin, ->] (C) to (G);
      \draw [blue,thin, ->] (E) to (H);
      \draw [blue,thin, <-] (F) to (I);
      \draw [blue,thin, ->] (H) to (J);
    \end{tikzpicture}
    \vspace{-1mm}
    \end{minipage} &
    $\left\{\begin{gathered}
     \ell\leq m\leq r \leq 0 \\
     s\leq 0 \\
     t\leq 0 \\
     s\leq n
     \end{gathered}\right\}$
    \\
    \hline
    $D_fM^\mfrak{g}_\mfrak{p}(\lambda)$ &  $T(v)$ &
    \begin{minipage}{0.31\textwidth}
    \vspace{2mm}
    \begin{tikzpicture}[x=7mm, y=7mm, >=latex, yscale=1, xscale=2, every node/.style={circle, draw=white, fill=black, line width=1.5mm, inner sep=1mm}]
      \node (A) at (0,0) {};
      \node (B) at (1,0) {};
      \node (C) at (2,0) {};
      \node (D) at (3,0) {};
      \node (E) at (0.5,-1) {};
      \node (F) at (1.5,-1) {};
      \node (G) at (2.5,-1) {};
      \node (H) at (1,-2) {};
      \node (I) at (2,-2) {};
      \node (J) at (1.5,-3) {};
      \draw [blue,thin, ->] (A) to (E);
      \draw [blue,thin, ->] (B) to (F);
      \draw [blue,thin, ->] (C) to (G);
      \draw [blue,thin, ->] (E) to (H);
      \draw [blue,thin, <-] (F) to (I);
    \end{tikzpicture}
    \vspace{-1mm}
    \end{minipage} &
    $\left\{\begin{gathered}
      m\leq r \leq 0 \\
      s\leq 0 \\
      t\leq 0 \\
      s\leq n \\
    \end{gathered}\right\}$
    \\
    \hline
    $T_fM^\mfrak{g}_\mfrak{p}(\lambda)$ &  $T(v+\delta^{1,1})$  &
    \begin{minipage}{0.31\textwidth}
    \vspace{2mm}
    \begin{tikzpicture}[x=7mm, y=7mm, >=latex, yscale=1, xscale=2, every node/.style={circle, draw=white, fill=black, line width=1.5mm, inner sep=1mm}]
      \node (A) at (0,0) {};
      \node (B) at (1,0) {};
      \node (C) at (2,0) {};
      \node (D) at (3,0) {};
      \node (E) at (0.5,-1) {};
      \node (F) at (1.5,-1) {};
      \node (G) at (2.5,-1) {};
      \node (H) at (1,-2) {};
      \node (I) at (2,-2) {};
      \node (J) at (1.5,-3) {};
      \draw [blue,thin, ->] (A) to (E);
      \draw [blue,thin, ->] (B) to (F);
      \draw [blue,thin, ->] (C) to (G);
      \draw [blue,thin, ->] (E) to (H);
      \draw [blue,thin, <-] (F) to (I);
      \draw [blue,thin, <-] (H) to (J);
    \end{tikzpicture}
    \vspace{-1mm}
    \end{minipage} &
    $\left\{\begin{gathered}
      m\leq r \leq 0 \\
      s\leq 0 \\
      t\leq 0 \\
      s\leq n \\
      m \leq \ell
    \end{gathered}\right\}$
    \\
    \hline
    $D^\nu_fM^\mfrak{g}_\mfrak{p}(\lambda)$ &  $T(v+\nu\delta^{1,1})$ &
    \begin{minipage}{0.31\textwidth}
    \vspace{2mm}
    \begin{tikzpicture}[x=7mm, y=7mm, >=latex, yscale=1, xscale=2, every node/.style={circle, draw=white, fill=black, line width=1.5mm, inner sep=1mm}]
      \node (A) at (0,0) {};
      \node (B) at (1,0) {};
      \node (C) at (2,0) {};
      \node (D) at (3,0) {};
      \node (E) at (0.5,-1) {};
      \node (F) at (1.5,-1) {};
      \node (G) at (2.5,-1) {};
      \node (H) at (1,-2) {};
      \node (I) at (2,-2) {};
      \node (J) at (1.5,-3) {};
      \draw [blue,thin, ->] (A) to (E);
      \draw [blue,thin, ->] (B) to (F);
      \draw [blue,thin, ->] (C) to (G);
      \draw [blue,thin, ->] (E) to (H);
      \draw [blue,thin, <-] (F) to (I);
    \end{tikzpicture}
    \vspace{-1mm}
    \end{minipage} &
    $\left\{\begin{gathered}
      m\leq r \leq 0 \\
      s\leq 0 \\
      t\leq 0 \\
      s\leq n \\
    \end{gathered}\right\}$
    \\
    \hline
\end{tabular}
\caption{Principal nilpotent orbit}
\label{tab:prinicpal orbit}
\vspace{-1mm}
\end{table}

Let us note that the parabolic subalgebra $\mfrak{p}$ in Table \ref{tab:prinicpal orbit} corresponds to the Borel subalgebra $\mfrak{b}$. In addition, we assume that $\nu\in \C\setminus \Z$ satisfies $\nu+ \langle \lambda+\rho, \alpha_{1,2}^\vee \rangle \notin \Z$. The $\mfrak{g}$-modules $M^\mfrak{g}_\mfrak{p}(\lambda)$, $T_fM^\mfrak{g}_\mfrak{p}(\lambda)$ and $D_f^\nu M^\mfrak{g}_\mfrak{p}(\lambda)$ are simple, since $\mcal{C}$ is the maximal set of relations satisfied by the tableau $T(w)$. All weight multiplicities are infinite.

%%%%%%%%%%%%%%%%%%%%%%%%%%%%%%%%%%%%%%%%%%%%%%%%%%%%%%%%%%%%%%%%%%%%%%%%%%%%%%%%%%%%%%%%%%%

\subsection{Subregular nilpotent orbit}

For $\lambda \in \Lambda_k(\mfrak{p}_{\alpha_3}^{\rm min}) \subset \widebar{{\rm Pr}}_k^{\smash{\mcal{O}_{\rm subreg}}}$, we have
\begin{align*}
 \lambda=\big(\lambda_1- {\textstyle {p \over q}}a,\lambda_2-{\textstyle {p \over q}}b,\lambda_3\big)\ \text{with}\ a,b \in \N,\ a+b \leq q-1.
\end{align*}
We set $v_1-v_3=\langle \lambda+\rho,\alpha_1^\vee \rangle \notin\Z$, $v_3-v_2=\langle \lambda+\rho,\alpha_2^\vee \rangle \notin \Z$,  $v_2-v_4=\langle \lambda+\rho,\alpha_3^\vee \rangle \in\N$ such that $v_1+v_2+v_3+v_4=-6$.
\begin{align*}
T(v) = \begin{aligned}
\begin{tikzpicture}[x=15mm, y=7mm, every node/.style={draw, minimum width=14.3mm, minimum height=6.3mm}]
    \node at (0,0) {$v_1$};
    \node at (1,0) {$v_2$};
    \node at (2,0) {$v_4$};
    \node at (3,0) {$v_3$};
    \node at (0.5,-1) {$v_1$};
    \node at (1.5,-1) {$v_2$};
    \node at (2.5,-1) {$v_3$};
    \node at (1,-2) {$v_1$};
    \node at (2,-2) {$v_2+1$};
    \node at (1.5,-3) {$v_1$};
\end{tikzpicture}
\end{aligned}
\end{align*}

\vspace{-5mm}

\begin{table}[h]
\centering
\renewcommand{\arraystretch}{2.5}
\begin{tabular}{|c|c|c|c|}
    \hline
    $V_\mcal{C}(T(w))$ & $T(w)$ & $\mcal{C}$ & $\mcal{B}_\mcal{C}(T(w))$ \\
    \hline
    $M^\mfrak{g}_\mfrak{p}(\lambda)$ & $T(v)$ &
    \begin{minipage}{0.31\textwidth}
    \vspace{2mm}
    \begin{tikzpicture}[x=7mm, y=7mm, >=latex, yscale=1, xscale=2, every node/.style={circle, draw=white, fill=black, line width=1.5mm, inner sep=1mm}]
      \node (A) at (0,0) {};
      \node (B) at (1,0) {};
      \node (C) at (2,0) {};
      \node (D) at (3,0) {};
      \node (E) at (0.5,-1) {};
      \node (F) at (1.5,-1) {};
      \node (G) at (2.5,-1) {};
      \node (H) at (1,-2) {};
      \node (I) at (2,-2) {};
      \node (J) at (1.5,-3) {};
      \draw [blue,thin, ->] (A) to (E);
      \draw [blue,thin, ->] (B) to (F);
      \draw [blue,thin, <-] (C) to (F);
      \draw [blue,thin, ->] (D) to (G);
      \draw [blue,thin, ->] (E) to (H);
      \draw [blue,thin, <-] (F) to (I);
      \draw [blue,thin, ->] (H) to (J);
    \end{tikzpicture}
    \vspace{-1mm}
    \end{minipage} &
    $\left\{\begin{gathered}
     \ell\leq m\leq r \leq 0 \\
     -\lambda_3 \leq s\leq 0 \\
      t\leq 0 \\
      s\leq n
    \end{gathered}\right\}$
    \\
    \hline
    $D_fM^\mfrak{g}_\mfrak{p}(\lambda)$ &  $T(v)$ &
    \begin{minipage}{0.31\textwidth}
    \vspace{2mm}
    \begin{tikzpicture}[x=7mm, y=7mm, >=latex, yscale=1, xscale=2, every node/.style={circle, draw=white, fill=black, line width=1.5mm, inner sep=1mm}]
      \node (A) at (0,0) {};
      \node (B) at (1,0) {};
      \node (C) at (2,0) {};
      \node (D) at (3,0) {};
      \node (E) at (0.5,-1) {};
      \node (F) at (1.5,-1) {};
      \node (G) at (2.5,-1) {};
      \node (H) at (1,-2) {};
      \node (I) at (2,-2) {};
      \node (J) at (1.5,-3) {};
      \draw [blue,thin, ->] (A) to (E);
      \draw [blue,thin, ->] (B) to (F);
      \draw [blue,thin, <-] (C) to (F);
      \draw [blue,thin, ->] (D) to (G);
      \draw [blue,thin, ->] (E) to (H);
      \draw [blue,thin, <-] (F) to (I);
    \end{tikzpicture}
    \vspace{-1mm}
    \end{minipage} &
    $\left\{\begin{gathered}
      m\leq r \leq 0 \\
      -\lambda_3 \leq  s\leq 0 \\
      t\leq 0 \\
      s\leq n \\
    \end{gathered}\right\}$
    \\
    \hline
    $T_fM^\mfrak{g}_\mfrak{p}(\lambda)$ &  $T(v+\delta^{1,1})$  &
    \begin{minipage}{0.31\textwidth}
    \vspace{2mm}
    \begin{tikzpicture}[x=7mm, y=7mm, >=latex, yscale=1, xscale=2, every node/.style={circle, draw=white, fill=black, line width=1.5mm, inner sep=1mm}]
      \node (A) at (0,0) {};
      \node (B) at (1,0) {};
      \node (C) at (2,0) {};
      \node (D) at (3,0) {};
      \node (E) at (0.5,-1) {};
      \node (F) at (1.5,-1) {};
      \node (G) at (2.5,-1) {};
      \node (H) at (1,-2) {};
      \node (I) at (2,-2) {};
      \node (J) at (1.5,-3) {};
      \draw [blue,thin, ->] (A) to (E);
      \draw [blue,thin, ->] (B) to (F);
      \draw [blue,thin, <-] (C) to (F);
      \draw [blue,thin, ->] (D) to (G);
      \draw [blue,thin, ->] (E) to (H);
      \draw [blue,thin, <-] (F) to (I);
      \draw [blue,thin, <-] (H) to (J);
    \end{tikzpicture}
    \vspace{-1mm}
    \end{minipage} &
    $\left\{\begin{gathered}
     m\leq r \leq 0 \\
     -\lambda_3 \leq s\leq 0 \\
     t\leq 0 \\
     s\leq n \\
     m \leq \ell
    \end{gathered}\right\}$
    \\
    \hline
    $D^\nu_fM^\mfrak{g}_\mfrak{p}(\lambda)$ &  $T(v+\nu\delta^{1,1})$ &
    \begin{minipage}{0.31\textwidth}
    \vspace{2mm}
    \begin{tikzpicture}[x=7mm, y=7mm, >=latex, yscale=1, xscale=2, every node/.style={circle, draw=white, fill=black, line width=1.5mm, inner sep=1mm}]
      \node (A) at (0,0) {};
      \node (B) at (1,0) {};
      \node (C) at (2,0) {};
      \node (D) at (3,0) {};
      \node (E) at (0.5,-1) {};
      \node (F) at (1.5,-1) {};
      \node (G) at (2.5,-1) {};
      \node (H) at (1,-2) {};
      \node (I) at (2,-2) {};
      \node (J) at (1.5,-3) {};
      \draw [blue,thin, ->] (A) to (E);
      \draw [blue,thin, ->] (B) to (F);
      \draw [blue,thin, <-] (C) to (F);
      \draw [blue,thin, ->] (D) to (G);
      \draw [blue,thin, ->] (E) to (H);
      \draw [blue,thin, <-] (F) to (I);
    \end{tikzpicture}
    \vspace{-1mm}
    \end{minipage} &
    $\left\{\begin{gathered}
      m\leq r \leq 0 \\
      -\lambda_3 \leq  s\leq 0 \\
      t\leq 0 \\
      s\leq n \\
    \end{gathered}\right\}$
    \\
    \hline
\end{tabular}
\caption{Subregular nilpotent orbit}
\label{tab:subregular orbit}
\vspace{-1mm}
\end{table}

Let us note that the parabolic subalgebra $\mfrak{p}$ in Table \ref{tab:subregular orbit} corresponds to the minimal parabolic subalgebra $\mfrak{p}_{\alpha_3}^{\rm min}$. In addition, we assume that $\nu\in \C\setminus \Z$ satisfies $\nu+ \langle \lambda+\rho, \alpha_{1,2}^\vee \rangle \notin \Z$. The $\mfrak{g}$-modules $M^\mfrak{g}_\mfrak{p}(\lambda)$, $T_fM^\mfrak{g}_\mfrak{p}(\lambda)$ and $D_f^\nu M^\mfrak{g}_\mfrak{p}(\lambda)$ are simple, since $\mcal{C}$ is the maximal set of relations satisfied by the tableau $T(w)$. All weight multiplicities are infinite.

%%%%%%%%%%%%%%%%%%%%%%%%%%%%%%%%%%%%%%%%%%%%%%%%%%%%%%%%%%%%%%%%%%%%%%%%%%%%%%%%%%%%%%%%%%%

\subsection{Rectangular nilpotent orbit}

For $\lambda \in \Lambda_k(\mfrak{p}_{\alpha_2}^{\rm max}) \subset \widebar{{\rm Pr}}_k^{\smash{\mcal{O}_{\rm rect}}}$, we have
\begin{align*}
 \lambda=\big(\lambda_1,\lambda_2-{\textstyle {p \over q}}a,\lambda_3\big)\ \text{with}\ a \in \N,\ a \leq q-1.
\end{align*}
We set $v_1-v_3=\langle \lambda+\rho,\alpha_1^\vee \rangle \in\N$, $v_3-v_2=\langle \lambda+\rho,\alpha_2^\vee \rangle \notin \Z$,  $v_2-v_4=\langle \lambda+\rho,\alpha_3^\vee \rangle \in\N$ such that $v_1+v_2+v_3+v_4=-6$.
\begin{align*}
T(v) = \begin{aligned}
\begin{tikzpicture}[x=15mm, y=7mm, every node/.style={draw, minimum width=14.3mm, minimum height=6.3mm}]
    \node at (0,0) {$v_1$};
    \node at (1,0) {$v_3$};
    \node at (2,0) {$v_2$};
    \node at (3,0) {$v_4$};
    \node at (0.5,-1) {$v_1$};
    \node at (1.5,-1) {$v_3$};
    \node at (2.5,-1) {$v_2$};
    \node at (1,-2) {$v_1$};
    \node at (2,-2) {$v_2+1$};
    \node at (1.5,-3) {$v_1$};
\end{tikzpicture}
\end{aligned}
\end{align*}

\vspace{-5mm}

\begin{table}[h]
\centering
\renewcommand{\arraystretch}{2.5}
\begin{tabular}{|c|c|c|c|}
    \hline
    $V_\mcal{C}(T(w))$ & $T(w)$ & $\mcal{C}$ & $\mcal{B}_\mcal{C}(T(w))$ \\
    \hline
    $M^\mfrak{g}_\mfrak{p}(\lambda)$ & $T(v)$ &
    \begin{minipage}{0.31\textwidth}
    \vspace{2mm}
    \begin{tikzpicture}[x=7mm, y=7mm, >=latex, yscale=1, xscale=2, every node/.style={circle, draw=white, fill=black, line width=1.5mm, inner sep=1mm}]
      \node (A) at (0,0) {};
      \node (B) at (1,0) {};
      \node (C) at (2,0) {};
      \node (D) at (3,0) {};
      \node (E) at (0.5,-1) {};
      \node (F) at (1.5,-1) {};
      \node (G) at (2.5,-1) {};
      \node (H) at (1,-2) {};
      \node (I) at (2,-2) {};
      \node (J) at (1.5,-3) {};
      \draw [blue,thin, ->] (A) to (E);
      \draw [blue,thin, <-] (B) to (E);
      \draw [blue,thin, ->] (B) to (F);
      \draw [blue,thin, ->] (C) to (G);
      \draw [blue,thin, <-] (D) to (G);
      \draw [blue,thin, ->] (E) to (H);
      \draw [blue,thin, <-] (F) to (H);
      \draw [blue,thin, <-] (G) to (I);
      \draw [blue,thin, ->] (H) to (J);
    \end{tikzpicture}
    \vspace{-1mm}
    \end{minipage} &
    $\left\{\begin{gathered}
     t-\lambda_1\leq m\leq r \\
     -\lambda_3 \leq s\leq 0 \\
      t\leq 0 \\
      s\leq n \\
     -\lambda_1 \leq r \leq 0 \\
      \ell \leq m
    \end{gathered}\right\}$
    \\
    \hline
    $D_fM^\mfrak{g}_\mfrak{p}(\lambda)$ &  $T(v)$ &
    \begin{minipage}{0.31\textwidth}
    \vspace{2mm}
    \begin{tikzpicture}[x=7mm, y=7mm, >=latex, yscale=1, xscale=2, every node/.style={circle, draw=white, fill=black, line width=1.5mm, inner sep=1mm}]
      \node (A) at (0,0) {};
      \node (B) at (1,0) {};
      \node (C) at (2,0) {};
      \node (D) at (3,0) {};
      \node (E) at (0.5,-1) {};
      \node (F) at (1.5,-1) {};
      \node (G) at (2.5,-1) {};
      \node (H) at (1,-2) {};
      \node (I) at (2,-2) {};
      \node (J) at (1.5,-3) {};
      \draw [blue,thin, ->] (A) to (E);
      \draw [blue,thin, <-] (B) to (E);
      \draw [blue,thin, ->] (B) to (F);
      \draw [blue,thin, ->] (C) to (G);
      \draw [blue,thin, <-] (D) to (G);
      \draw [blue,thin, ->] (E) to (H);
      \draw [blue,thin, <-] (F) to (H);
      \draw [blue,thin, <-] (G) to (I);
    \end{tikzpicture}
    \vspace{-1mm}
    \end{minipage} &
    $\left\{\begin{gathered}
     t-\lambda_1\leq m\leq r \\
     -\lambda_3 \leq s\leq 0 \\
     t\leq 0 \\
     s\leq n \\
     -\lambda_1 \leq r \leq 0 \\
    \end{gathered}\right\}$
    \\
    \hline
    $T_fM^\mfrak{g}_\mfrak{p}(\lambda)$ &  $T(v+\delta^{1,1})$  &
    \begin{minipage}{0.31\textwidth}
    \vspace{2mm}
    \begin{tikzpicture}[x=7mm, y=7mm, >=latex, yscale=1, xscale=2, every node/.style={circle, draw=white, fill=black, line width=1.5mm, inner sep=1mm}]
      \node (A) at (0,0) {};
      \node (B) at (1,0) {};
      \node (C) at (2,0) {};
      \node (D) at (3,0) {};
      \node (E) at (0.5,-1) {};
      \node (F) at (1.5,-1) {};
      \node (G) at (2.5,-1) {};
      \node (H) at (1,-2) {};
      \node (I) at (2,-2) {};
      \node (J) at (1.5,-3) {};
      \draw [blue,thin, ->] (A) to (E);
      \draw [blue,thin, <-] (B) to (E);
      \draw [blue,thin, ->] (B) to (F);
      \draw [blue,thin, ->] (C) to (G);
      \draw [blue,thin, <-] (D) to (G);
      \draw [blue,thin, ->] (E) to (H);
      \draw [blue,thin, <-] (F) to (H);
      \draw [blue,thin, <-] (G) to (I);
      \draw [blue,thin, <-] (H) to (J);
    \end{tikzpicture}
    \vspace{-1mm}
    \end{minipage} &
    $\left\{\begin{gathered}
     t-\lambda_1\leq m\leq r \\
     -\lambda_3 \leq s\leq 0 \\
      t\leq 0 \\
     s\leq n \\
     -\lambda_1 \leq r \leq 0 \\
      m \leq \ell
    \end{gathered}\right\}$
    \\
    \hline
    $D^\nu_fM^\mfrak{g}_\mfrak{p}(\lambda)$ &  $T(v+\nu\delta^{1,1})$ &
    \begin{minipage}{0.31\textwidth}
    \vspace{2mm}
    \begin{tikzpicture}[x=7mm, y=7mm, >=latex, yscale=1, xscale=2, every node/.style={circle, draw=white, fill=black, line width=1.5mm, inner sep=1mm}]
      \node (A) at (0,0) {};
      \node (B) at (1,0) {};
      \node (C) at (2,0) {};
      \node (D) at (3,0) {};
      \node (E) at (0.5,-1) {};
      \node (F) at (1.5,-1) {};
      \node (G) at (2.5,-1) {};
      \node (H) at (1,-2) {};
      \node (I) at (2,-2) {};
      \node (J) at (1.5,-3) {};
      \draw [blue,thin, ->] (A) to (E);
      \draw [blue,thin, <-] (B) to (E);
      \draw [blue,thin, ->] (B) to (F);
      \draw [blue,thin, ->] (C) to (G);
      \draw [blue,thin, <-] (D) to (G);
      \draw [blue,thin, ->] (E) to (H);
      \draw [blue,thin, <-] (F) to (H);
      \draw [blue,thin, <-] (G) to (I);
    \end{tikzpicture}
    \vspace{-1mm}
    \end{minipage} &
    $\left\{\begin{gathered}
     t-\lambda_1\leq m\leq r \\
     -\lambda_3 \leq s\leq 0 \\
     t\leq 0 \\
     s\leq n \\
     -\lambda_1 \leq r \leq 0 \\
    \end{gathered}\right\}$
    \\
    \hline
\end{tabular}
\caption{Rectangular nilpotent orbit}
\label{tab:rectangular orbit}
\vspace{-1mm}
\end{table}

Let us note that the parabolic subalgebra $\mfrak{p}$ in Table \ref{tab:rectangular orbit} corresponds to the maximal parabolic subalgebra $\mfrak{p}_{\alpha_2}^{\rm max}$. In addition, we assume that $\nu\in \C\setminus \Z$ satisfies $\nu+ \langle \lambda+\rho, \alpha_{1,2}^\vee \rangle \notin \Z$. The $\mfrak{g}$-modules $M^\mfrak{g}_\mfrak{p}(\lambda)$, $T_fM^\mfrak{g}_\mfrak{p}(\lambda)$ and $D_f^\nu M^\mfrak{g}_\mfrak{p}(\lambda)$ are simple with unbounded weight multiplicities, since $\mcal{C}$ is the maximal set of relations satisfied by the tableau $T(w)$.

%%%%%%%%%%%%%%%%%%%%%%%%%%%%%%%%%%%%%%%%%%%%%%%%%%%%%%%%%%%%%%%%%%%%%%%%%%%%%%%%%%%%%%%%%%%

\subsection{Minimal nilpotent orbit}

For $\lambda \in \Lambda_k(\mfrak{p}_{\alpha_1}^{\rm max}) \subset \widebar{{\rm Pr}}_k^{\smash{\mcal{O}_{\rm min}}}$, we have
\begin{align*}
 \lambda=\big(\lambda_1-{\textstyle {p \over q}}a,\lambda_2, \lambda_3\big)\ \text{with}\ a \in \N,\ a \leq q-1.
\end{align*}
We set $v_1-v_3=\langle \lambda+\rho,\alpha_1^\vee \rangle \notin\Z$, $v_3-v_2=\langle \lambda+\rho,\alpha_2^\vee \rangle \in \N$,  $v_2-v_4=\langle \lambda+\rho,\alpha_3^\vee \rangle \in\N$ such that $v_1+v_2+v_3+v_4=-6$.
\begin{align*}
T(v) = \begin{aligned}
\begin{tikzpicture}[x=15mm, y=7mm, every node/.style={draw, minimum width=14.3mm, minimum height=6.3mm}]
    \node at (0,0) {$v_1$};
    \node at (1,0) {$v_3$};
    \node at (2,0) {$v_2$};
    \node at (3,0) {$v_4$};
    \node at (0.5,-1) {$v_1$};
    \node at (1.5,-1) {$v_3$};
    \node at (2.5,-1) {$v_2$};
    \node at (1,-2) {$v_1$};
    \node at (2,-2) {$v_2+1$};
    \node at (1.5,-3) {$v_1$};
\end{tikzpicture}
\end{aligned}
\end{align*}

\vspace{-5mm}

\begin{table}[h]
\centering
\renewcommand{\arraystretch}{2.5}
\begin{tabular}{|c|c|c|c|}
    \hline
    $V_\mcal{C}(T(w))$ & $T(w)$ & $\mcal{C}$ & $\mcal{B}_\mcal{C}(T(w))$ \\
    \hline
    $M^\mfrak{g}_\mfrak{p}(\lambda)$ & $T(v)$ &
    \begin{minipage}{0.31\textwidth}
    \vspace{2mm}
    \begin{tikzpicture}[x=7mm, y=7mm, >=latex, yscale=1, xscale=2, every node/.style={circle, draw=white, fill=black, line width=1.5mm, inner sep=1mm}]
      \node (A) at (0,0) {};
      \node (B) at (1,0) {};
      \node (C) at (2,0) {};
      \node (D) at (3,0) {};
      \node (E) at (0.5,-1) {};
      \node (F) at (1.5,-1) {};
      \node (G) at (2.5,-1) {};
      \node (H) at (1,-2) {};
      \node (I) at (2,-2) {};
      \node (J) at (1.5,-3) {};
      \draw [blue,thin, ->] (A) to (E);
      \draw [blue,thin, ->] (B) to (F);
      \draw [blue,thin, <-] (C) to (F);
      \draw [blue,thin, ->] (C) to (G);
      \draw [blue,thin, <-] (D) to (G);
      \draw [blue,thin, ->] (E) to (H);
      \draw [blue,thin, ->] (F) to (I);
      \draw [blue,thin, <-] (G) to (I);
      \draw [blue,thin, ->] (H) to (J);
    \end{tikzpicture}
    \vspace{-1mm}
    \end{minipage} &
    $\left\{\begin{gathered}
     \ell \leq m\leq r \leq 0 \\
     -\lambda_2  \leq t\leq 0 \\
     -\lambda_3  \leq s\leq 0 \\
     s\leq n \leq  t-\lambda_2
    \end{gathered}\right\}$
    \\
    \hline
    $D_fM^\mfrak{g}_\mfrak{p}(\lambda)$ &  $T(v)$ &
    \begin{minipage}{0.31\textwidth}
    \vspace{2mm}
    \begin{tikzpicture}[x=7mm, y=7mm, >=latex, yscale=1, xscale=2, every node/.style={circle, draw=white, fill=black, line width=1.5mm, inner sep=1mm}]
      \node (A) at (0,0) {};
      \node (B) at (1,0) {};
      \node (C) at (2,0) {};
      \node (D) at (3,0) {};
      \node (E) at (0.5,-1) {};
      \node (F) at (1.5,-1) {};
      \node (G) at (2.5,-1) {};
      \node (H) at (1,-2) {};
      \node (I) at (2,-2) {};
      \node (J) at (1.5,-3) {};
      \draw [blue,thin, ->] (A) to (E);
      \draw [blue,thin, ->] (B) to (F);
      \draw [blue,thin, <-] (C) to (F);
      \draw [blue,thin, ->] (C) to (G);
      \draw [blue,thin, <-] (D) to (G);
      \draw [blue,thin, ->] (E) to (H);
      \draw [blue,thin, ->] (F) to (I);
      \draw [blue,thin, <-] (G) to (I);
    \end{tikzpicture}
    \vspace{-1mm}
    \end{minipage} &
    $\left\{\begin{gathered}
     m\leq r \leq 0 \\
     -\lambda_2  \leq t\leq 0 \\
     -\lambda_3  \leq s\leq 0 \\
     s\leq n \leq  t-\lambda_2 \\
    \end{gathered}\right\}$
    \\
    \hline
    $T_fM^\mfrak{g}_\mfrak{p}(\lambda)$ &  $T(v+\delta^{1,1})$  &
    \begin{minipage}{0.31\textwidth}
    \vspace{2mm}
    \begin{tikzpicture}[x=7mm, y=7mm, >=latex, yscale=1, xscale=2, every node/.style={circle, draw=white, fill=black, line width=1.5mm, inner sep=1mm}]
      \node (A) at (0,0) {};
      \node (B) at (1,0) {};
      \node (C) at (2,0) {};
      \node (D) at (3,0) {};
      \node (E) at (0.5,-1) {};
      \node (F) at (1.5,-1) {};
      \node (G) at (2.5,-1) {};
      \node (H) at (1,-2) {};
      \node (I) at (2,-2) {};
      \node (J) at (1.5,-3) {};
      \draw [blue,thin, ->] (A) to (E);
      \draw [blue,thin, ->] (B) to (F);
      \draw [blue,thin, <-] (C) to (F);
      \draw [blue,thin, ->] (C) to (G);
      \draw [blue,thin, <-] (D) to (G);
      \draw [blue,thin, ->] (E) to (H);
      \draw [blue,thin, ->] (F) to (I);
      \draw [blue,thin, <-] (G) to (I);
      \draw [blue,thin, <-] (H) to (J);
    \end{tikzpicture}
    \vspace{-1mm}
    \end{minipage} &
    $\left\{\begin{gathered}
      m\leq r \leq 0 \\
      -\lambda_2  \leq t\leq 0 \\
      -\lambda_3  \leq s\leq 0 \\
      s\leq n \leq  t-\lambda_2 \\
      m \leq \ell
    \end{gathered}\right\}$
    \\
    \hline
    $D^\nu_fM^\mfrak{g}_\mfrak{p}(\lambda)$ &  $T(v+\nu\delta^{1,1})$ &
    \begin{minipage}{0.31\textwidth}
    \vspace{2mm}
    \begin{tikzpicture}[x=7mm, y=7mm, >=latex, yscale=1, xscale=2, every node/.style={circle, draw=white, fill=black, line width=1.5mm, inner sep=1mm}]
      \node (A) at (0,0) {};
      \node (B) at (1,0) {};
      \node (C) at (2,0) {};
      \node (D) at (3,0) {};
      \node (E) at (0.5,-1) {};
      \node (F) at (1.5,-1) {};
      \node (G) at (2.5,-1) {};
      \node (H) at (1,-2) {};
      \node (I) at (2,-2) {};
      \node (J) at (1.5,-3) {};
      \draw [blue,thin, ->] (A) to (E);
      \draw [blue,thin, ->] (B) to (F);
      \draw [blue,thin, <-] (C) to (F);
      \draw [blue,thin, ->] (C) to (G);
      \draw [blue,thin, <-] (D) to (G);
      \draw [blue,thin, ->] (E) to (H);
      \draw [blue,thin, ->] (F) to (I);
      \draw [blue,thin, <-] (G) to (I);
    \end{tikzpicture}
    \vspace{-1mm}
    \end{minipage} &
    $\left\{\begin{gathered}
     m\leq r \leq 0 \\
     -\lambda_2  \leq t\leq 0 \\
     -\lambda_3  \leq s\leq 0 \\
     s\leq n \leq  t-\lambda_2 \\
    \end{gathered}\right\}$
    \\
    \hline
\end{tabular}
\caption{Minimal nilpotent orbit}
\label{tab:minimal orbit}
\vspace{-1mm}
\end{table}

Let us note that the parabolic subalgebra $\mfrak{p}$ in Table \ref{tab:minimal orbit} corresponds to the maximal parabolic subalgebra $\mfrak{p}_{\alpha_1}^{\rm max}$. In addition, we assume that $\nu\in \C\setminus \Z$ satisfies $\nu+ \langle \lambda+\rho, \alpha_{1,2}^\vee \rangle \notin \Z$. The $\mfrak{g}$-modules $M^\mfrak{g}_\mfrak{p}(\lambda)$, $T_fM^\mfrak{g}_\mfrak{p}(\lambda)$ and $D_f^\nu M^\mfrak{g}_\mfrak{p}(\lambda)$ are simple with bounded weight multiplicities, since $\mcal{C}$ is the maximal set of relations satisfied by the tableau $T(w)$.

%%%%%%%%%%%%%%%%%%%%%%%%%%%%%%%%%%%%%%%%%%%%%%%%%%%%%%%%%%%%%%%%%%%%%%%%%%%%%%%%%%%%%%%%%%%
%%%%%%%%%%%%%%%%%%%%%%%%%%%%%%%%%%%%%%%%%%%%%%%%%%%%%%%%%%%%%%%%%%%%%%%%%%%%%%%%%%%%%%%%%%%

\section*{Acknowledgments}

V.\,F.\ is supported in part by the CNPq (304467/2017-0) and by the Fapesp (2018/23690-6). V.\,F.\ and L.\,K.\ are gratefully acknowledge the hospitality and excellent working conditions of the International Center for Mathematics of SUSTech (Shenzhen, China) where part of this work was done. O.\,A.\,H.\,M.\ is supported by the Coordena\c{c}\~ao de Aperfei\c{c}oamento de Pessoal de N\'ivel Superior -- Brasil (CAPES) -- Finance Code 001

%%%%%%%%%%%%%%%%%%%%%%%%%%%%%%%%%%%%%%%%%%%%%%%%%%%%%%%%%%%%%%%%%%%%%%%%%%%%%%%%%%%%%%%%%%
%%%%%%%%%%%%%%%%%%%%%%%%%%%%%%%%%%%%%%%%%%%%%%%%%%%%%%%%%%%%%%%%%%%%%%%%%%%%%%%%%%%%%%%%%%

%\bibliographystyle{amsalpha}
%\bibliographystyle{amsplain}
%\bibliography{C:/Temp/TeX/Clanek/Reference/reference}

\providecommand{\bysame}{\leavevmode\hbox to3em{\hrulefill}\thinspace}
\providecommand{\MR}{\relax\ifhmode\unskip\space\fi MR }
\providecommand{\MRhref}[2]{%
  \href{http://www.ams.org/mathscinet-getitem?mr=#1}{#2}
}
\providecommand{\href}[2]{#2}

\end{document}